\documentclass{article}

\usepackage{amsmath}
\usepackage{amssymb}
\usepackage[standard,amsmath,thmmarks,thref,hyperref]{ntheorem}
\usepackage{hyperref}
\usepackage{stmaryrd}
\usepackage{tikz}
	\usetikzlibrary{positioning}
\usepackage[all]{xy}
\usepackage{url}

\theoremstyle{plain}
\theoremheaderfont{\bfseries}
\theorembodyfont{\normalfont}
\theoremseparator{.}
\newtheorem{thm}{Theorem}[section]
\newtheorem{prop}[thm]{Proposition}
\newtheorem{rmk}[thm]{Remark}

\newtheorem{lem}[thm]{Lemma}
\newtheorem{nota}[thm]{Notation}
\newtheorem{cor}[thm]{Corollary}
\newtheorem{defn}[thm]{Definition}
\newtheorem{ex}[thm]{Example}

\numberwithin{equation}{section}

\theoremstyle{nonumberplain}
\theoremheaderfont{\itshape}
\theorembodyfont{\normalfont}
\theoremseparator{.}
\theoremsymbol{\ensuremath{\Box}}
\newtheorem{pf}{Proof}

\newcommand{\calO}{\mathcal{O}}
\newcommand{\calE}{\mathcal{E}}
\newcommand{\Z}{\mathbb{Z}}
\newcommand{\Q}{\mathbb{Q}}
\newcommand{\N}{\mathbb{N}}
\newcommand{\F}{\mathbb{F}}
\newcommand{\tcalO}{\widetilde{\calO}}
\newcommand{\tphi}{\widetilde{\phi}}
\newcommand{\alg}{^\mathrm{alg}}
\newcommand{\Set}[2]{\left\{ #1 \,\middle|\, #2 \right\}}
\newcommand{\abs}[2][]{|#2|_{#1}}
\newcommand{\df}[1]{\emph{#1}}
\newcommand{\set}[2]{\{ #1 \,|\, #2 \}}
\newcommand{\Perf}{\mathbf{Perf}}
\newcommand{\pfg}{\oldpi_1^\mathrm{prof}}
\newcommand{\setcat}{\mathbf{Set}}
\newcommand{\calo}{\calO}
\newcommand{\Abs}[2][]{\left|#2\right|_{#1}}
\newcommand{\nr}{^\mathrm{nr}}

\newcommand{\sep}{\mathrm{sep}}
\newcommand{\perf}{^\mathrm{perf}}
\newcommand{\tcalE}{\widetilde{\calE}}
\newcommand{\R}{\mathbb{R}}

\newcommand{\GQpD}{G_{\mathbb{Q}_p,\Delta}}
\newcommand{\GKD}{G_{K,\Delta}}
\newcommand{\GFD}{G_{F,\Delta}}
\newcommand{\GQpa}{G_{\mathbb{Q}_p,\alpha}}
\newcommand{\HQpD}{H_{\mathbb{Q}_p,\Delta}}
\newcommand{\HKD}{H_{K,\Delta}}
\newcommand{\HFD}{H_{F,\Delta}}
\newcommand{\HQpa}{H_{\mathbb{Q}_p,\alpha}}
\newcommand{\HKa}{H_{K,\alpha}}
\newcommand{\Zp}{\mathbb{Z}_p}
\newcommand{\Fp}{\mathbb{F}_p}

\newcommand{\Qp}{\mathbb{Q}_p}

\newcommand{\OED}{\mathcal{O}_{\mathcal{E}_\Delta}}
\newcommand{\TD}{\mathbf{1}_\Delta}
\newcommand{\bs}{\llbracket}
\newcommand{\js}{\rrbracket}

\renewcommand{\phi}{\varphi}
\let\oldpi\pi
\renewcommand{\pi}{\varpi}

\DeclareMathOperator{\Hom}{Hom}
\DeclareMathOperator{\GL}{GL}
\DeclareMathOperator{\Gal}{Gal}
\DeclareMathOperator{\Spa}{Spa}
\DeclareMathOperator{\Spd}{Spd}
\DeclareMathOperator{\FEt}{\mathbf{FEt}}
\DeclareMathOperator{\id}{id}
\DeclareMathOperator{\Spec}{Spec}
\DeclareMathOperator{\Frac}{Frac}
\DeclareMathOperator{\Ind}{Ind}
\DeclareMathOperator{\Ker}{Ker}
\DeclareMathOperator{\Tr}{Tr}
\DeclareMathOperator{\res}{res}

\title{Drinfeld's lemma for perfectoid spaces and overconvergence of multivariate $(\phi, \Gamma)$-modules}
\author{Annie Carter \and Kiran S.\ Kedlaya \and Gergely Z\'abr\'adi}
\date{\today}

\begin{document}

\maketitle

\begin{abstract}
Let $p$ be a prime, let $K$ be a finite extension of $\mathbb{Q}_p$, and let $n$ be a positive integer. We construct equivalences of categories between continuous $p$-adic representations of the $n$-fold product of the absolute Galois group $G_K$ and $(\phi, \Gamma)$-modules over one of several rings of $n$-variable power series.
The case $n=1$ recovers the original construction of Fontaine and the subsequent refinement by Cherbonnier--Colmez;
for general $n$, the case $K = \mathbb{Q}_p$ had been previously treated by the third author. To handle general $K$ uniformly,
we use a form of Drinfeld's lemma on the profinite fundamental groups of products of spaces in characteristic $p$, but for perfectoid spaces instead of schemes. We also construct the multivariate analogue of the Herr complex to compute Galois cohomology; the case $K = \mathbb{Q}_p$ had been previously treated by Pal and the third author, and we reduce to this case using a form of Shapiro's lemma.
\end{abstract}
\maketitle

\section{Introduction}

\subsection{Overview}

Throughout this paper, fix a prime number $p$.
The theory of \df{$(\varphi, \Gamma)$-modules} was introduced by Fontaine \cite{Fontaine}
as a tool for describing and classifying continuous representations of the Galois group of a finite 
extension of $\Q_p$ on a finite-dimensional $\Q_p$-vector space. Thanks to subsequent refinements,
notably the work of Cherbonnier--Colmez \cite{CherbonnierColmez} and Berger \cite{Berger1,Berger2}, it has become clear that essentially all of $p$-adic Hodge theory can be formulated in terms of $(\varphi, \Gamma)$-modules; moreover, this formulation has driven much recent progress in the subject and powered some notable applications in arithmetic geometry. See \cite{KedlayaNewMethods} for a quick introduction to this circle of ideas or
\cite{Schneider} for a more in-depth treatment.

The goal of this paper is to initiate a systematic development of \df{multivariate} $(\varphi, \Gamma)$-modules, founded upon the theory of perfectoid spaces,
as a tool for studying representations of \df{products} of Galois groups of $p$-adic fields.
The relevance of such representations may not be immediately clear from general considerations of arithmetic geometry; however, products of Galois groups occur naturally in the approach to geometric Langlands 
developed for $\GL_2$ by Drinfeld \cite{DrinfeldICM} and extended to $\GL_n$ by L. Lafforgue
\cite{LLafforgue} and to other reductive groups by V. Lafforgue \cite{VLafforgue1, VLafforgue2}.

The relationship between multivariate $(\varphi, \Gamma)$-modules and Galois representations
was previously explored by the third author \cite{PalZabradi, Zabradi} from a slightly different point of view:
this line of inquiry emerged as part of a program to extend Colmez's construction of the $p$-adic local Langlands correspondence for $\GL_2(\Q_p)$ \cite{Colmez} by exhibiting analogues of $(\varphi, \Gamma)$-modules obtained from higher-rank groups  \cite{SchneiderVignerasZabradi,ZabradiFunctor}. 

One motivation for consolidating the theory of multivariate $(\varphi, \Gamma)$-modules is to prepare for an eventual unification of this program with the work of V. Lafforgue described above; however, such a unification lies far beyond the scope of the present work. Another motivation is to flesh out the point of view suggested in the last paragraph of the introduction of \cite{Zabradi}, which relates multivariate $(\varphi, \Gamma)$-modules to 
a form of Drinfeld's lemma for perfectoid spaces (more on which below).

\subsection{Main results}

Before diving into the weeds of perfectoid spaces, we give a brief indication of our main results
(and recall that the case $K = \Q_p$ was treated in \cite{PalZabradi, Zabradi}).
Let $K$ denote a finite extension of $\Q_p$ with absolute Galois group $G_K$. Let $G_{K,\Delta}$ be the Cartesian power of $G_K$ indexed by the finite set $\Delta$. 
(One can also consider products $G_{K_1} \times \cdots \times G_{K_n}$ where $K_1,\dots,K_n$ are possibly distinct finite extensions of $\Q_p$; to keep notation under control, we suppress this level of generality until the end of the paper.)

\begin{thm}[see \thref{full-equivalence}] \label{mainthm1}
The category of continuous representations of $G_{K,\Delta}$ on finite free $\Z_p$-modules is canonically equivalent to the category of projective \'etale $(\phi_\Delta, \Gamma_\Delta)$-modules over 
each of the rings 
\[
\calO_{\calE_{\Delta}(K)}, \quad \tcalO_{\calE_{\Delta}(K)}, \quad \calO^\dagger_{\calE_{\Delta}(K)}, \quad \tcalO^\dagger_{\calE_{\Delta}(K)},
\]
described below. A similar statement also holds for representations of $G_{K,\Delta}$ on finite-dimensional $\Q_p$-vector spaces; see \thref{full-equivalence-Qp}.
\end{thm}

\begin{thm}[see \thref{shapirophigamma}] \label{mainthm2}
In \thref{mainthm1}, the Galois cohomology of a representation of $G_{K,\Delta}$ is canonically isomorphic to the
continuous $(\varphi_\Delta, \Gamma_\Delta)$-cohomology of the corresponding $(\varphi_\Delta, \Gamma_\Delta)$-module (i.e., the cohomology of the complex of continuous $(\varphi_\Delta,\Gamma_\Delta)$-cochains valued in the module).
\end{thm}

To unpack this, let us start with the case where $\Delta$ is a singleton set.
In this case, $\calO_{\calE_\Delta(K)}$ is none other than Fontaine's ring $\calO_{\calE}$ which serves as his original base ring for $(\varphi, \Gamma)$-modules \cite{Fontaine}; it is the $p$-adic completion of a Laurent series ring in a variable $\varpi$ with coefficients in a certain finite \'etale extension of $\Z_p$,
carrying actions of a Frobenius lift $\varphi$ and a profinite group $\Gamma_K$ (the Galois group of the maximal cyclotomic extension of $K$), which commute with each other. The other rings are built from $\calO_{\calE_\Delta(K)}$ in such a way as to also carry actions of $\varphi$ and $\Gamma_K$: the tilde indicates an enlargement that makes the action of $\varphi$ bijective (\df{perfection}), while the dagger indicates passage to a subring satisfying a growth condition (\df{overconvergence}). A \df{projective \'etale $(\phi, \Gamma)$-module} over one of these rings is then a finite projective (hence free) module $M$ equipped with commuting semilinear actions of $\varphi$ and $\Gamma_K$, for which the linearization $\varphi^* M \to M$ of the $\varphi$-action is an isomorphism. The equivalence of categories stated in \thref{mainthm1} then incorporates Fontaine's original equivalence, together with its refinement by Cherbonnier--Colmez \cite{CherbonnierColmez} using the overconvergent subring.
The description of Galois cohomology stated in \thref{mainthm2} is due to Herr \cite{Herr1,Herr2}.

For general $\Delta$, the rings in question are certain topological Cartesian powers of the rings arising in the singleton case. In particular, for each $\alpha \in \Delta$, there will be an element $\varpi_\alpha$ arising from the factor of the product indexed by $\alpha$; moreover, there will be a partial Frobenius lift $\varphi_\alpha$ and a profinite group $\Gamma_{K,\alpha}$ which act on $\varpi_\alpha$ but not on the other variables.
(The use of the symbols $\Delta$ and $\alpha$ is meant to suggest roots of a Lie algebra; the relevance of this
will not be apparent herein, but can be seen more directly in earlier work on the subject, especially
\cite{ZabradiFunctor}.)

\subsection{Perfectoid spaces and Drinfeld's lemma}

We next explain what perfectoid spaces and Drinfeld's lemma have to do with each, and with the aforementioned results.

The theory of \df{perfectoid spaces}, while having notched diverse achievements since its promulgation in the early 2010s, is at its heart a geometric reinterpretation of the core ideas of $p$-adic Hodge theory
(as in \cite{BhattMorrowScholze,KedlayaLiu, KedlayaLiu2,ScholzeHodge}). It begins with a vast generalization of the ``field of norms'' isomorphism between the absolute Galois groups of $\Q_p(\mu_{p^\infty})$ and $\F_p\llparenthesis t \rrparenthesis$ introduced by
Fontaine--Wintenberger \cite{FontaineWintenberger2,FontaineWintenberger1} and underpinning the classical theory of $(\varphi, \Gamma)$-modules. In this generalization, every field $K$ of characteristic 0 which is complete for a nonarchimedean absolute value and ``sufficiently large'' (i.e., \df{perfectoid}) has associated with it a corresponding field $K^\flat$ of characteristic $p$ (the \df{tilt} of $K$) with a canonically isomorphic Galois group. Following the model of the \df{almost purity theorem} of Faltings, one then spreads out this correspondence to obtain a similar correspondence of spaces that matches up \'etale topoi; however, this takes place in the realm of analytic rather than algebraic geometry, specifically in Huber's category of \df{adic spaces}.

The term \df{Drinfeld's lemma} refers collectively to a statement used by Drinfeld
\cite[Theorem~2.1]{DrinfeldICM}, \cite[Proposition~6.1]{DrinfeldCohomology} in his study of the Langlands correspondence for $\GL_2$ over global function fields of characteristic $p$, together with subsequent generalizations
\cite[IV.2, Th\'eor\`eme~4]{LLafforgue},
\cite[Lemma~8.1.2]{Lau}, \cite[Theorem~17.2.4]{Scholze}, \cite[Theorem~4.2.12]{KedlayaSSS}.
These results address the behavior of (profinite) \'etale fundamental groups of schemes under formation of products, and in particular the significant discrepancy between this behavior in characteristic $0$ and in characteristic $p$.
In characteristic 0, \'etale fundamental groups are the profinite completions of topological fundamental groups,
and so their formation commutes with taking fiber products over an algebraically closed field. By contrast, it is easy to construct examples in characteristic $p$ where this commutativity fails; however, one obtains a similar statement by taking fiber products over $\F_p$, but with all of the objects divided by Frobenius (in the natural stack-theoretic sense).

It was first observed by Scholze \cite[Lecture~17]{Scholze} that Drinfeld's lemma might be related to the geometric simple connectivity of Fargues--Fontaine curves.
These appear in \cite{FarguesFontaine} as geometric objects whose vector bundles are closely related to $p$-adic Galois representations and $(\varphi, \Gamma)$-modules. The geometric simple connectivity property of the ``basic'' Fargues--Fontaine curve (the one associated to a completed algebraic closure of $\Q_p$) was established independently by Fargues--Fontaine \cite[Chapter~8]{FarguesFontaine} and  Weinstein \cite{WeinsteinGQp}; this has subsequently been extended to the curves associated to arbitrary algebraically closed perfectoid fields by the second author \cite{KedlayaSimple}. This result may be interpreted as the analogue of Drinfeld's lemma for the product of two geometric points; by emulating some of the steps
in the case of schemes, one obtains a full analogue of Drinfeld's lemma for perfectoid spaces \cite[Theorem~4.3.14]{KedlayaSSS}.
(It is natural to state the latter result in Scholze's language of \df{diamonds} \cite{Scholze}; we do so here, but diamonds are not essential for our present work.)

Using Drinfeld's lemma for perfectoid spaces, it is almost but not entirely straightforward to recover the multivariate analogue of Fontaine's original construction of $(\varphi, \Gamma)$-modules. The one difficulty is that taking a product of perfectoid fields in the sense of Drinfeld's lemma does not quite yield the adic space associated to the base ring of multivariate $(\varphi, \Gamma)$-modules, but rather a large open subspace thereof. To bridge the gap, we need an argument about bounded functions on certain non-quasicompact perfectoid spaces, which can be viewed as an application of the \df{perfectoid Riemann extension theorem (Hebbarkeitssatz)} appearing in the work of Scholze \cite{Scholze_TorsionCohomology} on torsion Galois representations associated to automorphic forms, and in the proofs of the direct summand conjecture by Andr\'e \cite{Andre2,Andre1} and Bhatt \cite{Bhatt}. With this in place, we can then exhibit the analogue of the Cherbonnier-Colmez refinement; for this, we prefer the simplified approach of \cite{KedlayaNewMethods} which avoids any use of Tate--Sen formalism.

As in \cite{PalZabradi,Zabradi} for the case $K = \Q_p$, one can extend various results about univariate $(\varphi, \Gamma)$-modules to the multivariate case, such as Herr's description of Galois cohomology \cite{Herr1,Herr2}. While these could be derived from scratch, we instead follow the approach of Liu \cite{LiuHerr} of reduction
to the case $K = \Q_p$ using a form of Shapiro's lemma for $(\varphi, \Gamma)$-modules.

\subsection{Followup questions}

At the end of the paper, we discuss a number of followup questions that emerge naturally from this line of investigation.
One of these is to extend the correspondence between Galois representations and $(\varphi, \Gamma)$-modules to products of \'etale fundamental groups (in the sense of de Jong \cite{deJong}) of rigid analytic spaces and perfectoid spaces, in the style of \cite{KedlayaLiu}.
Another is to modify the construction to reproduce some other examples of multivariate $(\varphi, \Gamma)$-modules in the literature, such as those exhibited by Berger \cite{BergerMulti1,BergerMulti2} using Lubin--Tate towers;
in particular, as suggested in \cite{KedlayaFrobMod}, it may be possible to establish an analogue of Cherbonnier--Colmez via this approach.

\subsection*{Acknowledgments}
Carter was a postdoctoral fellow of the UCSD Research Training Group in Algebra, Algebraic Geometry, and Number Theory (NSF grant DMS-1502651).
Kedlaya was supported by NSF (grants DMS-1501214, DMS-1802161) and by UC San Diego (Warschaw\-ski Professorship).
Z\'abr\'adi was supported by the J\'anos Bolyai Scholarship of the Hungarian Academy of Sciences, by an NKFIH Research grant FK-127906, by the MTA Alfr\'ed R\'enyi Institute of Mathematics Lend\"ulet Automorphic Research Group, and by the Thematic Excellence Programme, Industry and Digitization Subprogramme, NRDI Office, 2019.

\section{Notation}
\label{sec:Notation}

Fix a finite extension $K  / \Q_p$, an algebraic closure $K\alg$ of $K$, and a finite set $\Delta$ with $n$ elements. We are interested in continuous $\Z_p$-representations of the group
	\[ G_{K, \Delta} := \prod_{\alpha \in \Delta} \Gal(K\alg  / K). \]
Along the way, we will also encounter the groups
	\begin{align*}
	H_{K, \Delta} &:= \prod_{\alpha \in \Delta} \Gal(K\alg  / K(\mu_{p^\infty})) \\
	\Gamma_{K, \Delta} &:= \prod_{\alpha \in \Delta} \Gal(K(\mu_{p^\infty})  / K).
	\end{align*}
As mentioned in the introduction, one can also handle products of Galois groups of possibly distinct finite extensions of $\Q_p$, but to simplify notation we postpone discussion of this case until Section~\ref{distinct-factors}.

\subsection{The univariate case}

We begin by defining the suite of rings used to describe continuous $\Z_p$-represen\-tations of 
$G_K := \Gal(K\alg  / K)$. This will amount to treating the case where $\Delta$ is a singleton set; 
we will then go back and define various product constructions to handle general $\Delta$.
See \cite{BergerIntro} for a more detailed treatment.

\begin{nota}
Let $\calO_K$ be the valuation ring of $K$ and $k$ its residue field. Let $K_0 := \Frac W(k)$, identified with a subfield of $K$. 
Notice that the completion $\widehat{K_0(\mu_{p^\infty})}$ is a perfectoid field; let $E_0$ be its \df{tilt}
in the sense of the general theory of perfectoid rings (see Subsection \ref{subsec:adic-perfectoid-diamond}). For the moment, we recall that $E_0$ is set-theoretically the inverse limit of $\widehat{K_0(\mu_{p^\infty})}$ under the $p$-power map,
so we may choose $\epsilon := (1, \zeta_p, \zeta_{p^2}, \ldots) \in E_0$ where $\zeta_{p^n}$ denotes a primitive $p^n$-th root of unity; we then have $E_0 \simeq \widehat{k \llparenthesis \bar{\pi} \rrparenthesis\perf}$
for $\bar{\pi} := \epsilon - 1$.

Let $\pi := [\epsilon] - 1 \in W(E_0)$. Let $\calO_{\calE_0}$ be the $p$-adic completion of $\calO_{K_0} \llparenthesis \pi \rrparenthesis \subseteq W(E_0)$, i.e.
	\[ \calO_{\calE_0} = \Set{\sum_{n=-\infty}^\infty a_n \pi^n}{a_n \in \calO_{K_0},\ \text{$a_n \to 0$ as $n \to -\infty$}}. \] 
The ring $\calO_{\calE_0}$ is a complete discrete valuation ring, with maximal ideal generated by $p$ and residue field $k \llparenthesis \bar{\pi} \rrparenthesis$. Let $\calE_0 := \Frac \calO_{\calE_0} = \calO_{\calE_0}[p^{-1}]$. Denote by $\phi$ the unique ring homomorphism on $\calO_{K_0}$ lifting the absolute Frobenius on $k$. The rings $\calE_0$ and $\calO_{\calE_0}$ have commuting, $\calO_{K_0}$-semilinear actions of the map $\phi$ and the group $\Gamma_{K_0}$, defined by
	\[ \phi(\pi) = (1 + \pi)^p - 1 \qquad \text{and} \qquad \gamma(\pi) = (1 + \pi)^{\gamma} - 1, \]
where we identify an element $\gamma \in \Gamma_{K_0}$ with an element of $\Z_p^\times$ via the cyclotomic character.
(The action of $\phi$ on coefficients is via the Witt vector Frobenius; the action of $\Gamma_{K_0}$ on coefficients is trivial.)

Let $\tcalO_{\calE_0} := W(E_0)$; we may represent elements of $\tcalO_{\calE_0}$ as series
	\[ \sum_{n \in \Z[1/p]} a_n [\bar{\pi}]^n \]
with coefficients $a_n \in \calO_{K_0}$ such that $a_n \to 0$ as $n \to -\infty$ and, for each $c > 0$ and $r > 0$, there are at most finitely many coefficients $a_n$ with $\abs[p]{a_n} \geq c$ and $n \leq r$, since the ring of all such series satisfies the universal property for Witt vectors: it is $p$-adically complete and separated, and its residue field consists of elements \[ \sum_{n \in \Z[1/p]} \bar{a}_n \bar{\pi}^n \] such that for each $r > 0$ there are at most finitely many nonzero $\bar{a}_n$ with $n < r$; such elements constitute the ring $E_0$. The point is that we have a natural inclusion $\calO_{\calE_0} \subseteq \tcalO_{\calE_0}$ (in fact, $\tcalO_{\calE_0}$ is the completion of $\calO_{\calE_0}$ with respect to the \emph{weak topology}, defined below), but note that $\pi \neq [\bar{\pi}]$.
Let $\tcalE_0 := \Frac \tcalO_{\calE_0} = \tcalO_{\calE_0}[p^{-1}]$.
(Note the typographical convention whereby we write $\tcalO_{\calE}$ instead of the more logical $\calO_{\tcalE}$.
\end{nota}

So far, everything we have constructed depends only on $K_0$. We now introduce corresponding constructions depending on $K$, for which the previous constructions amount to the special case $K_0 = K$.
\begin{nota}
By the properties of the tilting construction, we have canonical isomorphisms 
\[
H_{K_0} \simeq \Gal(E_0^\sep  / E_0)\simeq \Gal(\calE_0\nr  / \calE_0) \simeq \Gal(\tcalE_0\nr  / \tcalE_0),
\]
where $\calE_0\nr$ and $\tcalE_0\nr$ denote the maximal unramified extensions of $\calE_0$ and $\tcalE_0$, respectively.
Define
\begin{gather*}
E := (E_0^\sep)^{H_K}, \quad
\calO_{\calE} := (\calO_{\calE_0^{\nr}})^{H_K}, \quad
\calE := (\calE_0\nr)^{H_K} = \calO_{\calE}[p^{-1}], \\
\tcalO_{\calE} := (\tcalO_{\calE_0^{\nr}})^{H_K}, \quad
\tcalE := (\tcalE_0\nr)^{H_K} = \tcalO_{\calE}[p^{-1}],
\end{gather*}
so that 
\[
H_K \simeq \Gal(E^\sep  / E) \simeq \Gal(\calE\nr  / \calE). 
\]
The rings $\calE$, $\calO_\calE$, $\tcalE$, and $\tcalO_\calE$ are stable under the actions of $\phi$ and $\Gamma_K$.

The ring $\calO_\calE$ is again a complete discrete valuation ring with maximal ideal generated by $p$. Its residue field is a finite separable extension of $k \llparenthesis \bar{\pi} \rrparenthesis$, which can itself be (noncanonically) identified with $k' \llparenthesis \bar{\pi}_K \rrparenthesis$ where $k'$ is the residue field of $K(\mu_{p^\infty})$. The completed perfect closure of the residue field of $\calO_\calE$ is canonically isomorphic to $E$, the residue field of $\tcalO_{\calE}$.
\label{tilted-extensions}
\end{nota}

\begin{defn}
A \df{series parameter} in $\calO_\calE$ is an element $\varpi_K$ which maps to $\bar{\pi}_K$ under some isomorphism $\calO_{\calE}/p\calO_\calE \simeq k' \llparenthesis \bar{\pi}_K \rrparenthesis$. 
(By the Cohen structure theorem, this just means that $\varpi_K$ maps to a uniformizer of the complete discretely valued field $\calO_{\calE}/p\calO_\calE$.)
For any such element, we may write 	\[ \calO_{\calE} = \Set{\sum_{n=-\infty}^\infty a_n \pi_K^n}{a_n \in W(k'),\ \text{$a_n \to 0$ as $n \to -\infty$}}. \]
An \df{overconvergent series parameter} in $\calO_\calE$ is a series parameter $\varpi_K$ satisfying the following additional condition: there exists a positive integer $c$ such that for each positive integer $n$, the element $[\epsilon - 1]^{cn} \varpi_K \in W(E)$ is congruent modulo $p^n$ to some element of $W(\calO_E)$.
The existence of such a series parameter follows from the discussion in \cite[Definition~2.1.4]{KedlayaNewMethods}.
\end{defn}

\subsection{Product constructions}

We now adapt the preceding constructions to products over the finite set $\Delta$. 
In the notation, we suppress the field $K$ for visual clarity; when it is necessary to specify $K$ explicitly,
we will write $\calE_\Delta(K)$ in place of $\calE_\Delta$ in all of the notations.
It will also be convenient to label the elements of $\Delta$ as $\alpha_1,\dots,\alpha_n$, but nothing will depend in an essential way on this ordering.

\begin{nota}
For each $\alpha \in \Delta$,
let $\tcalO_{\calE_\alpha}$ be a copy of $\tcalO_{\calE}$,
let $\calO_{\calE_\alpha}$ denote the corresponding copy of $\calO_{\calE}$ inside $\tcalO_{\calE_\alpha}$,
and let $\varpi_\alpha$ be an overconvergent series parameter in $\calO_{\calE_\alpha}$. The choice of $\varpi_\alpha$ is needed in order to articulate the subsequent definitions, but again nothing will depend in an essential way on this choice. Further, let $\bar{\pi}_\alpha$ denote the image of $\pi_\alpha$ in $R_\alpha:=\tcalO_{\calE_\alpha}/p \tcalO_{\calE_\alpha}$.

Let	\[ R_\Delta = R_\Delta(K) := (R_{\alpha_1} \mathbin{\widehat{\otimes}_{\F_p}} \cdots \mathbin{\widehat{\otimes}_{\F_p}} R_{\alpha_n})[\bar{\pi}_{\alpha_1}^{-1}, \ldots, \bar{\pi}_{\alpha_n}^{-1}], \]
where we identify $\bar{\pi}_{\alpha_1}$ with $\bar{\pi}_{\alpha_1} \otimes 1 \otimes \cdots \otimes 1$, etc., and the hats in the tensor product denote completion with respect to the $(\bar{\pi}_{\alpha_1}, \ldots, \bar{\pi}_{\alpha_n})$-adic topology.

 Let $\tcalO_{\calE_\Delta} := W(R_\Delta)$. Then 
\[
\tcalO_{\calE_\Delta} = \varprojlim_{m}
(\tcalO_{\calE_{\alpha_1}}/p^m \tcalO_{\calE_{\alpha_1}})  \mathbin{\widehat{\otimes}}_{\Z_p} \cdots \mathbin{\widehat{\otimes}}_{\Z_p} (\tcalO_{\calE_{\alpha_n}}/p^m \tcalO_{\calE_{\alpha_n}})
\]
where similarly the hats denote completion with respect to the $(\pi_{\alpha_1}, \ldots, \pi_{\alpha_n})$-adic topology.
Define 
\[
\calO_{\calE_\Delta} := \varprojlim_{m}
(\calO_{\calE_{\alpha_1}}/p^m \calO_{\calE_{\alpha_1}})  \mathbin{\widehat{\otimes}}_{\Z_p} \cdots \mathbin{\widehat{\otimes}}_{\Z_p} (\calO_{\calE_{\alpha_n}}/p^m \calO_{\calE_{\alpha_n}}),
\]
viewed as a subring of $\tcalO_{\calE_\Delta}$. 

For any of the above rings, let $\phi_\alpha$ and $\Gamma_{K, \alpha}$ denote the actions of $\phi$ and $\Gamma_K$ on the factor indexed by $\alpha$ in the product, fixing the other factors. Denote by $\phi_\Delta$ the monoid generated by the $\phi_\alpha$ for $\alpha \in \Delta$. 
\label{notation-for-product-rings-1}
\end{nota}

\begin{rmk}
The ring $\calO_{\calE_\Delta}$ is noetherian, but the ring $\tcalO_{\calE_\Delta}$ is not (because $R_\Delta$ is not).
\end{rmk}

We define a family of ``Gauss norms'' on $\calO_{\calE_{\Delta}}$ as follows.
\begin{nota}
Let $e$ be the ramification index of $K(\mu_{p^\infty})$ over $K_0(\mu_{p^\infty})$ (or equivalently, of $E$ over $E_0$).
For $j=1, \ldots, n$ and real number $r>0$, define the submultiplicative norm $|\ |_{j,r}$ on $\calO_{\calE_{\Delta}}$ by
	\[ \Abs[j,r]{\sum_{i_1,\ldots,i_n} (a_{i_1}\pi_{\alpha_1}^{i_1}) \otimes \cdots \otimes (a_{i_n} \pi_{\alpha_n}^{i_n})} = \sup_{i_1, \ldots, i_n} \{ p^{-r i_j p / (e(p-1))} \abs[p]{a_{i_1} \cdots a_{i_n}} \}. \]
For $r>0$, define the submultiplicative norm $|\ |_r$ on $\calO_{\calE_\Delta}$ by
	\[ \Abs[r]{\sum_{i_1, \ldots, i_n} (a_{i_1}\pi_{\alpha_1}^{i_1}) \otimes \cdots \otimes (a_{i_n} \pi_{\alpha_n}^{i_n})} = \sup_{i_1, \ldots, i_n} \{ p^{-rp \min\{i_j\} / (e(p-1))} \abs[p]{a_{i_1} \cdots a_{i_n}} \}. \]
Let $\calO^{j,r-}_{\calE_{\Delta}}$ denote those $Y \in \calO_{\calE_{\Delta}}$ with $|Y|_{j,r}$ finite, and let $\calO^{r-}_{\calE_{\Delta}}$ denote those $Y \in \calO_{\calE_{\Delta}}$ with $|Y|_r$ finite. Notice that $|Y|_r = \max_j |Y|_{j,r}$; in particular, $Y \in \calO^{r-}_{\calE_{\Delta}}$ if and only if $|Y|_{j,r}$ is finite for each $j$.
\label{notation-for-Gauss-norms}
\end{nota}

For $\tcalO_{\calE_\Delta}$, the following related construction is more useful.
\begin{nota}
For $j=1,\dots,n$, define the submultiplicative norm $|\ |'_j$ on $R_\Delta$ by
	\[ \Abs[j]{\sum_{i_1,\dots,i_n} (\bar{a}_{i_1} \bar{\pi}_{\alpha_1}^{i_1}) \otimes \cdots \otimes (\bar{a}_{i_n}\bar{\pi}_{\alpha_n}^{i_n})}' = (p^{-p/(e(p-1))})^{\min\set{i_j}{\bar{a}_{i_1} \cdots \bar{a}_{i_n} \neq 0}} \]
	and define $|\ |' := \max_j \{|\ |'_j\}$.
	For $r>0$, define submultiplicative norms $|\ |_{j,r}$ and $|\ |_r$ on $\tcalO_{\calE_\Delta}$ as follows: for
	$x = \sum_{m=0}^\infty p^m [\bar{x}_m] \in \tcalO_{\calE_\Delta}$, set
\[
	\Abs[j,r]{x} = \sup_m \{ p^{-m} \Abs[j]{\bar{x}_m}^{\prime r}\}, \quad
	\Abs[r]{x} = \sup_m \{ p^{-m} \Abs{\bar{x}_m}^{\prime r}\}.
\]
\label{notation-for-norm-R-Delta}
\end{nota}

In the sequel we equip $R_\Delta$ with the topology induced by $|\ |'$ which we call the \df{perfectoid topology} (see \thref{notation-for-perfectoid-product-rings}). Note that $R_\Delta$ is complete with respect to $|\ |'$,
but the induced topology is coarser than the $(\pi_{\alpha_1}, \ldots, \pi_{\alpha_n})$-adic topology.

Similarly, the rings $\tcalO_{\calE_\Delta}$ and $\calO_{\calE_\Delta}$ are complete with respect to each of the following topologies:
\begin{itemize}
\item
the $p$-adic topology;
\item
the \df{weak topology}:
the inverse limit topology induced by the $(\pi_{\alpha_1}, \ldots, \pi_{\alpha_n})$-adic topology modulo each power of $p$;
\item
the \df{perfectoid topology}: the inverse limit topology induced by the perfectoid topology modulo each power of $p$.
\end{itemize}

\begin{prop}
For $Y \in \bigcup_{r > 0} \calO^{r-}_{\calE_{\Delta}}$, let $\bar{Y}$ denote the reduction of $Y$ modulo $p$. If $\bar{Y} \neq 0$, we have \[\lim_{r \to 0^+} |Y|_r = 1.\] If $\bar{Y} = 0$, then \[\limsup_{r \to 0^+} |Y|_r \leq p^{-1}.\] These statements also hold with $|\ |_r$ replaced with $|\ |_{j,r}$ and/or with $\calO$ replaced by $\tcalO$.
\label{abs-val-lim}
\end{prop}
\begin{pf}
See \cite[Remark 1.7.3]{KedlayaNewMethods}.
\end{pf}

\begin{prop}
For $j=1,\dots,n$, for $r$ sufficiently small (depending on the choice of $\varpi_{\alpha_j}$),
the definitions of $|\, |_{j,r}$ on $\calO_{\calE_\Delta}$ and $\tcalO_{\calE_\Delta}$ agree.
Consequently, for $r$ sufficiently small (depending on the choices of all of the $\varpi_\alpha$)
the definitions of $|\, |_{r}$ on $\calO_{\calE_\Delta}$ and $\tcalO_{\calE_\Delta}$ agree.
\label{match-ring-norms}
\end{prop}
\begin{pf}
See \cite[Remark 2.2.8]{KedlayaNewMethods}.
\end{pf}

\begin{nota}
We may now define the following rings:
	\begin{align*}
	\calO^{j, \dagger}_{\calE_{\Delta}} &:= \bigcup_{r > 0} \calO^{j, r-}_{\calE_{\Delta}} & \calO^\dagger_{\calE_{\Delta}} &:= \bigcup_{r > 0} \calO^{r-}_{\calE_{\Delta}} \\
	\tcalO^{j, \dagger}_{\calE_{\Delta}} &:= \bigcup_{r > 0} \tcalO^{j, r-}_{\calE_{\Delta}} & \tcalO^\dagger_{\calE_{\Delta}} &:= \bigcup_{r > 0} \tcalO^{r-}_{\calE_{\Delta}}.
	\end{align*}
We also define $\calE^\dagger_\Delta := \calO^\dagger_{\calE_\Delta}[p^{-1}]$ and $\tcalE^\dagger_\Delta := \tcalO^\dagger_{\calE_\Delta}[p^{-1}]$. (Note again the typographical choice to write $\calO^\dagger_{\calE}$ instead of the more logical $\calO_{\calE^\dagger}$, and so on.)
The rings $\calO^\dagger_{\calE_{\Delta}}$, $\tcalO^\dagger_{\calE_{\Delta}}$, $\calE^\dagger_\Delta$, and $\tcalE^\dagger_\Delta$ are preserved by the actions of $\phi_\Delta$ and $\Gamma_{K, \Delta}$. By \thref{match-ring-norms}, within $\tcalE_\Delta$ we have the equalities
\[
\calO_{\calE_\Delta} \cap \tcalO^\dagger_{\calE_\Delta} = \calO^\dagger_{\calE_\Delta},
\quad
\calE_\Delta \cap \tcalE^\dagger_{\Delta} = \calE^\dagger_\Delta.
\]
\label{notation-for-product-rings-2}
\end{nota}

\subsection{$(\varphi, \Gamma)$-modules}

We now give the definition of $(\varphi, \Gamma)$-modules over the various rings we have constructed; this is formally similar to the univariate case.

\begin{defn}
Let $\calO$ be a ring with commuting actions of $\phi_\Delta$ and $\Gamma_{K, \Delta}$ (such as $\calO_{\calE_{\Delta}}$ or $\tcalO_{\calE_{\Delta}}$). A \df{$\phi_\Delta$-module} over $\calO$ is a finitely presented $\calO$-module with commuting, semilinear actions of the $\phi_\alpha$. A \df{$(\phi_\Delta, \Gamma_{K, \Delta})$-module} over $\calO$ is a finitely presented $\calO$-module $M$ with commuting semilinear actions of the $\phi_\alpha$ and the $\Gamma_{K,\alpha}$. We apply additional ring-theoretic modifiers (such as ``torsion'' or ``projective'') by passing them through to the underlying $\calO$-module.
\end{defn}

\begin{defn}
Let $\calO$ be one of the rings $\calO_{\calE_\Delta}$, $\tcalO_{\calE_\Delta}$, $\calO^\dagger_{\calE_\Delta}$, or $\tcalO^\dagger_{\calE_\Delta}$.
A $\varphi_\Delta$-module $M$ over $\calO$ is \df{\'etale} if the induced maps
	\[ \phi_\alpha^\ast M \to M, \qquad a \otimes x \mapsto a \phi_\alpha(x) \]
are isomorphisms for all $\alpha \in \Delta$; here $\phi_\alpha^\ast M$ denotes the module $\calO \otimes_{\phi, \calO} M$, in which $a \otimes bx = a \phi_\alpha(b) \otimes x$ for $a, b \in \calO$, $x \in M$. In the case when $M$ is a free module, this condition holds if and only if for each $\alpha$, $\phi_\alpha$ maps some basis of $M$ to another basis of $M$; it then maps every basis of $M$ to another basis.
A $(\varphi_\Delta, \Gamma_{K, \Delta})$-module over $\calO$ is \df{\'etale} if its underlying $\varphi_\Delta$-module is \'etale.
\end{defn}

\begin{defn}
Let $\calO$ be a ring of the form $\calO_0[p^{-1}]$ where $\calO_0$ is one of the rings $\calO_{\calE_\Delta}$, $\tcalO_{\calE_\Delta}$, $\calO^\dagger_{\calE_\Delta}$, or $\tcalO^\dagger_{\calE_\Delta}$; that is, $\calO$ is one of the rings $\calE_\Delta$, $\tcalE_\Delta$, $\calE^\dagger_{\Delta}$, or $\tcalE^\dagger_{\Delta}$.
A $\varphi_\Delta$-module or $(\varphi_\Delta, \Gamma_{K, \Delta})$-module 
$M$ over $\calO$ is \df{\'etale} if it has the form $M_0[p^{-1}]$
for some projective \'etale $\varphi_\Delta$-module  or $(\varphi_\Delta, \Gamma_{K, \Delta})$-module
$M_0$ over $\calO_0$;
in particular, in this setup we are insisting that $M$ be projective.
\label{phiGamma-rational}
\end{defn}

\begin{rmk}
We point out two subtleties in the previous definitions which do not have much impact on our work here, but may become relevant when comparing with other literature.
On the one hand, in \thref{phiGamma-rational}, our definition of \'etale 
$(\varphi_\Delta, \Gamma_{K, \Delta})$-modules imposes the \'etale condition not just on the action of $\varphi_\Delta$, but also on $\Gamma_{K, \Delta}$. In some settings, it might be preferable to have a definition in which the \'etale condition is described solely in terms of $\varphi_\Delta$.

On the other hand, to avoid imposing the \'etale condition on $\Gamma_{K, \Delta}$, one probably has to replace it with a condition asserting that the action map $\Gamma_{K, \Delta} \times M \to M$ is continuous for some topology on $M$. That topology should be induced by some topology on the base ring $\calO$ for which the action of $\Gamma_{K, \Delta}$ is itself continuous (e.g., the $p$-adic topology or the weak topology); since $M$ appears in both the source and target of the action map, the continuity conditions for different topologies on $\calO$ are not immediately comparable even if the topologies themselves are comparable.

In our setup, for \'etale $(\varphi_\Delta, \Gamma_{K, \Delta})$-modules, any continuity condition of this form will follow \df{a posteriori} from the comparison between these objects and representations of $G_{K, \Delta}$, and so no such condition needs to be included in either definitions or theorem statements. For an example of the tradeoff when we try to weaken the \'etale condition, see
\thref{weaker-phiGamma-rational}.
\label{remarks-on-rational-phiGamma-def}
\end{rmk}

\begin{rmk}
Suppose $M$ and $N$ are \'etale $\phi_\Delta$-modules (or \'etale $(\phi_\Delta, \Gamma_{K, \Delta})$-modules) over $\calO$. Then provided that $\Hom_{\calO}(M,N)$ is finitely presented, we may view it as an \'etale $\phi_\Delta$-module (or \'etale $(\phi_\Delta, \Gamma_{K, \Delta})$-module) by requiring, as appropriate,
	\begin{align*}
	\phi_\alpha(f)(\phi_\alpha(\mathbf{e})) &= \phi_\alpha(f(\mathbf{e})) \\
	\gamma_\alpha(f)(\gamma_\alpha(\mathbf{e})) &= \gamma_\alpha(f(\mathbf{e})) 
	\end{align*}
for $\alpha \in \Delta$. The morphisms $M \to N$ of $\phi_\Delta$-modules (or of $(\phi_\Delta, \Gamma_{K, \Delta})$-modules) are exactly those $\calO$-module homomorphisms fixed by $\phi_\Delta$ (or by $\phi_\Delta$ and $\Gamma_{K, \Delta}$).

We point out two key cases in which the finite presentation condition on $\Hom_{\calO}(M,N)$ is always satisfied: the case where $M$ is projective, and the case where $\calO$ is noetherian. As noted above, the latter holds for $\calO = \calO_{\calE_\Delta}$; it also holds for $\calO^\dagger_{\calE_\Delta}$, but we will not need this fact.
\label{hom-phi-module}
\end{rmk}

\begin{rmk}
In connection with the previous remark, we note that by
\cite[Proposition~3.2.13]{KedlayaLiu}, a $\varphi_\Delta$-module over $\tcalO_{\calE_\Delta}$ or $\tcalO^\dagger_{\calE_\Delta}$
which is flat over $\Z/p^n \Z$ is a finite projective $\tcalO_{\calE_\Delta}/p^n \tcalO_{\calE_\Delta}$-module.
\label{phi-module-flat}
\end{rmk}

\section{Drinfeld's lemma for diamonds}

In this section, we give a very brief summary of Drinfeld's lemma for perfectoid spaces and diamonds,
following \cite[Lecture~4]{KedlayaSSS}.

\subsection{Adic spaces, perfectoid spaces, and diamonds}
\label{subsec:adic-perfectoid-diamond}

We begin by recalling some terminology and results about adic spaces, perfectoid spaces, and diamonds.
 Good introductions to this material can be found in \cite{Huber}, \cite{Weinstein}, and \cite{KedlayaSSS}; see also \cite{Scholze}.

\begin{defn}
We say that $(A, A^+)$ is a \df{Huber pair} if
	\begin{enumerate}
	\item $A$ is a \df{Huber ring}, i.e.\ it is a topological ring which contains an open subring $A_0$ whose topology is the $I$-adic topology for some finitely generated ideal $I$ of $A_0$ (the ring $A_0$ is called a \df{ring of definition} and the ideal $I$ an \df{ideal of definition});
	\item $A^+$ is a \df{ring of integral elements}, i.e.\ $A^+$ is an open, integrally closed subring contained in the subring $A^\circ$ of power-bounded elements of $A$.
	\end{enumerate}
\end{defn}

\begin{defn}
A Huber ring $A$ (or a Huber pair $(A, A^+)$) is \df{uniform} if the subring $A^\circ$ is bounded. It is \df{analytic} if the topologically nilpotent elements of $A$ generate the unit ideal. It is \df{Tate} if $A$ contains a \df{pseudo-uniformizer}, i.e.\ a topologically nilpotent unit.
\end{defn}

\begin{rmk} \label{remark-plus-ring-in-uniform-case}
Beware that the condition that a Huber pair $(A,A^+)$ be uniform is only a condition on $A$. While it does imply that $(A,A^\circ)$ is also a Huber pair, it does not force the inclusion $A^+ \subseteq A^\circ$ to be an equality. On the other hand, the difference between the two is not large: the topologically nilpotent elements of $A^\circ$ form an ideal of $A^\circ$, and this ideal is itself contained in $A^+$ (because $A^+$ is open and integrally closed).
\end{rmk}

\begin{defn}
The \df{adic spectrum} of a Huber pair $(A, A^+)$ is the set $\Spa(A, A^+)$ of equivalence classes of continuous valuations $v$ on $A$ such that $v(f) \leq 1$ for all $f \in A^+$.
\end{defn}

\begin{defn}[{{\cite[Definition 1.2.1]{KedlayaSSS}}}]
Given a Huber pair $(A, A^+)$, a \df{rational subspace} of $X = \Spa(A, A^+)$ is a set of the form
	\[ X \left(\frac{f_1, \ldots, f_n}{g}\right) := \set{v \in X}{v(f_i) \leq v(g) \neq 0\ \forall\ i} \]
for some collection of elements $f_1, \ldots, f_n, g$ which generate an open ideal in $A$. Rational subspaces provide a basis for a topology on $X$. The \df{rational localization} corresponding to a rational subspace $Y$ of $X$ is a morphism $(A, A^+) \to (B, B^+)$ of complete Huber pairs which is initial among morphisms for which $\Spa(B, B^+)$ maps into $Y$. This morphism induces a homeomorphism $\Spa(B, B^+) \simeq  Y$ which moreover identifies rational subspaces of $\Spa(B, B^+)$ with rational subspaces of $X$ contained in $Y$.

In the case that  $\{f_1,\dots,f_n,g\}$ generates the unit ideal, $B$ is canonically isomorphic to the quotient of $A \langle T_1,\dots,T_n \rangle$
by the closure of the ideal generated by $\{g T_1 - f_1, \dots, gT_n - f_n\}$. This is always the case when $A$ is analytic. (Otherwise, one must also invert $g$.)
\end{defn}

\begin{defn}
Let $(A, A^+)$ be a Huber pair. We define the \df{structure presheaf} of $X := \Spa(A, A^+)$ as follows: If $U \subseteq X$ is an open subspace, let
	\begin{align*}
	\calO(U) &:= \varprojlim B \\
	\calO^+(U) &:= \varprojlim B^+
	\end{align*}
where the limits run over rational localizations $(A, A^+) \to (B, B^+)$ such that $\Spa(B, B^+) \subseteq U$. In particular, if $U = \Spa(B, B^+)$, then $(\calO(U), \calO^+(U)) = (B, B^+)$. The structure presheaf has the property that
	\[ \calO^+(U) = \set{f \in \calO(U)}{v(f) \leq 1\ \forall\ v \in U}. \]
A Huber pair $(A, A^+)$ is \df{sheafy} if $\calO$ is a sheaf; in this case $\calO^+$ is also.

In general, the sheafiness condition in a Huber pair is highly nontrivial, which causes severe complications in setting up the theory of adic spaces. Fortunately, for perfectoid spaces this complication disappears; see \thref{remark:perfectoid}.
\end{defn}

\begin{defn}
An \df{adic space} is a topological space $X$ together with a sheaf of topological rings $\calO_X$ and a continuous valuation on the stalk $\calO_{X, x}$ for each $x \in X$, such that $X$ can be covered by open subsets of the form $\Spa(A, A^+)$ where each $(A, A^+)$ is a sheafy Huber pair. In particular, for any identification of an open subset of $X$ with a space $\Spa(A, A^+)$, and a corresponding identification of $\calO_X$ with the structure sheaf on $\Spa(A,A^+)$, the valuation on $\calO_{X, x}$ is the one whose valuation ring is the stalk of $\calO_X^+$ at $x$.
\end{defn}

\begin{defn}
A \df{perfectoid ring} is a uniform, analytic Huber ring $A$ which contains a ring $A^+$ of integral elements (which is then a ring of definition) and an ideal of definition $I \subseteq A^+$ such that $p \in I^p$ and the $p$-th power map $A^+ / I \to A^+ / I^p$ is surjective. A \df{perfectoid pair} is a Huber pair $(A, A^+)$ with $A$ perfectoid.
\end{defn}

\begin{rmk} \label{remark:perfect implies uniform}
An important special case is when $A$ is an analytic Huber ring of characteristic $p$ which is perfect (that is,
its Frobenius endomorphism is bijective). In this case, $A$ is automatically uniform (see \cite[Example~2.1.2]{KedlayaSSS}) and hence a perfectoid ring.
\end{rmk}

\begin{rmk}[{{\cite[Theorem 3.1.3]{Weinstein}}, {\cite[Theorem 2.5.3]{KedlayaSSS}}}] \label{remark:perfectoid}
Let $(A, A^+)$ be a perfectoid pair. Then
	\begin{enumerate}
	\item $(A, A^+)$ is sheafy;
	\item $X := \Spa(A, A^+)$ is an adic space;
	\item $\calO_X(U)$ is a perfectoid ring for all rational subsets $U \subseteq X$.
	\end{enumerate}
\end{rmk}

\begin{thm}[Tilting Correspondence, {{\cite[Theorem 2.3.9]{KedlayaSSS}}}]
There is an equivalence of categories
	\begin{gather*}
	(A, A^+) \mapsto (A^\flat, A^{\flat +}, I) = \left(\varprojlim_{x^p \mapsfrom x} A, \varprojlim_{x^p \mapsfrom x} A^+, \ker\left(W(A^{\flat +}) \to A^+\right)\right) \\
	\left(W^b(R) / I W^b(R), W(R^+) / I\right) \mapsfrom (R, R^+, I)
	\end{gather*}
between the category of perfectoid pairs $(A, A^+)$ and the category of characteristic $p$ perfectoid pairs $(R, R^+)$ together with a primitive ideal $I$ of $W(R^+)$. Here $W(R^+)$ denotes the ring of $p$-typical Witt vectors over $R^+$, and $W^b(R)$ denotes the ring of $p$-typical Witt vectors
	\[ \sum_{n=0}^\infty p^n [\bar{x}_n] \]
such that the set $\set{\bar{x}_n}{n \in \N}$ is bounded in $R$. An ideal of $W(R^+)$ is \df{primitive}
if it is principal on some generator $z = \sum_{n=0}^\infty p^n [\bar{z}_n]$ for which $\bar{z}_0$ is topologically nilpotent and $\bar{z}_1$ is a unit in $R^+$.
\end{thm}

\begin{thm}[{{\cite[Theorem 1.5.6]{KedlayaNewMethods}}}]
If $L$ is a perfectoid field, then the absolute Galois groups of $L$ and $L^\flat$ are isomorphic as topological groups.
\end{thm}

\begin{thm}[{{\cite[Theorem 2.5.1]{KedlayaSSS}}}]
Given a perfectoid pair $(A, A^+)$, there is a homeomorphism
	\begin{align*}
	\Spa(A, A^+) &\to \Spa(A^\flat, A^{\flat +}) \\
	v &\mapsto v^\flat,
	\end{align*}
where $v^\flat((f_n)_n) = v(f_0)$. If $f = (f_n)_n \in A^\flat$, we define $f^\sharp := f_0 \in A$; then we have $v^\flat(f) = v(f^\sharp)$.
\end{thm}

\begin{defn}
A \df{perfectoid space} is an adic space which is covered by open subspaces of the form $\Spa(A, A^+)$ with $A$ perfectoid. Any such subspace is called an \df{affinoid perfectoid space}.
\end{defn}

\begin{defn}
A morphism $(R, R^+) \to (S, S^+)$ of Huber pairs is \df{finite \'etale} if $S$ is a finite \'etale $R$-algebra with the induced topology and $S^+$ is the integral closure of $R^+$ in $S$.

A morphism $f\colon X \to Y$ of adic spaces is \df{finite \'etale} if there is a cover of $Y$ by open affinoids $V \subseteq Y$ such that the pre-image $U = f^{-1}(V)$ is affinoid and the associated morphism of Huber pairs
	\[ (\calO_Y(V), \calO^+_Y(V)) \to (\calO_X(U), \calO^+_X(U)) \]
is finite \'etale.

A morphism $f\colon X \to Y$ of adic spaces is \df{\'etale} if for all points $x \in X$, there exist open neighborhoods $U$ of $x$ and $V$ of $f(x)$ and a commutative diagram
	\begin{center}
	\begin{tikzpicture}
	\node (u) {$U$};
	\node[right=of u] (w) {$W$};
	\node[below=of w] (v) {$V$};
	\draw[->] (u) -- node[above]{$j$} (w);
	\draw[->] (w) -- node[right]{$p$} (v);
	\draw[->] (u) -- node[below left]{$f|_U$} (v);
	\end{tikzpicture}
	\end{center}
with $j$ an open embedding and $p$ finite \'etale.

A morphism $f\colon X \to Y$ of perfectoid spaces is \df{pro-\'etale} if locally on $X$ it is of the form $\Spa(A_\infty, A^+_\infty) \to \Spa(A, A^+)$, where $A$ and $A_\infty$ are perfectoid rings, and
	\[ (A_\infty, A^+_\infty) = \left[ \varinjlim (A_i, A^+_i) \right]^\wedge \]
is a filtered colimit of pairs $(A_i, A^+_i)$ such that $\Spa(A_i, A^+_i) \to \Spa(A, A^+)$ is \'etale.
\end{defn}

\begin{defn}
Let $\Perf$ denote the category of perfectoid spaces of characteristic $p$.
We hereafter view $\Perf$ as a site using the \df{pro-\'etale topology}, whose coverings (\df{pro-\'etale coverings}) are collections of morphisms $\set{f_i\colon X_i \to X}{i \in I}$ such that each $f_i$ is pro-\'etale and for all quasi-compact open $U \subseteq X$, there exists a collection of quasi-compact open sets $U_i \subseteq X_i$ indexed by a finite subset $I_U \subseteq I$ such that
	\[ U = \bigcup_{i \in I_U} f_i(U_i). \]
The pro-\'etale topology is subcanonical; that is, for $X \in \Perf$, the functor $h_X$ on $\Perf$ represented by $X$ is a sheaf (see \thref{perfectoid-spaces-as-diamonds}).
\end{defn}

\begin{defn}
A morphism of sheaves $\mathcal{F} \to \mathcal{G}$ on $\Perf$ is \df{pro-\'etale} if for all perfectoid spaces $X$ and maps $h_X \to \mathcal{G}$, the pullback $h_X \times_{\mathcal{G}} \mathcal{F}$ is representable by a perfectoid space $Y$, and the morphism $Y \to X$ corresponding to the map $h_Y = h_X \times_{\mathcal{G}} \mathcal{F} \to h_X$ is pro-\'etale.

Let $\mathcal{F}$ be a sheaf on $\Perf$. A \df{pro-\'etale equivalence relation} is a monomorphism $\mathcal{R} \to \mathcal{F} \times \mathcal{F}$ in the category of sheaves on $\Perf$ such that each projection $\mathcal{R} \to \mathcal{F}$ is pro-\'etale, and such that for all objects $S$ of $\Perf$, the image of the map $\mathcal{R}(S) \to \mathcal{F}(S) \times \mathcal{F}(S)$ is an equivalence relation on $\mathcal{F}(S)$.
\end{defn}

\begin{defn}
A \df{diamond} is a sheaf $\mathcal{F}$ on $\Perf$ which is the quotient of a perfectoid space by a pro-\'etale equivalence relation. More precisely, there exist a perfectoid space $X$ and a representable equivalence relation $\mathcal{R} \to h_X \times h_X$ such that the two projections $\mathcal{R} \to h_X$ are pro-\'etale.
(Compare \cite[Definition~11.1]{ScholzeDiamonds}.)
\end{defn}

\begin{defn}
	\begin{enumerate}
	\item Given a perfectoid space $X$, we denote by $X^\diamond$ the representable sheaf $h_{X^\flat}$.
	\item The \df{diamond spectrum} of a perfectoid ring $A$ is the sheaf $\Spd A = (\Spa(A, A^\circ))^\diamond$.
	\end{enumerate}
\end{defn}


\begin{thm}[{{\cite[Theorem 3.5.2]{Weinstein}, \cite[Theorem 3.8.2]{KedlayaSSS}}}] 
\label{perfectoid-spaces-as-diamonds}
If $X$ is a perfectoid space, then $X^\diamond$ is a diamond. Moreover, for any perfectoid space $S$, the functor $X \mapsto X^\diamond$ from perfectoid spaces over $S$ to diamonds over $S^\diamond$ is fully faithful.
\end{thm}

\subsection{Finite \'etale covers and profinite fundamental groups}

\begin{nota}
Let $X$ be a scheme, a perfectoid space, or a diamond. We denote by $\FEt(X)$ the category of finite \'etale coverings of $X$.
(For $X$ a diamond, a finite \'etale covering of $X$ is a representable morphism $Y \to X$ whose pullback to any perfectoid space is a finite \'etale covering.)

We use a similar notation for formal quotients by group actions: if we take some $X$ and formally quotient by a group $\Gamma$ of automorphisms of $X$, then $\FEt(X/\Gamma)$ is the category of objects of $\FEt(X)$ equipped with an action of $\Gamma$ for which the structure morphism is $\Gamma$-equivariant.
\end{nota}

\begin{nota}
Let $X_1, \ldots, X_n$ be diamonds, and let $X = X_1 \times \cdots \times X_n$. 
Let $\phi_i$ denote the absolute Frobenius of $X_i$ (induced by the $p$-th power map on rings), and let $\phi = \phi_1 \times \cdots \times \phi_n$ be the absolute Frobenius of $X$. 
As per \cite[Remark~4.2.14]{KedlayaSSS},
let $X/\Phi$ be the functor from perfectoid spaces to sets taking $Y$
to the set of tuples $(f, \beta_1,\dots,\beta_n)$ where $f: Y \to X$ is a morphism and
$\beta_i: Y \to \phi_i^* Y$ are morphisms which ``commute and compose to $\phi$" in the sense of 
 \cite[Definition~4.2.10]{KedlayaSSS}. That is, for any $i$ and $j$, the composition
 \[
 Y \stackrel{\beta_i}{\to} \phi_i^* Y \stackrel{\beta_j}{\to} \phi_i^* \phi_j^* Y = (\phi_i \circ \phi_j)^* Y = (\phi_j \circ \phi_i)^* Y
 \]
 is the same as the corresponding composition with $i$ and $j$ reversed; and the composition
 \[
 Y \stackrel{\beta_1}{\to} \cdots \stackrel{\beta_n}{\to} (\phi_1 \circ \cdots \circ \phi_n)^* Y = \phi^* Y
 \]
 is relative Frobenius for $Y/X$.
 
A more concrete, but less symmetric, description of $X/\Phi$ can be given by picking an index $j \in \{1,\dots,n\}$; one then has a canonical isomorphism
\[
X/\Phi \cong X / \langle \phi_1,\dots,\widehat{\phi_j},\dots,\phi_n \rangle
\]
given by discarding $\beta_j$ and instead recovering it from the other data.
In other words, $X/\Phi$ is the formal quotient of $X$ by the group $\Phi$ generated by 
$\langle \phi_1,\dots,\widehat{\phi_j},\dots,\phi_n \rangle$; we will refer to the group $\Phi$ on its own in various contexts where the choice of $j$ does not matter.

%
\end{nota}

\begin{defn}
Let $X$ and $\Phi$ be as above. Given a geometric point $\bar{x}$ of $X$, the \df{profinite fundamental group} $\pfg(X/\Phi, \bar{x})$ is the group of natural isomorphisms of the functor $\FEt(X/\Phi) \to \setcat$ taking a covering $Y$ to the underlying set of $Y \times_X \bar{x}$.
\end{defn}

\begin{thm}[{{\cite[Remark 4.1.4]{KedlayaSSS}}}]
Suppose that $F$ is a perfectoid field of characteristic $p$, and let $F\alg$ be an algebraic closure of $F$. Then
	\[ G_F = \Gal(F\alg  / F) \simeq \pfg(X, \bar{x}) \]
for $X := \Spd(F)$ and $\bar{x} := \Spd(F\alg)$.
\label{pfg-galois-group}
\end{thm}

\begin{thm}[Drinfeld's lemma, {{\cite[Theorem 4.3.14]{KedlayaSSS}}}]
Let $X_1, \ldots, X_n$ be connected spatial (in the sense of 
\cite[Definition~17.3.1]{Scholze}) diamonds and put $X :=  X_1 \times \cdots \times X_n$. Then $X/\Phi$ is a connected (that is, $X$ admits no $\Phi$-invariant disconnection) spatial diamond, and for any geometric point $\bar{x}$ of $X$, the map
	\[ \pfg(X/\Phi, \bar{x}) \to \prod_{i=1}^n \pfg(X_i, \bar{x}) \]
is an isomorphism of profinite groups.
\label{Drinfelds-Lemma}
\end{thm}

\begin{rmk}
The formulation of \thref{Drinfelds-Lemma} above uses the language of diamonds, but this is not strictly necessary for our purposes: we will be interested exclusively in the case where $X_1,\dots,X_n$ are perfectoid spaces, in which case  so are $X$ (compare \thref{notation-for-perfectoid-product-Y}) and $X/\Phi$ (because $\Phi$ acts properly discontinuously on $X$; compare \cite[Corollary 4.3.16]{KedlayaSSS}). We will see 
an explicit example of this phenomenon in Subsection~\ref{subsec:mod-p-essential-surjectivity}.
\end{rmk}

\section{$(\phi, \Gamma)$-modules and representations}

In this section, we use Drinfeld's lemma for diamonds to make an initial link between multivariate $(\phi, \Gamma)$-modules and representations of $G_{K, \Delta}$, culminating in the following result.

\begin{thm}[see \thref{tilde-equivalence}, \thref{o-to-o-tilde}]
The category of continuous $\Z_p$-representations of $G_{K, \Delta}$ is equivalent to the category of \'etale $(\phi_\Delta, \Gamma_{K, \Delta})$-modules over either $\calO_{\calE_\Delta(K)}$ or $\tcalO_{\calE_{\Delta}(K)}$. Moreover, these equivalences are exact.
\end{thm}
 This will then be refined in later sections to include the rings $\calO^\dagger_{\calE_\Delta(K)}$ and
$\tcalO^\dagger_{\calE_{\Delta}(K)}$.

\subsection{Product constructions for Huber rings}

We start by adapting the construction of the ring $R_\Delta$ (\thref{notation-for-product-rings-1}) to form topological products of arbitrary perfectoid rings.
\begin{defn}
For $(R,R^+)$ a perfect Huber pair of characteristic $p$, an $R^+$-module is said to be \df{almost zero} if it is annihilated by every topologically nilpotent element of $R^+$. For any fixed pseudo-uniformizer $\varpi$ of $R$, it is sufficient to require the module to be annihilated by $\varpi^{p^{-n}}$ for all nonnegative integers $n$.

For example, since $R$ is perfectoid (\thref{remark:perfect implies uniform}), 
$(R,R^\circ)$ is also a Huber pair (\thref{remark-plus-ring-in-uniform-case}),
but the inclusion $R^+ \subseteq R^\circ$ need not be an equality.
However, the quotient $R^\circ/R^+$ is almost zero.
\end{defn}

\begin{nota}
For $i = 1,\ldots,n$, let $(R_i, R_i^+)$ be a perfect Huber pair of characteristic $p$. For $i = 1,\ldots,n$, choose a pseudo-uniformizer $\varpi_i$ of $R_i$ and define
	\begin{align*}
	\varpi &:= \varpi_1 \cdots \varpi_n \\
	R_0 &:= (R_1^+ \otimes_{\F_p} \cdots \otimes_{\F_p} R_n^+)^{\wedge}_{(\varpi_1,\dots,\varpi_n)} \\
	R &:= R_0[\varpi^{-1}].
	\end{align*}
Note that $R_0$ is perfect; it is also complete for both the $(\varpi_1,\dots,\varpi_n)$-adic topology and for the finer
$(\varpi_1 \cdots \varpi_n)$-adic topology. By equipping $R_0$ with the latter topology, we may give $R$ the structure of a perfect Huber ring containing $R_0$ as an open subring; by \thref{remark:perfect implies uniform}, $R$ is also uniform and hence perfectoid.

The choice of the $\varpi_i$ has no ultimate effect on either $R_0$ or $R$; moreover, $R$ does not depend (either algebraically or topologically) on the choice of the $R_i^+$ within $R_i$. We refer to $R$ as the \df{completed tensor product} of $R_1,\dots,R_n$. As for ordinary tensor products, any continuous endomorphism of $R_i$ extends naturally to $R$ so as to fix $R_j$ for $j \neq i$.

In the setting of \thref{notation-for-product-rings-1}, $R_\Delta$ with the perfectoid topology coincides with the completed tensor product of $R_{\alpha_1},\dots,R_{\alpha_n}$.
\label{notation-for-perfectoid-product-rings}
\end{nota}

In \thref{notation-for-perfectoid-product-rings}, it is not immediately clear that $R_0$ is an integrally closed subring of $R$, which would then ensure that $(R,R_0)$ is a perfectoid Huber pair. For our purposes, it will be sufficient to check something slightly weaker.

\begin{lem} \label{L:circ subring of product}
In \thref{notation-for-perfectoid-product-rings}, 
the quotient $R^\circ/R_0$ is killed by $(\varpi_1 \cdots \varpi_n)^{p^{-m}}$ for every nonnegative integer $m$.
\end{lem}
\begin{pf}
If $x \in R^\circ$, then by definition the sequence $\{x^{p^m}\}_m$ is topologically bounded.
Since $R_0$ is open,
there exists a single nonnegative integer $h$ such that $(\varpi_1 \cdots \varpi_n)^h x^{p^m} \in R_0$
for all $m \geq 0$. But since $R_0$ is perfect, this implies that $(\varpi_1 \cdots \varpi_n)^{hp^{-m}} x \in R_0$ for all $m \geq 0$. Since $hp^{-m} \to 0$ as $m \to \infty$, this implies the claim.
\end{pf}

\subsection{Mod-$p$ representations: full faithfulness}
\label{mod-p-representations}

As in the usual theory of $(\varphi, \Gamma)$-modules, the desired statement about $\Z_p$-representations will be deduced from a corresponding statement about torsion representations, and most of the key ideas appear already in the study of $\F_p$-representations. We correspondingly start by proving the following theorem, which will occupy us for the entirety of this and the next two subsections.

\begin{thm}
The category of continuous $\F_p$-representations of $G_{K, \Delta}$ is equivalent to the category of \'etale $(\phi_\Delta, \Gamma_{K, \Delta})$-modules over $R_\Delta = \tcalO_{\calE_\Delta(K)} / p \tcalO_{\calE_\Delta(K)}$.
\label{tilde-equivalence-mod-p}
\end{thm}

In this subsection, we formulate a more general version of this statement (\thref{cor-4.3.16}), then
in that context produce a fully faithful functor from Galois representations to $(\phi, \Gamma)$-modules (\thref{D-fully-faithful}). Essential surjectivity of this functor will be established in Subsection~\ref{subsec:mod-p-essential-surjectivity},
modulo a technical result about perfectoid spaces which we defer to Subsection~\ref{subsec:H0 of Y}.

We make the following observations:
	\begin{enumerate}
	\item The completion of the field $K(\mu_{p^\infty})$ is perfectoid with tilt $E$
(see \thref{tilted-extensions}).
	\item The absolute Galois groups of $K(\mu_{p^\infty})$, its completion, and $E$ are isomorphic. 
    \item The action of $G_{K(\mu_{p^\infty})}\simeq G_E$ on $E\alg$ extends to an action of $G_K$ by functoriality of tilting. This leads to an action of $\Gamma_K$ on $E$.
	\end{enumerate}
Thus it suffices to establish an equivalence between the category of representations of $H_{K, \Delta}$ and the category of \'etale $\phi_\Delta$-modules over $R_\Delta$, for then the action of $\Gamma_{K, \Delta}$ allows us to recover the categories of $\F_p$-representations of $G_{K, \Delta}$ and \'etale $(\phi_\Delta, \Gamma_{K, \Delta})$-modules over $R_\Delta$. 
In fact we prove a slightly more general statement; the following is \cite[Corollary 4.3.16]{KedlayaSSS}, but here we fill in many details of the proof.

\begin{thm}
Let $F_1, \ldots, F_n$ be perfectoid fields of characteristic $p$. Let $\phi_{\alpha_i}$ act on $F_i$ via the absolute Frobenius. 
Let $R$ be the completed tensor product of $F_1,\dots,F_n$ as per \thref{notation-for-perfectoid-product-rings}. Then the functor $D$ defined in \thref{defn-D-functor} below yields an equivalence of categories between the category of continuous $\F_p$-representations of
\[ G_{F_\Delta} := G_{F_1} \times \cdots \times G_{F_n} \]
and the category of \'etale $\phi_\Delta$-modules over $R$ (that is, finite projective $R$-modules $M$ having commuting semilinear bijective actions of $\phi_{\alpha_1}, \ldots, \phi_{\alpha_n}$).
\label{cor-4.3.16}
\end{thm}

\begin{defn}
For each $i$, fix a completed algebraic closure $\bar{F}_i$ of $F_i$ and identify $G_{F_i}$ with $\Gal(\bar{F}_i  / F_i)$.
Let $\bar{R}$ be the completed tensor product of $\bar{F}_1, \dots, \bar{F}_n$ as per \thref{notation-for-perfectoid-product-rings}.
The ring $\bar{R}$ has a natural action of $G_{F_\Delta}$. Moreover, whenever $F_1=F_2=\dots=F_n$ is the tilt of $K(\mu_{p^\infty})$ for some finite extension $K  / \Qp$, then the action of $G_{F_\Delta}\simeq \HKD$ extends to an action of $\GKD$.
\label{defn-product-of-alg-closures}
\end{defn}

\begin{lem}
We have $\bar{R}^{G_{F_\Delta}} = R$.
\label{invariant-subring-field-product}
\end{lem}
\begin{pf}
Note that $\Spa(\bar{R}, \bar{R}^+) \to \Spa(R,R^+)$ is a pro-\'etale covering; the claim thus follows from the fact that the structure sheaf on $\Spa(R,R^+)$ is a sheaf also for the pro-\'etale topology
\cite[Theorem~3.5.5]{KedlayaLiu2}.
\end{pf}

\begin{defn}
We define a functor $D$ 
which maps an $\F_p$-representation $V$ of $G_{F_\Delta}$ to the $R$-module
	\[ D(V) := (V \otimes_{\F_p} \bar{R})^{G_{F_\Delta}}, \]
where $G_{F_\Delta}$ acts diagonally on the tensor product. For each $i$, we define an action of $\phi_{\alpha_i}$ on $D(V)$ by
	\[ \phi_{\alpha_i}(x \otimes a) = x \otimes \phi_{\alpha_i}(a). \]
This defines commuting semilinear bijective actions of the $\phi_{\alpha_i}$ on $D(V)$.
\label{defn-D-functor}
\end{defn}

\begin{prop}
For $V$ an $\F_p$-representation of $G_{F_\Delta}$, the module $D(V)$ is a finite, projective $R$-module.
\end{prop}
\begin{pf}
Fix a particular choice of representation $V$. Since $V$ is a finite-dimensional $\F_p$-vector space, and thus a finite set, the kernel of the action of $G_{F_\Delta}$ must be an open subgroup. Thus there exist finite, Galois extensions $E_i  / F_i$ in $\bar{F}_i$ such that, for $S$ the completed tensor product of $E_1,\dots,E_n$,
	\[ D(V) = (V \otimes_{\F_p} S)^{\Gal(E_1 / F_1) \times \cdots \times \Gal(E_n / F_n)}. \]
By the Normal Basis Theorem, $E_i \otimes_{F_i} E_i$ has an $F_i$-basis of elements $\sigma_{ij}(e_i) \otimes \sigma_{ik}(e_i)$ for $E_i \otimes_{F_i} E_i$, where $\Gal(E_i  / F_i) = \{\sigma_{i1}, \ldots, \sigma_{in}\}$. Since the maps $\sigma_{ij}$ are $F_i$-linear, we can assume that the element $e_i$ has norm at most 1. By taking $e_{ij} = \sigma_{ij}(e_i)$, it follows that the elements $e_{ij} \otimes \sigma_{ik}(e_{ij})$ form a basis of $E_i \otimes_{F_i} E_i$. Letting $s_j = e_{1j} \otimes \cdots \otimes e_{nj}$, we conclude that the elements $s_j \otimes \sigma(s_j)$ form an $R$-basis of $S \otimes_R S$, where $\sigma$ runs over the elements of
	\[\Gal(E_\Delta  / F_\Delta) = \Gal(E_1  / F_1) \times \cdots \times \Gal(E_n  / F_n).\]
The map given by
	\[ v \otimes s \otimes \sigma(s) \mapsto s \otimes \sigma(v) \otimes \sigma(s) \]
with respect to this basis is a descent datum
	\[ (V \otimes_{\F_p} S) \otimes_R S \to S \otimes_R (V \otimes_{\F_p} S) \]
with respect to the faithfully flat map $R \to S$ (faithfully flat since $S$ is a finite, free $R$-module).
It follows that
	\[ D(V) \otimes_R S \simeq V \otimes_{\F_p} S, \]
and thus that $D(V)$ is a finite, projective $R$-module.
\end{pf}

\begin{prop}
The functor $D$ is fully faithful.
\label{D-fully-faithful}
\end{prop}
\begin{pf}
Recall that for $\F_p$-representations $V$ and $W$, the set $\Hom_{\F_p}(V, W)$ is itself an $\F_p$-representation of $G_{F_\Delta}$, and by \thref{hom-phi-module}, $\Hom_R(D(V), D(W))$ is a $\phi_\Delta$-module over $R$. The morphisms $V \to W$ are those elements of $\Hom_{\F_p}(V, W)$ which are fixed by $G_{F_\Delta}$, while the morphisms $D(V) \to D(W)$ are those elements of $\Hom_R(D(V), D(W))$ that are fixed by $\phi_\Delta$. Using the fact that
	\[ D(\Hom_{\F_p}(V, W)) \simeq \Hom_R(D(V), D(W)), \]
we reduce the problem of showing that $D$ is fully faithful to showing that
	\[ V^{G_{F_\Delta}} \simeq D(V)^{\phi_\Delta} \]
for an arbitrary $\F_p$-representation $V$ playing the role of $\Hom_{\F_p}(V, W)$. We have
	\begin{align*}
	V^{G_{F_\Delta}} &= (V \otimes_{\F_p} \bar{R})^{\phi_\Delta, G_{F_\Delta}} \\
	&\simeq (D(V) \otimes_R \bar{R})^{G_{F_\Delta}, \phi_\Delta} \\
	&= D(V)^{\phi_\Delta}
	\end{align*}
by the proof of \thref{invariant-subring-field-product}, as desired.
\end{pf}

\subsection{Mod-$p$ representations: essential surjectivity}
\label{subsec:mod-p-essential-surjectivity}

Continuing with our study of mod-$p$ representations, we next address essential surjectivity of the functor in \thref{cor-4.3.16}. 
Let $D$ be a finitely presented (and hence projective by \thref{phi-module-flat}) $R$-module with commuting semilinear bijective actions of $\phi_{\alpha_1}, \ldots, \phi_{\alpha_n}$. Let $\phi := \phi_{\alpha_1} \circ \cdots \circ \phi_{\alpha_n}$. 

\begin{rmk}
We have seen in the proof of \thref{D-fully-faithful} that if $D \simeq D(V)$ for some $\F_p$-representation $V$ of $G_{F_\Delta}$, then at the level of \'etale sheaves we have
$V \simeq D^{\varphi_\Delta}$; it is thus natural to try to show that the expression on the right gives rise to a representable functor on $\FEt(R)$. As a first step, we consider $\varphi$-invariants instead of $\varphi_\Delta$-invariants in \thref{represented-mod-p} below.
\end{rmk}

\begin{prop}
The functor that maps a finite \'etale $R$-algebra $T$ to the $\F_p$-vector space $(D \otimes_R T)^{\phi}$ is represented by a finite \'etale $R$-algebra $S$---that is, there is a natural isomorphism
	\begin{equation}
	(D \otimes_R T)^{\phi} \simeq \Hom_R(S, T)
	\label{represents}
	\end{equation}
for all finite \'etale $R$-algebras $T$.
Moreover, $S$ carries natural actions of $\phi_{\alpha_1},\dots,\phi_{\alpha_n}$.
\label{represented-mod-p}
\end{prop}
\begin{pf}
We show this in analogy with \cite[Lemma 3.2.6]{KedlayaLiu} (but correcting some errors therein). Supposing first that $D$ is a free $R$-module, let $e_1, \ldots, e_r$ be a basis for $D$. Let $A$ be the matrix of the action of $\phi$ on $D$ with respect to this basis; since $D$ is \'etale 
the matrix $A$ is invertible. We observe that each element of $D \otimes_R T$, where $T$ is any finite \'etale $R$-algebra, can be written in the form $t_1 e_1 + \cdots + t_r e_r$ for some elements $t_i \in T$, and such an element belongs to $(D \otimes_R T)^{\phi}$ if and only if
	\[ t_i = \sum_{j=1}^r t_j^p a_{ij} \]
for all $i$, that is, if
	\[ \begin{pmatrix}
	t_1 \\ \vdots \\ t_r
	\end{pmatrix}
	= A
	\begin{pmatrix}
	t_1^p \\ \vdots \\ t_r^p
	\end{pmatrix}. \]
Let $S := R[T_1, \ldots, T_r] / ((T_i^p)_i - A^{-1}(T_i)_i)$. Then $S$ is a finite $R$-algebra of degree $p^r$, and the Jacobian matrix $(\partial f_i / \partial T_j)_{i,j}$ has full rank, where the $f_i$ are the polynomials defining the ideal $((T_i^p)_i - A^{-1}(T_i)_i)$. Thus by the Jacobian criterion, $S$ is \'etale over $R$. Since this construction is independent of the choice of basis, we can glue 
to obtain $S$ in the general case.

Now each partial Frobenius $\phi_{\alpha_i}$ induces an $R$-module isomorphism
	\[ D \otimes_{R, \phi_{\alpha_i}} R \to D, \]
so it follows that
	\begin{equation}
	\Hom_R(S, T) \simeq (D \otimes_{R, \phi_{\alpha_i}} R \otimes_R T)^{\phi} \simeq (D \otimes_R (\phi_{\alpha_i}^\ast T))^{\phi}.
	\label{pullback-represents}
	\end{equation}
Combining \eqref{represents} and \eqref{pullback-represents}, we can conclude that
	\[ \Hom_R(S, \phi_{\alpha_i}^\ast T) \simeq \Hom_R(S, T) \]
for all \'etale $R$-algebras $T$, and in particular $\Hom_R(S, \phi_{\alpha_i}^\ast S) \simeq \Hom_R(S, S)$. Let $f_i\colon S \to \phi_{\alpha_i}^\ast S$ be the map corresponding to $\id_S$ under this correspondence. As the composition
	\[ S \to \phi_{\alpha_1}^\ast S \to \phi_{\alpha_2}^\ast \phi_{\alpha_1}^\ast S \to \cdots \to \phi_{\alpha_n}^\ast \cdots \phi_{\alpha_1}^\ast S \]
is just the isomorphism $S \to \phi^\ast S$, it follows that each $f_i$ is an isomorphism. The inverse of $f_i$ thus gives us a semilinear action of $\phi_{\alpha_i}$ on $S$, for which the composition $\phi_{\alpha_1} \circ \cdots \circ \phi_{\alpha_n}$ is the absolute Frobenius.
\end{pf}

We now tie this construction to Drinfeld's lemma, discovering a hitch along the way whose resolution we postpone to Subsection~\ref{subsec:H0 of Y}.

\begin{defn}
 Let $X_i := \Spd(F_i)$ and $\bar{x}_i := \Spd(\bar{F}_i)$. Let $X := X_1 \times \cdots \times X_n$, and let $\bar{x}$ be a geometric point of $X$ lying over each $\bar{x}_i$. By Drinfeld's Lemma (\thref{Drinfelds-Lemma}), we have
	\[ \pfg(X/\Phi, \bar{x}) \simeq \prod_{i=1}^n \pfg(X_i, \bar{x}) \]
as profinite groups, and by \thref{pfg-galois-group}, the right-hand side is isomorphic to $G_{F_\Delta}$.

If it were the case that we could identify $X$ with the diamond associated to $\Spa(R,R^\circ)$, we could then identify the finite \'etale $R$-algebra $S$ with a finite \'etale cover of $X$ carrying natural actions of $\varphi_{\alpha_1},\dots,\varphi_{\alpha_n}$, then use Drinfeld's lemma to see that such a cover must be dominated by a cover formed by taking a product of finite extensions of the $F_i$. Unfortunately, the situation turns out to be a bit more subtle than this.
\label{defn:drinfelds-lemma-application}
\end{defn}

\begin{prop}
We have a natural identification $X = Y^\diamond$ for
	\[ Y := \set{\abs{\ } \in \Spa(R, R^\circ)}{\abs{\bar{\pi}_i}^m \to 0, m \to \infty \ \forall\ i} \]
	viewed as an open subspace of $\Spa(R, R^\circ)$ by writing it 
as the increasing union of the affinoid subspaces
\begin{equation} \label{eq:presentation of Y as union}
	U_m := \set{\abs{\ } \in \Spa(R, R^\circ)}{\abs{\bar{\pi}_i}^m \leq \abs{\bar{\pi}_j}\ \forall\ i, j}.
\end{equation}
\label{prop-X-within-Spa-R}
\end{prop}
\begin{pf}
We show that $X = Y^\diamond$ by showing that every morphism $\Spa(T, T^+) \to \Spa(R, R^\circ)$ with image in $Y$ (for $\Spa(T, T^+)$ an affinoid adic space) corresponds to a tuple of morphisms $\Spa(T, T^+) \to \Spa(F_i, F_i^\circ)$. That is to say, every morphism to $\Spa(R, R^\circ)$ with image in $Y$ factors through $\Spa(F_1, F_1^\circ) \times \cdots \times \Spa(F_n, F_n^\circ)$, identifying $Y^\diamond$ with $\Spd(F_1) \times \cdots \times \Spd(F_n) = X$.

Supoose that $f\colon \Spa(T, T^+) \to \Spa(R, R^\circ)$ has image in $Y$. Then for each $i$, $f$ maps $\bar{\pi}_i$ to a topologically nilpotent element of $T$
and so the composition
$F_i \to R \to T$ is continuous (even though $F_i \to R$ is not). This gives us a family of morphisms $f_i\colon \Spa(T, T^+) \to \Spa(F_i, F_i^\circ)$.

 Conversely, suppose that we have maps $f_i\colon \Spa(T, T^+) \to \Spa(F_i, F_i^\circ)$ for all $i$.
 Because the elements $f_i(\bar{\pi}_i)$ are all topologically nilpotent,
 the induced map $F_1^\circ \otimes_{\F_p} \cdots \otimes_{\F_p} F_n^\circ \to T$ extends to $R_0$. We thus recover a continuous map $R \to T$
 and thus a map  $f\colon \Spa(T, T^+) \to \Spa(R, R^\circ)$.
\end{pf}

\begin{rmk}
Proposition~\ref{prop-X-within-Spa-R} implies that the complement of $Y$ in $\Spa(R,R^\circ)$ is substantial; it is a sort of ``boundary'' consisting of points at which the power-bounded elements $\bar{\varpi}_1,\dots,\bar{\varpi}_n$ are not all individually topologically nilpotent, but at least one of them is (ensuring that their product is also).
For example, when $n=2$, the boundary disconnects into two pieces depending on whether $\bar{\varpi}_1$ or $\bar{\varpi}_2$ is topologically nilpotent.

The identification of $Y$ with a subspace of $\Spa(R,R^\circ)$ induces a morphism $R \to H^0(Y, \calO)$ which is injective but not surjective. For example, 
\[
\sum_{n=1}^\infty \bar{\varpi}_1^{p^n} \bar{\varpi}_2^{-n}
\]
is not contained in the image of $R$.

This has the following effect on the proof of essential surjectivity in \thref{cor-4.3.16}. 
Given the finite \'etale $R$-algebra $S$ from Subsection~\ref{represented-mod-p}, we can pull it back to a finite \'etale cover $Y$ of $X/\Phi$,
then apply Drinfeld's lemma to choose finite extensions $E_i$ of $F_i$ with the property that for $T$ the completed product of the $E_i$,
we can pull back the morphism $\Spa(T,T^\circ) \to \Spa(R,R^\circ)$ to $X/\Phi$ to obtain a finite \'etale cover that dominates $Y$.
This is enough to produce a candidate representation $V$, but not enough to construct an isomorphism $D(V) \simeq D$; for this, we need to show that
$\Spa(T,T^\circ)$ dominates $\Spa(S,S^\circ)$.
\label{remark-boundary-of-X}
\end{rmk}

Our workaround for the preceding issue is given by the following statement, whose proof we defer to
Subsection~\ref{subsec:H0 of Y}.
\begin{prop}[see \thref{o-plus-of-Y}]
The natural map $R^\circ \to H^0(Y, \calO^+)$ is injective and almost surjective (that is, its cokernel is killed by  $(\bar{\pi}_1 \cdots \bar{\pi}_n)^{p^{-m}}$ for all nonnegative integers $m$).
\end{prop}

\begin{cor}
We have
\begin{equation}
R^{\varphi_\Delta} = \Fp.
\label{RmodPhi-connected}
\end{equation}
\end{cor}
\begin{pf}
Any element of $R^{\varphi_\Delta}$ maps to a $\varphi_\Delta$-invariant locally constant function $|Y| \to \F_p$.
By \thref{prop-X-within-Spa-R}, if such an element were not in $\F_p$, then it would imply the existence of a disconnection of $X/\Phi$, which would contradict \thref{Drinfelds-Lemma}.
\end{pf}

Let $Z$ be the diamond corresponding to the adic space $Y \times_{\Spa(R, R^\circ)} \Spa(S, S^\circ)$. The set $V := Z \times_X \bar{x}$ carries an action of $\pfg(X/\Phi, \bar{x}) \simeq G_{F_\Delta}$, by definition, and in fact we can say more: 

\begin{prop}
The set $V$ has the structure of an $\F_p$-representation of $G_{F_\Delta}$.
\label{prop-vector-space-structure}
\end{prop}
\begin{pf}
Let $e \in (D \otimes_R S)^{\phi} \simeq \Hom_R(S, S)$ be the element corresponding to the identity map, and let $\iota_1, \iota_2$ be the two natural inclusions of $(D \otimes_R S)^{\phi}$ into $(D \otimes_R S \otimes_R S)^{\phi}$. Then the element
	\[ \iota_1(e) + \iota_2(e) \in (D \otimes_R S \otimes_R S)^{\phi} \simeq \Hom_R(S, S \otimes_R S) \]
induces an addition law
	\[ \Spa(S, S^\circ) \times_{\Spa(R, R^\circ)} \Spa(S, S^\circ) = \Spa(S \otimes_R S, (S \otimes_R S)^+) \to \Spa(S, S^\circ) \]
(where $(S \otimes_R S)^+$ denotes the integral closure of $S^\circ \otimes_{R^\circ} S^\circ$ in $S \otimes_R S$). This in turn induces an addition law on $\Spd(S \otimes_R \bar{k})$, where $\bar{k}$ is an algebraically closed field such that $\bar{x} = \Spec(\bar{k})$. It remains to observe that $Z \times_X \bar{x} \simeq \Spd(S \otimes_R \bar{k})$ and that this addition law induces an $\F_p$-vector space structure, which by virtue of its naturality commutes with the action of $G_{F_\Delta}$. (Note that $p$-fold addition corresponds to the map $S \to S$ obtained by post-composing the $p$-fold coaddition law $S \to S \otimes_R \cdots \otimes_R S$ with the multiplication map $S \otimes_R \cdots \otimes_R S \to S$; this composition corresponds to the element $pe = 0$ in $(D \otimes_R S)^{\phi}$.)
\end{pf}

As described in \thref{remark-boundary-of-X}, it remains to confirm that $D(V) \simeq D$.
The crucial case is when $V$ is a trivial representation of $G_{F_\Delta}$.

\begin{prop}
If $V$ is a trivial representation of $G_{F_\Delta}$, then $D$ is a trivial $\phi_\Delta$-module of rank $\dim_{\Fp}V$, and in particular is isomorphic to $D(V)$.
\label{trivial-trivial}
\end{prop}
\begin{pf}
In this case, we have
	\[ Y \times_{\Spa(R, R^\circ)} \Spa(S, S^\circ) \simeq \coprod Y. \]
That is, the map $\Spa(S, S^\circ) \to \Spa(R, R^\circ)$ splits completely after pullback to $Y$.

With notation as in \eqref{eq:presentation of Y as union}, let $R_m := H^0(U_m, \calO)$ and $R_m^+ := H^0(U_m, \calO^+)$; then by \thref{o-plus-of-Y}, $R^\circ \to \varprojlim_m R^+_m$ is an almost isomorphism and $R \to \varprojlim_m R_m$ is injective. Let $\widetilde{S}$ be the pushforward of the structure sheaf on $\Spa(S, S^\circ)$ to $\Spa(R, R^+)$, and let $S_m := H^0(U_m, \widetilde{S}) = S \otimes_R R_m$; since $S$ is a finite, \'etale $R$-algebra, $S_m$ is a finite, \'etale $R_m$-algebra. Consider the $\Phi$-invariant idempotents corresponding to the decomposition
	\[ \calO(\Spa(S, S^\circ)) \times_{\calO(\Spa(R, R^\circ))} \calO(Y) \simeq \bigoplus \calO (Y). \]
The restrictions from $Y$ to $U_m$ induce isomorphisms
	\[ S_m \simeq S \otimes_R R_m \simeq \bigoplus R_m. \]
The given idempotents induce compatible systems of idempotents in $S_m$, and thus the idempotents in question belong to $\varprojlim_m S_m$.

Fix a presentation of $S$ as a direct summand of a finite, free $R$-module. This gives a choice of coordinates in $R$ for each element in $S$, which in turn gives a choice of coordinates in $\varprojlim_m R_m$ for each element of $\varprojlim_m S_m$. We remark that an element in $\varprojlim_m S_m$ belongs to $S$ if an only if each of its coordinates belongs to $R$. Now, let $U$ be a rational subspace of $\Spa(R, R^\circ)$ in $Y$ containing a fundamental domain for the action of $\Phi$. As $U$ is quasicompact, the restriction of each idempotent to $H^0(U, \widetilde{S})$ must be an element with bounded coordinates, so there exists an $m$ such that the coordinates belong to $\bar{\pi}_n^{-m} H^0(U, \calO^+)$. (Note that as there is no preferred choice of valuation on $H^0(U, \calO)$, we measure elements against powers of the single element $\bar{\pi}_n$, which is topologically nilpotent on $U$.) On the other hand, since $\phi = \phi_{\alpha_1} \circ \cdots \circ \phi_{\alpha_n}$ acts trivially on $Y$, the group $\Phi$ acts on $Y$ via its quotient modulo $\phi$, which can be generated by the classes of $\phi_{\alpha_1}, \ldots, \phi_{\alpha_{n-1}}$; in particular, the coordinates of our idempotents when restricted to $\psi(U)$ will still belong to $\bar{\pi}_n^{-m} H^0(\psi(U), \calO^+)$ for any $\psi \in \Phi$. Thus after glueing we see that the coordinates of our idempotents belong to $\bar{\pi}_n^{-m} (\bar{\pi}_1 \cdots \bar{\pi}_n)^{-1} R^\circ \subseteq R$. It follows that our $\Phi$-invariant idempotents belong to $S$. This shows that $S \simeq \bigoplus R$; by \eqref{RmodPhi-connected}, this is a decomposition of $S$ into $\Phi$-connected components.

The components in this decomposition correspond to various $\varphi_\Delta$-equivariant homomorphisms $S \to R$. Via 
the natural isomorphism \eqref{represents} (with $T=R$), these in turn correspond to elements of $D^{\varphi_\Delta} \subseteq D^{\phi}$. That is, we have a natural injective morphism $V \to D^{\varphi_\Delta}$ of $\Fp$-vector spaces.

Let $r$ be the rank of $D$.
By the naturality of the previous construction, we have a commutative diagram
\[
\xymatrix{
\Hom_{\Fp}(\Fp^{\oplus r}, V) \times \Hom_{\Fp}(V, \Fp^{\oplus r}) \ar[r] \ar[d] & \Hom_{\Fp}(\Fp^{\oplus r}, \Fp^{\oplus r}) \ar[d] \\
\Hom_R(R^{\oplus r}, D)^{\phi_\Delta} \times \Hom_R(D, R^{\oplus r})^{\phi_\Delta} \ar[r] & \Hom_R(R^{\oplus r}, R^{\oplus r})^{\phi_\Delta},
}
\]
in which the vertical arrows are injective and the horizontal arrows denote composition: $(f,g) \mapsto g \circ f$.
Consequently, if we choose any isomorphism $f\colon \Fp^{\oplus r} \simeq V$ of $\Fp$-vector spaces and take $g$ to be its inverse,
applying the left vertical arrow to $(f,g)$ yields a pair of $\Phi$-equivariant morphisms between $R^{\oplus r}$ and $D$ whose composition is the identity on $R^{\oplus r}$. Since both $R^{\oplus r}$ and $D$ are finite projective $R$-modules of rank $r$, this implies that the two morphisms are indeed inverses on both sides, so $D$ is a trivial $\phi_\Delta$-module as claimed.
\end{pf}

We finish by returning to the case when $V$ is not necessarily a trivial $\F_p$-representation of $G_{F_\Delta}$. 

\begin{prop}
The module $D$ arises from the representation $V$ in general.
\label{D-from-diamond}
\end{prop}
\begin{pf}
Let $E_i  / F_i$ be finite extensions such that $V$ is a trivial representation of $G_{E_\Delta}$. Let
$R_E$ be the completed tensor product of $E_1,\dots,E_n$.
By \thref{D-fully-faithful} and \thref{trivial-trivial}, there is a canonical $\phi_\Delta$-equivariant and $G_{E_\Delta}$-equivariant isomorphism 
\[
D \otimes_R R_E \simeq V \otimes_{\Fp} R_E;
\]
that is, $D \otimes_R R_E \simeq D(V|_{E_\Delta})$.
By canonicality, this isomorphism is also $G_{F_\Delta}$-equivariant for the diagonal action on both sides; by taking $G_{F_\Delta}$-invariants, we obtain an isomorphism $D \simeq D(V)$.
%
%
\end{pf}
This completes the proof of \thref{cor-4.3.16}, and thus of \thref{tilde-equivalence-mod-p}, modulo the comparison between
$R^\circ$ and $H^0(Y, \calO^+)$; we give this next.

\subsection{$H^0(Y, \calO^+)$: a geometric interlude}
\label{subsec:H0 of Y}

In this subsection, we complete the comparison between $R^\circ$ and $H^0(Y, \calO^+)$ alluded to above (\thref{o-plus-of-Y}); this will in turn complete the proof of \thref{cor-4.3.16}, and will appear at a corresponding point
in the study of mod-$p^n$ representations (\thref{trivial-trivial-modpn}). For this calculation, it is convenient to work more generally, by taking products not of fields but of arbitrary perfect, analytic Huber rings. We note in passing that much of the previous section can be carried over to this level of generality
(as in the comparison of $\varphi$-modules and $\Zp$-local systems given in \cite[\S 8]{KedlayaLiu}), but that will not be necessary for our present purposes.

We start by describing a special case in detail, in order to illustrate why such an assertion is reasonable to expect. This example was suggested in \cite[Remark 4.3.18]{KedlayaSSS}.

\begin{ex}
\label{2-case}
Let $\calo_{F_1}$ be the $T_1$-adic completion of $\F_p\llbracket T_1 \rrbracket [T_1^{p^{-\infty}}]$, which we denote by $\F_p\llbracket T_1^{p^{-\infty}} \rrbracket$. Let $F_1 := \calo_{F_1}[T_1^{-1}]$. Let
	\begin{align*}
	A^+ &:= \calo_{F_1}\langle T_2^{p^{-\infty}} \rangle \\
		&= \left\{\sum_{j \in (1/p^\infty)\Z_{\geq 0}} \left( \sum_{i \in (1/p^\infty)\Z_{\geq 0}} a_{ij} T_1^i \right) T_2^j \, \middle| \, \text{$\forall\ k$, $\exists$ finitely many nonzero } \right. \\
			&\qquad \left. \text{$a_{ij}\in\Fp$ with $i + j \leq k$; $\Abs[T_1]{\sum_{i \in (1/p^\infty)\Z_{\geq 0}} a_{ij} T_1^i} \to 0$ as $j \to \infty$} \right\} \\
	A &:= F_1 \langle T_2^{p^{-\infty}} \rangle \\
		&= \left\{\sum_{j \in (1/p^\infty)\Z_{\geq 0}} \left( \sum_{i \in (1/p^\infty)\Z_{\geq m_0}} a_{ij} T_1^i \right) T_2^j\, \middle| \, m_0 \in \Z; \text{ $\forall\ k$, $\exists$ finitely many} \right. \\
			&\qquad \left. \text{nonzero $a_{ij}$ with $i + j \leq k$; $\Abs[T_1]{\sum_{i \in (1/p^\infty)\Z_{\geq m_0}} a_{ij} T_1^i} \to 0$ as $j \to \infty$} \right\}
	\end{align*}
Consider 
	\[ Y := \set{\abs{\ } \in \Spa(A, A^+)}{0 < \abs{T_2} < 1} \]
as an adic space by identifying it with the union of the rational open subspaces
	\begin{align*}
	U_m &:= \set{\abs{\ } \in \Spa(A, A^+)}{\abs{T_1^m} \leq \abs{T_2} \leq \abs{T_1^{1/m}}} \\
		&= U(T_1^m/T_2, T_2/T_1^{1/m})
	\end{align*}
as $m$ runs over the powers of $p$. Let $(\calO, \calO^+)$ denote the structure sheaves of $\Spa(A, A^+)$. Then $\calO^+(U_m)$ is the completion of the integral closure of $A^+[T_1^m/T_2, T_2/T_1^{1/m}]$. That is,
	\begin{align*}
	\calO^+(U_m) &= A \langle T_1^m/T_2, T_2/T_1^{1/m} \rangle^+ \\
		&= \left\{ \sum_{r, s \in (1/p^\infty)\Z_{\geq 0}} \left( \sum_{i, j \in (1/p^\infty)\Z_{\geq 0}} a_{ijrs} T_1^i T_2^j \right) \left(\frac{T_1^m}{T_2}\right)^r \left(\frac{T_2}{T_1^{1/m}} \right)^s \, \right| \\
			&\quad \text{$\forall\ k$, $\exists$ finitely many nonzero $a_{ijrs}\in\Fp$ with $i + j + r(m - 1) +$} \\
		&\quad \left. \text{$ s(m - 1) \leq k$; $\Abs[T_1]{\sum_{i \in (1/p^\infty)\Z_{\geq 0}} a_{ijrs} T_1^i} \to 0$ as $j \to \infty$;} \right. \\
		    &\quad \left. \text{$\Abs[T_1]{\sum_{i,j \in (1/p^\infty)\Z_{\geq 0}} a_{ijrs} T_1^i T_2^j} \to 0$ as $\max\{r, s\} \to \infty$} \right\}
	\end{align*}
Now, for any choices of coefficients $a_{ij}$ such that for each $k$, at most finitely many $a_{ij}$ with $i + j \leq k$ are nonzero, we have 
	\begin{align*}
	\sum_{i, j \in (1/p^\infty)\Z_{\geq 0}} a_{ij} T_1^i T_2^j &= \sum_{i, j \in (1/p^\infty)\Z_{\geq 0}} a_{ij} T_1^{i + \lfloor j \rfloor / m} T_2^{j - \lfloor j \rfloor} (T_2/T_1^{1/m})^{\lfloor j \rfloor} \\
		&= \sum_{s \in \Z_{\geq 0}} \left(\sum_{u, v \in (1/p^\infty)\Z_{\geq 0}, v < 1} a_{(u - s/m)v} T_1^u T_2^v\right) \left(\frac{T_2}{T_1^{1/m}}\right)^s \\
		&\in \calO^+(U_m),
	\end{align*}
so
	\begin{align*}
	\F_p&\llbracket T_1^{p^{-\infty}}, T_2^{p^{-\infty}}\rrbracket \subseteq \calO^+(U_m) \\ &\subseteq \left\{ \sum_{i, j \in (1/p^\infty)\Z} a_{ij} T_1^i T_2^j\; \middle|\; i + mj \geq 0, mi + j \geq 0, \right. 
		\left. \vphantom{\sum_{i, j \in (1/p^\infty)\Z}} \forall\ k\ \exists\ \text{finitely many} \right. \\
		&\quad \left. \vphantom{\sum_{i, j \in (1/p^\infty)\Z}}\text{nonzero $a_{ij}$ with $i + j \leq k$}\right\}.
	\end{align*}
	
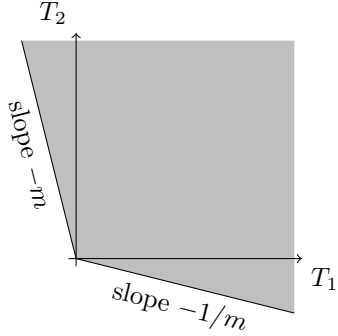
	\begin{center}
	\begin{tikzpicture}
	\fill[lightgray] (0,0) -- (2.9, -.725) -- (2.9,2.9) -- (-.725, 2.9) -- cycle;
	\draw (-.725, 2.9) --node[below, sloped]{slope $-m$} (0,0) --node[below, sloped]{slope $-1/m$} (2.9,-.725);
	\draw[->] (-.1,0) -- (3,0) node[below right] {$T_1$};
	\draw[->] (0,-.1) -- (0,3) node[above left] {$T_2$};
	\end{tikzpicture}
	\end{center}
It follows that
	\begin{align*}
	\calO^+(Y) &= \calO^+\left(\bigcup_{m=p^k} U_m\right) \\
		&= \bigcap_{m=p^k} \calO^+(U_m) \\
		&= \F_p\llbracket T_1^{p^{-\infty}}, T_2^{p^{-\infty}} \rrbracket.
	\end{align*}
In particular, in this case the map $R^\circ \to H^0(Y, \calO^+)$ is a genuine isomorphism, not just an almost isomorphism.
\end{ex}

We now set notation so as to work with a more general completed tensor product of perfectoid rings.

\begin{nota}
With notation as in \thref{notation-for-perfectoid-product-rings}, 
for $i = 1,\ldots,n$, let $\varphi_i\colon R \to R$ be the ring homomorphism obtained by tensoring the Frobenius morphism on $R_i$ with the identity morphism on $R_j$ for $j \neq i$. Put 
	\[ Y := \set{v \in \Spa(R,R^\circ)}{v(\varpi_1)^m, \dots, v(\varpi_n)^m \to 0, m \to \infty}, \]
	viewed as an open subspace of $\Spa(R,R^+)$ as in \thref{prop-X-within-Spa-R}; by the same proof, we have
	$\Spd(R_1) \times \cdots \times \Spd(R_n) = Y^\diamond$.
	\label{notation-for-perfectoid-product-Y}
\end{nota}

In order to simulate arguments using ordinary power series, we recall the following construction.
\begin{defn}
Let $k$ be a ring and let $\Gamma$ be a totally ordered abelian group (written multiplicatively). The ring $k\llparenthesis \Gamma \rrparenthesis$ of \df{Hahn--Mal'cev--Neumann generalized power series} is the set of functions $\Gamma \to k$ whose support is a well-ordered subset of $\Gamma$. (In view of the multiplicative notation for $\Gamma$, a well-ordered subset here must be taken to be one containing no infinite \df{increasing} sequence.) We represent the function $\gamma \mapsto c_\gamma$ as a formal sum $\sum_{\gamma \in \Gamma} c_\gamma [\gamma]$; then $k\llparenthesis \Gamma \rrparenthesis$ forms a ring with respect to the operations
	\begin{align*}
	\sum_{\gamma} c_{\gamma} [\gamma] + \sum_{\gamma} d_{\gamma} [\gamma] &= 
	\sum_{\gamma} (c_\gamma + d_\gamma) [\gamma] \\
	\sum_{\gamma} c_{\gamma} [\gamma] \times \sum_{\gamma} d_{\gamma} [\gamma] &= 
	\sum_\gamma \left( \sum_{\gamma' \gamma'' = \gamma} c_{\gamma'} d_{\gamma''} \right) [\gamma].
	\end{align*}
The condition on well-ordered supports is needed to establish that multiplication is well-defined; one must first check that the sum over $\gamma', \gamma''$ is always finite, and then that the support of the resulting sum is again well-ordered. See for example \cite[\S 4]{kaplansky1}.

From its construction, the ring $k\llparenthesis \Gamma \rrparenthesis$ comes equipped with a natural valuation: the function assigning to every nonzero formal sum the maximal element of its support (which exists by the well-ordered condition).
\end{defn}

\begin{lem}[Kaplansky] \label{L:Kaplansky}
Let $F$ be a nonarchimedean field of equal characteristics, with value group $\Gamma$ and residue field $k$, which is algebraically closed and \df{maximally complete}. (The latter condition means that there is no nontrivial extension of $F$ with the same value group and residue field as $F$.) Then there exists an isomorphism $F \simeq k\llparenthesis \Gamma \rrparenthesis$ of fields with valuation.
\end{lem}
\begin{pf}
See \cite[Theorem~7]{kaplansky1}.
\end{pf}

\begin{lem} \label{L:series calculation}
Let $\R^+$ denote the multiplicative group of positive real numbers.
Suppose that $n=2$ and there exists an isomorphism $R_1 \simeq k_1\llparenthesis \R^+ \rrparenthesis$ for some perfect ring $k_1$ (which need not be a field). Then the map $R^\circ \to H^0(Y, \calO^+)$ is an isomorphism.
\end{lem}
\begin{pf}
Fix a power-multiplicative norm $\left|\,  \right|$ defining the topology on $R_2$.
We extend this norm to $R_2 \otimes_{\F_p} k_1$ as follows: choose a basis of $k_1$ as an $\F_p$-vector space, use it to view $R_2 \otimes_{\F_p} k_1$ as a direct sum of copies of $R_2$, and take the supremum norm of the coordinates. This does not depend on the choice of the basis of $k_1$.

We may explicitly describe $R$ as a certain set of formal sums $\sum_{\gamma \in \R^+} c_\gamma [\gamma]$ with coefficients in $R_2 \otimes_{\F_p} k_1$. These sums must satisfy the following conditions.
\begin{itemize}
\item[(i)]
For any neighborhood $U$ of $0$ in $R_2$, the set of $\gamma \in \R^+$ such that $c_{\gamma} \notin U$ is well-ordered; and for each $\gamma_0 \in \R^+$, there exists a finitely generated $\F_p$-submodule $M$ of $R_2$ such that the quantities $c_\gamma$ for $\gamma \geq \gamma_0$ all belong to $(U + M) \otimes_{\F_p} k_1$.
\item[(ii)]
The set $\set{\gamma \in \R^+}{c_\gamma \neq 0}$ is bounded above in $\R^+$ (but not necessarily well-ordered), and the set $\set{c_\gamma}{\gamma \in \R^+}$ is bounded in $R_2 \otimes_{\F_p} k_1$.
\end{itemize}

In terms of these formal sums, for each $r > 0$ we may define a norm $v_r$ on $R$ by the formula
	\[
	v_r\left( \sum_{\gamma} c_\gamma [\gamma] \right) = \max\set{\gamma^r \abs{c_\gamma}}{\gamma \in \R^+}.
	\]
Now observe that the following two families of seminorms are mutually cofinal; that is, any member of one family is eventually dominated by all members of the other family.
\begin{itemize}
\item[(a)]
The functions $\max\{v_r, v_{1/r}\}$ as $r \to 0^+$.
\item[(b)]
The supremum norms over the subspaces 
$$
\set{|\ | \in Y}{\abs{[\gamma]} \leq \abs{\pi_2} \leq \abs{[\gamma^{-1}]}}
$$
as $\gamma \to 0^+$.
\end{itemize}
In (b), the subspaces in question have union $Y$; therefore, taking the inverse limit of the completions of $R$ with respect to these seminorms yields $H^0(Y, \calO)$. By the cofinal property, the same is true of (a); consequently, $H^0(Y, \calO)$ may be identified with the set of formal sums $\sum_{\gamma \in \R^+} c_\gamma [\gamma]$ with coefficients in $R_2 \otimes_{\F_p} k_1$ subject to the same condition (i) and a slightly weaker version of (ii):
\begin{itemize}
\item[(ii)$'$]
For any $r>0$, the set $\set{\gamma^r \abs{c_\gamma}}{ \gamma \in \R^+}$ is bounded above in $\R$.
\end{itemize}
As for $H^0(Y, \calO^+)$, we may characterize as the set of $s \in H^0(Y, \calO)$ for which each of the supremum norms in (b) is at most 1; by cofinality, this is equivalent to requiring that $v_r(s) \leq 1$ for all $r>0$.

From these descriptions, it is immediate that $R^\circ \to H^0(Y, \calO^+)$ is injective.
To check surjectivity, note that if $s = \sum_{\gamma} c_\gamma [\gamma] \in H^0(Y, \calO^+)$, then for each $\gamma \in \R^+$ for which $c_\gamma \neq 0$, we must have $\gamma^r \abs{c_\gamma} \leq 1$. In particular, if $\gamma > 1$ then $c_\gamma = 0$, while if $\gamma \leq 1$ then $c_\gamma \in R_2^\circ$. From this it is apparent that $s \in R^\circ$, so $R^\circ \to H^0(Y, \calO^+)$ is surjective.
(See \cite[Lemma~5.1]{KedlayaSimple} for a related argument.)
\end{pf}

\begin{prop}
With notation as in 	\thref{notation-for-perfectoid-product-Y},
the map $R^\circ \to H^0(Y, \calO^+)$ is injective and almost surjective.
\label{o-plus-of-Y}
\end{prop}
\begin{pf}
Since $R_i$ is allowed to be a perfect ring, not necessarily a field, we may deduce the general case by repeatedly applying the case $n = 2$; we thus assume $n = 2$ hereafter. Since $R_1$ is uniform, it embeds as a closed subring of a product of perfectoid fields; by \thref{L:Kaplansky}, each of these fields can be embedded into a ring of the form $k\llparenthesis \R^+\rrparenthesis$ for some perfect field $k$. Let $R_1'$ be the set of bounded elements (for the supremum norm) in the product of these rings
and put $R_1'' = R_1' \widehat{\otimes}_{R_1} R_1'$. Let $R'$ (resp. $R''$) be the completed tensor product of $R_1'$ (resp. $R_1''$) with $R_2$.
Let $Y'$ (resp. $Y''$) be the subspace of the spectrum of $R'$ (resp. $R''$) defined as in \thref{notation-for-perfectoid-product-Y}.
In the commutative diagram
\[
\xymatrix{
R^\circ \ar[r] \ar[d] & R_1^{\prime \circ} \widehat{\otimes}_{\F_p} R_2^\circ \ar@<-2pt>[r] \ar@<+2pt>[r] \ar[d] & R_1^{\prime \prime \circ} \widehat{\otimes}_{\F_p} R_2^\circ \ar[d] \\
R^\circ \ar[r] \ar[d] & R^{\prime \circ} \ar@<-2pt>[r] \ar@<+2pt>[r] \ar[d] & R^{\prime \prime \circ} \ar[d] \\
H^0(Y, \calO^+) \ar[r] & H^0(Y', \calO^+) \ar@<-2pt>[r] \ar@<+2pt>[r] & H^0(Y'', \calO^+)
}
\]
the vertical arrows between the first and second rows are almost isomorphisms by \thref{L:circ subring of product};
the top and bottom rows are almost equalizer diagrams by the fact that $\calO^+$ induces an almost acyclic sheaf for the  \df{v-topology}  \cite[Definition~3.8.5]{KedlayaSSS}
on any affinoid perfectoid space \cite[Theorem~3.5.5]{KedlayaLiu2};
and
the lower vertical arrow in the second column is injective and almost surjective by \thref{L:series calculation}
(applied to each factor of $R_1'$).

From this, we may first deduce that the lower vertical arrow in the first column is injective.
By similar logic, the lower vertical arrow in the third column is injective;
we may thus deduce that the lower vertical arrow in the first column is almost surjective.
\end{pf}

\begin{rmk}
The preceding proposition is conceptually related to the \df{perfectoid Riemann extension theorem}. The first such statement is due to Scholze \cite{Scholze_TorsionCohomology}; similar statements are used in the proofs of the direct summand conjecture by Andr\'e \cite{Andre2,Andre1} and Bhatt \cite{Bhatt}.

To make this link explicit, we describe an alternate proof of \thref{o-plus-of-Y} in the case $n=2$ (which again suffices for the general case), in which we appeal directly to the extension theorem as formulated by Bhatt.
Consider the covering of $X := \Spa(R,R^\circ)$ by the two rational subspaces
\[
U_1 := \set{|\ | \in X}{\abs{\pi_1} \leq \abs{\pi_2}}, \qquad
U_2 := \set{|\ | \in X}{\abs{\pi_2} \leq \abs{\pi_1}}.
\]
Since $R$ is perfectoid, it is sheafy, so it will suffice to check that for $i=1,2$ the map
\[
H^0(U_i, \calO^+) \to H^0(U_i \cap Y, \calO^+)
\]
is injective and an almost isomorphism. By symmetry, we may restrict to the case $i=1$. In this case,
we define $R' := H^0(U_1, \calO) = R \langle T \rangle/(\pi_2 T - \pi_1)$,
so that $R^{\prime \circ}$ is almost isomorphic to $H^0(U_1, \calO^+)$. We may then make the identification
\[
U_1 \cap Y \simeq \set{|\ | \in \Spa(R', R^{\prime \circ})}{\abs{\pi_2} > 0}
\]
and apply \cite[Theorem~4.2]{Bhatt} to deduce that $R^{\prime \circ} \to H^0(U_i \cap Y, \calO^+)$
is injective and its cokernel is almost zero. (As an aside, we note that the proof of \cite[Theorem~4.2]{Bhatt} can be conceptualized rather well in terms of the pictures from \thref{2-case}.)
\end{rmk}

\subsection{Mod-$p^m$ representations}

We return to the notation of Subsection \ref{mod-p-representations}
and prove the following.

\begin{thm}
The category of continuous $\Z_p$-representations of $G_{K, \Delta}$ is equivalent to the category of \'etale $(\phi_\Delta, \Gamma_{K, \Delta})$-modules over $\tcalO_{\calE_\Delta(K)}$.
\label{tilde-equivalence}
\end{thm}

As in the proof of \thref{tilde-equivalence-mod-p}, 
this reduces to the following statement.
\begin{thm}
With notation as in \thref{cor-4.3.16},
there exists an equivalence of categories between the category of continuous $\Z_p$-representations of $G_{F_\Delta}$
and the category of \'etale $\varphi_\Delta$-modules over $W(R)$.
\label{cor-4.3.16-integral}
\end{thm}

Moreover, it suffices to check that the analogous equivalence holds modulo $p^m$ on both sides for each $m$, since we can then obtain the desired result by taking limits. That is, we need to prove the following.

\begin{thm}
With notation as in \thref{cor-4.3.16},
for each positive integer $m$,
the functor $D$ defined in Definition~\ref{defn-functor-D-mod-p^m} below yields
an equivalence of categories between the category of continuous representations of $G_{F_\Delta}$ on finite free $\Z/p^m \Z$-modules
and the category of \'etale $\varphi_\Delta$-modules over $W(R)/p^m W(R)$.
\label{cor-4.3.16-modpn}
\end{thm}

Our method of proof is similar to that of the case $m=1$; consequently, we summarize a few points that  do not differ substantially from the previous case.

\begin{defn}
Define a functor $D$ from the category of continuous representations of $G_{F_\Delta}$ on finite projective $\Z / p^m \Z$-modules to the category of $W(R) / p^m W(R)$-modules with commuting semilinear bijective actions of $\phi_{\alpha_1}, \ldots, \phi_{\alpha_n}$ by
	\[ D(V) := (V \otimes_{\Z_p} W(\bar{R}))^{G_{F_\Delta}}, \]
where $\bar{R}$ is as in 
\thref{defn-product-of-alg-closures}.
\label{defn-functor-D-mod-p^m}
\end{defn}

\begin{prop}
The module $D(V)$ is a finite, projective $W(R) / p^m W(R)$-module.
\end{prop}
\begin{pf}
As in the proof of \thref{tilde-equivalence-mod-p}, for any fixed $V$ we can find finite, Galois extensions $E_i  / F_i$ such that the action of $G_{F_\Delta}$ on $V$ factors through $\Gal(E_\Delta  / F_\Delta)$ and, for $S$ the completed tensor product of $E_1,\dots,E_n$,
	\[ D(V / p^m V) = (V / p^m V \otimes_{\Z_p} W(S))^{\Gal(E_\Delta  / F_\Delta)}, \]
As $S$ is free over $R$, it is in particular faithfully finite flat over $R$, and in fact it is faithfully finite \'etale, 
so by \cite[Proposition~5.5.4]{KedlayaLiu}, the extension $W(R) \to W(S)$ is faithfully finite \'etale. It follows that the extension $W(R) / p^m W(R) \to W(S) / p^m W(S)$ is faithfully finite \'etale, and in particular faithfully flat. The map
	\begin{gather*}
	(V \otimes_{\Z_p} S) \otimes_R S \to S \otimes_R (V \otimes_{\Z_p} S) \\
	v \otimes s \otimes \sigma(s) \mapsto s \otimes \sigma(v) \otimes \sigma(s)
	\end{gather*}
induces a map
	\[ (V / p^m V \otimes_{\Z_p} W(S)) \otimes_{W(R)} W(S) \to W(S) \otimes_{W(R)} (V / p^m V \otimes_{\Z_p} W(S)), \]
which is then a descent datum with respect to the faithfully flat extension $W(R) / p^m W(R) \to W(S) / p^m W(S)$. It follows by faithfully flat descent for modules that
	\[ D(V) \otimes_{W(R)} W(S) \simeq V \otimes_{\Z_p} W(S), \]
and in particular $D(V)$ is a finite, projective $W(R) / p^m W(R)$-module for each $m$.
\end{pf}

\begin{prop}
The functor $D$ is fully faithful.
\label{ff-modpn}
\end{prop}
\begin{pf}
The proof is analogous to the proof of \thref{D-fully-faithful}.
\end{pf}

We now turn to essential surjectivity. Let $D$ be a finitely presented (hence projective by \thref{phi-module-flat}) module over $W(R) / p^m W(R)$ carrying commuting, semilinear actions of $\phi_{\alpha_1}, \ldots, \phi_{\alpha_n}$, and let $\phi := \phi_{\alpha_1} \circ \cdots \circ \phi_{\alpha_n}$.

\begin{prop}
There exist a finite \'etale $R$-algebra $T_m$ such that $D \otimes_{W(R)} W(T_m)$ admits  a $\phi$-invariant basis
over $W(T_m)/pW(T_m)$. (Note that we are not currently looking for a $\phi_\Delta$-invariant basis.)
\label{extension-with-basis}
\end{prop}
\begin{pf}
We proceed by induction on $m$, the case $m=1$ being included in \thref{tilde-equivalence-mod-p}. 
Suppose that the induction hypothesis holds with $m$ replaced by $m-1$;
then there exists a finite \'etale $R$-algebra $T_{m-1}$ such that  $D \otimes_{W(R)} W(T_{m-1}) / p^{m-1} W(T_{m-1})$
admits a $\phi$-invariant basis $e_1,\dots,e_r$. Choose lifts of these elements to $D \otimes_{W(R)} W(T_{m-1})$; by Nakayama's lemma, they still form a basis over $W(T_{m-1}) / p^{m} W(T_{m-1})$. Let 
\[
F := (f_{ij})_{i,j} \in M_r(W(T_{m-1}) / p^{m} W(T_{m-1}))
\]
be the matrix of the action of $\phi$ on $D \otimes_{W(R)} W(T_{m-1})$ with respect to this basis, i.e.\ $\phi(e_j) = \sum_{i=1}^r f_{ij} e_i$ for each $j$. As this basis is $\phi$-invariant modulo $p^{m-1}$, we have $F = I + p^{m-1} A$ for some matrix $A$ uniquely determined modulo $p$. We want to find a matrix $B$ (with coefficients potentially in a larger ring) which conjugates the basis $(e_i)_i$ to a $\phi$-invariant basis; such a matrix would satisfy $B^{-1} F \phi(B) = I$, in other words $(I + p^{m-1} A) \phi(B) = B$. Since this is satisfied modulo $p^{m-1}$ for $B = I$, we can look for a matrix of the form $I + p^{m-1} C$, where the entries of $C = ([\bar{c}_{ij}])_{i,j}$ are Teichm\"uller lifts. Thus we are looking for elements $\bar{c}_{ij}$  in some finite \'etale $T_{m-1}$-algebra $T_m$  with the property that
	\[(I + p^{m-1} A)(I + p^{m-1}([\bar{c}_{ij}^p])_{i,j}) = I + p^{m-1} ([\bar{c}_{ij}])_{i,j}\]
in $M_r(W(T_m) / p^{m} W(T_m))$. Expanding the left-hand side, and noting that $p^{2(m-1)} = 0$ in this ring,  we reduce to finding elements $\bar{c}_{ij} \in T_m$ such that
	\[
	A + ([\bar{c}_{ij}^p])_{i,j} = ([\bar{c}_{ij}])_{i,j}.
	\]
In other words, we need only adjoin to $T_{m-1}$ roots of the equations $x^p - x - \bar{a}_{ij}$ for each entry $a_{ij}$ of $A$. The resulting ring $T_{m}$ is then finite \'etale over $T_{m-1}$, and consequently over $R$.
\end{pf}

\begin{prop}
There exist a finite \'etale $R$-algebra $S_m$ which represents the functor $F_m$ defined as follows:
	\[ F_m(T) := (D \otimes_{W(R)} W(T))^{\phi}. \] 
\end{prop}
\begin{pf}
Consider $T_m$ as in \thref{extension-with-basis} and define $D_m := D \otimes_{W(R)} W(T_m)$. 
By construction, $D_m$ admits a $\phi$-invariant basis $e_1, \ldots, e_r$ over $W(T_m) / p^m W(T_m)$. We first show that there exists a finite \'etale $R$-algebra $Q_m$ which represents the functor that maps a finite \'etale $T_m$-algebra $T$ to the $\Z / p^m \Z$-module $(D_m \otimes_{W(T_m)} W(T))^{\phi}$. We can write
	\[ D_m = (W(T_m) / p^m W(T_m)) e_1 + \cdots + (W(T_m) / p^m W(T_m)) e_r. \]
Let $T$ be a finite \'etale $T_m$-algebra. Then
	\begin{multline*}
	    (D_m \otimes_{W(T_m)} W(T) / p^m W(T))^{\phi} \\= (W(T) / p^m W(T))^{\phi} e_1 + \cdots + (W(T) / p^m W(T))^{\phi} e_r;
	\end{multline*}
the ring $(W(T) / p^m W(T))^{\phi}$ consists of one copy of $\Z / p^m \Z$ per connected component of $T$.

Let $Q_m := T_m[\bar{x}_{ij}\ |\ i = 1, \ldots r, \quad j = 1, \ldots, m] / ( \cdots )$, where $\cdots$ are the relations necessary to ensure
	\[ [\bar{x}_{i1}] + p [\bar{x}_{i2}] + p^2 [\bar{x}_{i3}] + \cdots + p^{m-1} [\bar{x}_{im}] \]
is a $\phi$-invariant element in $W(Q_m) / p^m W(Q_m)$ for each $i$. Now
	\[ t_1 e_1 + \cdots + t_r e_r \in (D_m \otimes_{W(T_m)} W(T))^{\phi} \]
if and only if each element $t_i = \sum_{j=0}^{m-1} p^j [\bar{t}_{ij}]$ is fixed by $\phi$, so that such an element corresponds to the mapping
	\[ W(Q_m) / p^m W(Q_m) \to W(T) / p^m W(T) \]
sending the element
	\[ [\bar{x}_{i1}] + p [\bar{x}_{i2}] + p^2 [\bar{x}_{i3}] + \cdots + p^{m-1} [\bar{x}_{im}] \]
to $t_i$ for each $i$. Thus we have
	\begin{align*}
	(D_m &\otimes_{W(T_m)} W(T))^{\phi_R} \simeq \Hom_{T_m}(Q_m, T) \\
	&\simeq \Hom_{W(T_m) / p^m W(T_m)}(W(Q_m) / p^m W(Q_m), W(T) / p^m W(T)),
	\end{align*}
as desired.

To summarize, so far we have a finite \'etale $T_m$-algebra $Q_m$ which represents the functor
	\[ T \mapsto (D \otimes_{W(T_m)} W(T) / p^m W(T))^{\phi} \]
on finite \'etale $T_m$-algebras.
We want to show that $Q_m$ has the form $S_m \otimes_r T_m$, where $S_m$ is a finite \'etale $R$-algebra that represents the functor $F_m$ on finite \'etale $R$-algebras. On the same category, define the functors
	\begin{align*}
	G_1(T) &:= \Hom_{T_m}(Q_m \otimes_{T_m, 1} (T_m \otimes_R T_m), T) \\
		&\simeq \Hom_{T_m}(Q_m, T) \times \Hom_{T_m, 1}(T_m \otimes_R T_m, T) \\
		&\simeq F_{T_m}(T) \times \Hom_{T_m, 1}(T_m \otimes_R T_m, T),
	\end{align*}
and similarly define
	\begin{align*}
	G_2(T) &:= \Hom_{T_m}(Q_n \otimes_{T_m, 2} (T_m \otimes_R T_m), T) \\
		&\simeq F_{T_m}(T) \times \Hom_{T_m, 2}(T_m \otimes_R T_m, T),
	\end{align*}
where the 1 and 2 indicate whether $T_m \otimes_R T_m$ is considered as a $T_m$-algebra via its first or second component. In the case when $D$ is a trivial $\phi$-module, the functor which maps a finite \'etale $R$-module $T$ to $(D \otimes_{W(R)} W(T))^{\phi}$ is represented by $R^{p^{rm}}$, where $r$ is the rank of $D$, and so we have a canonical identification $Q_m \simeq R^{p^{rm}} \otimes_R T_m$. In this case there is a canonical descent datum
	\[ (R^{p^{rm}} \otimes_R T_m) \otimes_R T_m \to T_m \otimes_R (R^{p^{rm}} \otimes_R T_m). \]
This induces a natural isomorphism $G_2 \to G_1$ (by identifying $Q_m \otimes_{T_m, 1} (T_m \otimes_R T_m)$ with $Q_m \otimes_R T_m$ and $Q_m \otimes_{T_m, 2} (T_m \otimes_R T_m)$ with $T_m \otimes_R Q_m$), which in turn induces natural isomorphisms $F_{T_m} \to F_{T_m}$ and $\Hom_{T_m, 2}(T_m \otimes_R T_m, \ast) \to \Hom_{T_m, 1}(T_m \otimes_R T_m, \ast)$. But these two functors do not depend on $Q_m$; thus in the case when $D$ need not be trivial, these isomorphisms may be used to construct a natural isomorphism $G_2 \to G_1$, which by Yoneda's lemma arises from an isomorphism 
\[
Q_m \otimes_{T_m, 1} (T_m \otimes_R T_m) \to Q_m \otimes_{T_m, 2} (T_m \otimes_R T_m).
\]
This isomorphism provides a descent datum to which we may apply faithfully flat descent 
to deduce that $Q_m = S_m \otimes_R T_m$ for some finite \'etale $R$-algebra $S_m$. We show that $S_m$ represents the functor $F_m$. Consider the equalizer diagram
	\[ T \to T \otimes_R T_m \rightrightarrows T \otimes_R T_m \otimes_R T_m. \]
This induces equalizer diagrams
	\[ \Hom_{T_m}(Q_m, T) \to \Hom_{T_m}(Q_m, T \otimes_R T_m) \rightrightarrows \Hom_{T_m}(Q_m, T \otimes_R T_m \otimes_R T_m) \]
and
	\begin{align*}
	(D &\otimes_{W(R)} W(T))^{\phi} \to (D \otimes_{W(R)} W(T) \otimes_{W(R)} W(T_m) )^{\phi} \\
		&\rightrightarrows (D \otimes_{W(R)} W(T)  \otimes_{W(R)} W(T_m)  \otimes_{W(R)} W(T_m) )^{\phi}.
	\end{align*}
Since $Q_m$ represents the functor $F_{T_m}$, the right two objects in each equalizer diagram are isomorphic. It follows that the left objects are also isomorphic. But
	\[ \Hom_{T_m}(Q_m, T) \simeq \Hom_{T_m}(S_m \otimes_R T_m, T) \simeq \Hom_R(S_m, T), \]
as desired.
\end{pf}

Define $X$ and $\bar{x}$ as in \thref
{defn:drinfelds-lemma-application}, so that $\pfg(X/\Phi, \bar{x}) \simeq G_{F_\Delta}$. 
Define the perfectoid space $Y$ as in Proposition~\ref{prop-X-within-Spa-R}, so that $X \simeq Y^\diamond$.
With notation as in the proof of \thref{trivial-trivial}, let $Z$ be the diamond corresponding to the adic space $Y \times_{\Spa(R, R^\circ)} \Spa(S_m, S_m^+)$.
As in \thref{prop-vector-space-structure}, the set $V := Z \times_X \bar{x}$ has an action of the group $\pfg(X/\Phi, \bar{x}) \simeq G_{F_\Delta}$ by definition, and starting from the element
	\begin{align*}
	e \in &(D \otimes_{W(R) } W(T_m) )^{\phi} \\
		&\simeq \Hom_{W(R) / p^m W(R)}(W(S_m) / p^m W(S_m), \\
		&\quad W(S_m) / p^m W(S_m) \otimes_{W(R)} W(T_m) / p^m W(T_m))
	\end{align*}
we can construct an addition law on $\Spd(S_m \otimes_R \bar{k}) \simeq V$, under which the $p^m$-fold sum of $e$ corresponds to the element $p^m e = 0$ in $(D \otimes_{W(R)} W(T_m) )^{\phi}$. Thus $V$ carries the structure of a $\Z / p^m \Z$-module with an action of $G_{F_\Delta}$. 

\begin{prop}
If $V$ is a trivial $\Z / p^m \Z$-representation of $G_{F_\Delta}$, then $D$ is a trivial $\phi_\Delta$-module, and in particular is isomorphic to $D(V)$.
\label{trivial-trivial-modpn}
\end{prop}
\begin{pf}
In this case, we have 
	\[ Y \times_{\Spa(R, R^+)} \Spa(S_m, S^+_m) \simeq \coprod Y. \]
As in the proof of \thref{trivial-trivial}, we have $S_m = \bigoplus R$; by \eqref{RmodPhi-connected}, this is a decomposition of $S_m$ into $\Phi$-connected components.
Again as in the proof of \thref{trivial-trivial}, we have a natural injective morphism 
$V \to D^{\varphi_\Delta}$ of $\Z/p^m \Z$-modules; using naturality of the construction, we may deduce that $D$ is a trivial $\phi_\Delta$-module.
\end{pf}

\begin{prop}
The $\phi$-module $D$ arises from $V$ in general.
\label{es-modpn}
\end{prop}
\begin{pf}
The proof is analogous to the proof of \thref{D-from-diamond}.
\end{pf}

This completes the proof of \thref{cor-4.3.16-modpn}:
the functor $D$ is fully faithful by \thref{ff-modpn}
and essentially surjective by \thref{es-modpn}.
This in turn implies \thref{cor-4.3.16-integral}
and \thref{tilde-equivalence}.

\subsection{Descent from $\tcalO$ to $\calO$}

Returning to the case where $F_1=F_2=\dots=F_n$ is the tilt of $K(\mu_{p^\infty})$ for some finite extension $K  / \Qp$, denote by $\tcalO_{\calE_{\Delta}^\mathrm{nr}}$ the ring $W(\bar{R})$.  In the case $K = \Q_p$, the following result was proved previously (by another method) by the third author \cite{Zabradi}.

\begin{thm}
The category of continuous $\Z_p$-representations of $G_{K, \Delta}$ is equivalent to the category of \'etale $(\phi_\Delta, \Gamma_{K, \Delta})$-modules over $\calO_{\calE_\Delta(K)}$.
In particular, by \thref{tilde-equivalence},
base extension of $(\phi_\Delta, \Gamma_{K,\Delta})$-modules
from $\calO_{\calE_\Delta(K)}$ to $\tcalO_{\calE_{\Delta}(K)}$ 
is an equivalence of categories.
\label{o-to-o-tilde}
\end{thm}
\begin{pf}
By \thref{tilde-equivalence}, the functor
	\[ \widetilde{M}(V) := (V \otimes_{\Z_p} \tcalO_{\calE\nr_{\Delta}})^{H_{K, \Delta}} \]
is an equivalence of categories from the category of continuous $\Z_p$-representations of $G_{K, \Delta}$ to the category of \'etale $(\phi_\Delta, \Gamma_{K, \Delta})$-modules over $\tcalO_{\calE_\Delta(K)}$. This functor is the composition of the two functors
	\begin{align*}
	M(V) &:= (V \otimes_{\Z_p} \calO_{\calE\nr_{\Delta}})^{H_{K, \Delta}} \\
	B(M) &:= M \otimes_{\calO_{\calE_\Delta(K)}} \tcalO_{\calE_\Delta(K)},
	\end{align*}
the first taking a continuous $\Z_p$-representation of $G_{K, \Delta}$ to an \'etale $(\phi_\Delta, \Gamma_{K, \Delta})$-module over $\calO_{\calE_\Delta(K)}$, and the second base-changing to $\tcalO_{\calE_\Delta(K)}$. From the fact that $\widetilde{M}$ is an equivalence of categories, it follows that the functor $B$ is essentially surjective; furthermore $B$ is fully faithful, as we assume all $(\phi_\Delta, \Gamma_{K, \Delta})$-modules $M$ are finite projective. Thus $B$ is an equivalence of categories, and it follows that the category of continuous $\Z_p$-representations of $G_{K, \Delta}$ is equivalent to the category of \'etale $(\phi_\Delta, \Gamma_{K, \Delta})$-modules over $\calO_{\calE_\Delta(K)}$.
\end{pf}

\section{Descent for overconvergent Witt vectors}

In this section we show that base extension gives an equivalence between the categories of \'etale $(\phi_\Delta, \Gamma_{K, \Delta})$-modules over $\tcalO_{\calE_\Delta(K)}$ and over $\tcalO^\dagger_{\calE_\Delta(K)}$. 
Our arguments follow closely those in \cite[Section 2.4]{KedlayaNewMethods} with a twist. In that setting, one first descends the action of $\varphi$ from $\tcalO_{\calE_\Delta(K)}$ to $\tcalO^\dagger_{\calE_\Delta(K)}$
and then shows that this automatically causes the action of $\Gamma_K$ to descend. Here, we use a similar method to show that once the action of $\varphi$ is descended, this automatically causes the actions of both $\varphi_\Delta$ and $\Gamma_{K,\Delta}$ to descend.

\begin{thm}[see \thref{tilde-dagger-to-tilde}]
Base extension of \'etale $(\phi_\Delta, \Gamma_{K, \Delta})$-modules from $\tcalO^\dagger_{\calE_\Delta(K)}$ to $\tcalO_{\calE_\Delta(K)}$ is an equivalence of categories. Consequently, both categories are equivalent to the category of continuous representations of $G_{K, \Delta}$ on finitely generated $\Z_p$-modules.
\end{thm}

\begin{rmk}
Note that the rings $\tcalO^\dagger_{\calE_\Delta(K)}$ and $\tcalO_{\calE_\Delta(K)}$ have the same $p^n$-torsion quotients for every positive integer $n$. Consequently, the main content in the equivalence between
\'etale $(\phi_\Delta, \Gamma_{K, \Delta})$-modules over $\tcalO^\dagger_{\calE_\Delta(K)}$ and $\tcalO_{\calE_\Delta(K)}$ 
is at the level of $p$-torsion-free modules.
\label{remark-same-torsion-quotient}
\end{rmk}

\begin{rmk}
With \thref{tilde-equivalence} proven, perfectoid spaces will play no further role in the remainder of the paper.
That said, the arguments given here also yield a corresponding refinement of \thref{cor-4.3.16-modpn} in terms of the ring $W^\dagger(R)$ of \df{overconvergent Witt vectors}: base extension of \'etale $\varphi_\Delta$-modules from $W^\dagger(R)$ to $W(R)$ is an equivalence of categories. We omit further details here.
\label{remark-perfectoid-overconvergent-descent}
\end{rmk}

\subsection{Full faithfulness}

We start by showing full faithfulness of the base extension from  $\tcalO^\dagger_{\calE_\Delta(K)}$ to $\tcalO_{\calE_\Delta(K)}$. We first justify the introduction of a convenient additional hypothesis.

\begin{rmk}
Let $M$ be an \'etale $(\varphi_\Delta, \Gamma_{K,\Delta})$-module over either $\tcalO^\dagger_{\calE_\Delta(K)}$ or $\tcalO_{\calE_\Delta(K)}$. By \thref{tilde-equivalence-mod-p}, $M/pM$ corresponds to a continuous $\F_p$-representation of $G_{K,\Delta}$, and we can find a finite extension $L$ of $K$ such that this representation restricts trivially to
$G_{L,\Delta}$. This means that after base extension from $K$ to $L$, $M$ acquires a basis which modulo $p$ is fixed by $\varphi_\Delta$ and $\Gamma_{L,\Delta}$. Our approach will be to prove everything under the assumption of the existence of such a basis, then return to the general case at the end via faithfully flat descent.
\end{rmk}

\begin{lem}
Let $M$ be an \'etale $\phi$-module over $\tcalO^\dagger_{\calE_\Delta(K)}$ such that the action of $\phi$ is trivial modulo $p$. Then $M^{\phi} = (M \otimes_{\tcalO^\dagger_{\calE_\Delta(K)}} \tcalO_{\calE_\Delta(K)})^{\phi}$.
\label{tilde-fixed-points-mod-p}
\end{lem}
\begin{pf}
The proof is analogous to the proof of \cite[Lemma~2.4.2]{KedlayaNewMethods}.
\end{pf}

\begin{cor}
Let $M$ be an \'etale $\phi_\Delta$-module over $\tcalO^\dagger_{\calE_\Delta(K)}$ such that the action of $\phi_\Delta$ is trivial modulo $p$. Then $M^{\phi_\Delta} = (M \otimes_{\tcalO^\dagger_{\calE_\Delta(K)}} \tcalO_{\calE_\Delta(K)})^{\phi_\Delta}$.
\label{tilde-delta-fixed-points-mod-p}
\end{cor}
\begin{pf}
We note that an element is fixed by $\phi_\Delta$ if and only if it is fixed by $\phi_\alpha$ for all $\alpha \in \Delta$. The result then follows immediately from \thref{tilde-fixed-points-mod-p}.
\end{pf}

\begin{cor}
Base extension of \'etale $\phi_\Delta$-modules which are trivial modulo $p$ from $\tcalO_{\calE_\Delta(K)}^\dagger$ to $\tcalO_{\calE_\Delta(K)}$ is fully faithful.
\label{tilde-base-ext-mod-p-fully-faithful}
\end{cor}
\begin{pf}
This follows directly from \thref{tilde-delta-fixed-points-mod-p}, by letting $\Hom_{\tcalO_{\calE_\Delta(K)}^\dagger}(M,N)$ play the role of $M$, as in the proof of \thref{D-fully-faithful}, since by \thref{hom-phi-module} homomorphisms of $\phi_\Delta$-modules are exactly those elements of $\Hom_{\tcalO_{\calE_\Delta(K)}^\dagger}(M,N)$ that are fixed by $\phi_\Delta$.
\end{pf}

\subsection{Compatibility of descent with extra structures}

If a $\phi$-module descends from $\tcalO_{\calE_\Delta(K)}$ to $\tcalO^\dagger_{\calE_\Delta(K)}$, the following lemma allows us to conclude that the action of some commuting operator (say $\phi_\alpha$ or $\gamma \in \Gamma_\alpha$ for $\alpha \in \Delta$) also descends, assuming that the actions of both $\phi$ and the commuting operator are trivial modulo $p$.

\begin{lem}
Let $\nu$ be an endomorphism of $\tcalO_{\calE_\Delta(K)}$ which commutes with $\phi$ and sends
$\tcalO^\dagger_{\calE_\Delta(K)}$ into itself. 
Let $M$ be an \'etale $\phi$-module over $\tcalO_{\calE_\Delta(K)}$ with a commuting semilinear action of $\nu$.
Suppose that the action of $\phi$ on $M$ is trivial modulo $p$, and that $M$ has a basis $e_1, \ldots, e_d$ with respect to which the matrix $F$ of the action of $\phi$ has entries in $\tcalO^\dagger_{\calE_\Delta(K)}$. Then the matrix of the action of $\nu$ with respect to $e_1,\dots,e_d$ also has entries in $\tcalO^\dagger_{\calE_\Delta(K)}$.
\label{nu-commutes}
\end{lem}
\begin{pf}
By changing the basis $e_1,\dots,e_d$ via a suitable invertible matrix over $\tcalO^\dagger_{\calE_\Delta(K)}$, we can ensure that it is congruent mod $p$ to a basis on which $\phi$ acts trivially mod $p$ (whose existence was hypothesized). That is, we may assume that 
$F$ is congruent to the identity matrix modulo $p$.

Define $\tphi\colon M_d(\tcalO_{\calE_\Delta(K)}) \to M_d(\tcalO_{\calE_\Delta(K)})$ by
	\[ \tphi(B) = F \phi(B) \nu(F)^{-1}, \]
where $\phi$ and $\nu$ are applied to the matrices componentwise. The operator $\tphi$ is $\phi$-semilinear and thus induces the structure of a $\phi$-module on the $\tcalO_{\calE_\Delta(K)}$-module $M_d(\tcalO_{\calE_\Delta(K)})$. As $F$ and $\nu(F)^{-1}$ have entries in $\tcalO^\dagger_{\calE_\Delta(K)}$, the matrix of $\tphi$ on the standard basis of $M_d(\tcalO_{\calE_\Delta(K)})$ does also. Thus $\tphi$ induces a $\phi$-module structure on the $\tcalO^\dagger_{\calE_\Delta(K)}$-module $M_d(\tcalO^\dagger_{\calE_\Delta(K)})$.

Let $N$ be the matrix of the action of $\nu$ on the basis of the $e_i$. As $\nu$ commutes with $\phi$, we have
	\[ F \phi(N) = N \nu(F), \]
and thus $\tphi(N) = N$, which is to say $N \in M_d(\tcalO_{\calE_\Delta(K)})^{\tphi}$. Since 
$F$ is congruent to the identity matrix modulo $p$, the action of $\tphi$ is trivial modulo $p$ on the standard basis of $M_d(\tcalO_{\calE_\Delta(K)})$, so by \thref{tilde-fixed-points-mod-p},
	\[ M_d(\tcalO^\dagger_{\calE_\Delta(K)})^{\tphi} = (M_d(\tcalO^\dagger_{\calE_\Delta(K)}) \otimes_{\tcalO^\dagger_{\calE_\Delta(K)}} \tcalO_{\calE_\Delta(K)})^{\tphi} = M_d(\tcalO_{\calE_\Delta(K)})^{\tphi}, \]
and thus $N$ has entries in $\tcalO^\dagger_{\calE_\Delta(K)}$.
\end{pf}

\subsection{Completion of the proof}

To complete the proof, we must first prove essential surjectivity of base extension from $\tcalO^\dagger_{\calE_\Delta(K)}$ to $\tcalO_{\calE_\Delta(K)}$ for \'etale $\phi_\Delta$-modules which are trivial modulo $p$
(full faithfulness having been proved in \thref{tilde-base-ext-mod-p-fully-faithful}),
then perform a descent to eliminate the extra hypothesis.

\begin{prop}
Let $M$ be an \'etale $\phi_\Delta$-module over $\tcalO_{\calE_\Delta(K)}$ with a basis $e_1, \ldots, e_d$ that is fixed by $\phi$ modulo $p$. Then there exists a basis $e'_1, \ldots, e'_d$ of $M$ on which $\phi_\Delta$ acts via invertible matrices over $\tcalO^\dagger_{\calE_\Delta(K)}$.
\label{base-ext-mod-p-ess-surj}
\end{prop}
\begin{pf}
Let $\phi := \phi_{\alpha_1} \circ \cdots \circ \phi_{\alpha_n}$. Let $F \in \GL_d(\tcalO_{\calE_\Delta(K)})$ be the matrix of the action of $\phi$ on the $e_i$, i.e.\ $\phi(e_j) = \sum_{i=1}^d F_{ij} e_i$. It follows from (the proof of) \cite[Theorem~8.5.3]{KedlayaLiu} that there exists a basis $e'_1, \ldots, e'_n$ of $M$ on which $\phi$ acts via a matrix with entries in $\tcalO^\dagger_{\calE_\Delta(K)}$. Since each of the maps $\phi_\alpha$ commutes with $\phi$, and the action of $\phi$ is trivial modulo $p$, it follows from \thref{nu-commutes} that the matrix of the action of each $\phi_\alpha$ on the $e'_i$ also has entries in $\tcalO^\dagger_{\calE_\Delta(K)}$.
\end{pf}

\begin{thm}
Base extension of \'etale $(\phi_\Delta, \Gamma_{K, \Delta})$-modules from $\tcalO^\dagger_{\calE_\Delta(K)}$ to $\tcalO_{\calE_\Delta(K)}$ is an equivalence of categories. Consequently, both categories are equivalent to the category of continuous representations of $G_{K, \Delta}$ on finitely generated $\Z_p$-modules.
\label{tilde-dagger-to-tilde}
\end{thm}
\begin{pf}
To show that the base extension functor is fully faithful, let $M$ be an \'etale $(\phi_\Delta, \Gamma_{K, \Delta})$-module over $\tcalO^\dagger_{\calE_\Delta(K)}$ (playing the role of $\Hom_{\tcalO^\dagger_{\calE_\Delta(K)}}(M, N)$ as in the proof of \thref{tilde-base-ext-mod-p-fully-faithful}). We need to show that \[M^{\phi_\Delta, \Gamma_{K, \Delta}} = (M \otimes_{\tcalO^\dagger_{\calE_\Delta(K)}} \tcalO_{\calE_\Delta(K)})^{\phi_\Delta, \Gamma_{K, \Delta}}.\]

The space $M \otimes_{\tcalO^\dagger_{\calE_\Delta(K)}} \tcalO_{\calE_\Delta(K)}$ is an \'etale $(\phi_\Delta, \Gamma_{K, \Delta})$-module over $\tcalO_{\calE_\Delta(K)}$, and thus corresponds to a continuous representation of $G_{K, \Delta}$ on a finitely generated $\Z_p$-module $T$. The action of each $G_{K, \alpha}$ on $T/pT$ factors through $\Gal(L  / K)$ for each $\alpha \in \Delta$ for some finite extension $L  / K$. Since $G_L$ acts trivially on $T/pT$, the action of $G_{L,\Delta}$ on $T$ corresponds to an \'etale $(\phi_\Delta, \Gamma_{K, \Delta})$-module over $\tcalO_{\calE_{\Delta}(L)}$ which is trivial modulo $p$. Now since
	\[ (\tcalO_{\calE_{\Delta}(L)})^{\phi_\Delta, \Gamma_{K, \Delta}} = \Z_p \subseteq (\tcalO^\dagger_{\calE_{\Delta}(L)})^{\phi_\Delta, \Gamma_{K, \Delta}}, \]
we have
    \begin{multline*}
        (M \otimes_{\tcalO^\dagger_{\calE_\Delta(K)}} \tcalO_{\calE_\Delta(K)})^{\phi_\Delta, \Gamma_{K, \Delta}} \subseteq (M \otimes_{\tcalO^\dagger_{\calE_\Delta(K)}} \tcalO_{\calE_{\Delta}(L)})^{\phi_\Delta, \Gamma_{K, \Delta}} \\= (M \otimes_{\tcalO^\dagger_{\calE_\Delta(K)}} \tcalO^\dagger_{\calE_{\Delta}(L)})^{\phi_\Delta, \Gamma_{K, \Delta}}.
    \end{multline*}
But as
	\[ \tcalO_{\calE_\Delta(K)} \cap \tcalO^\dagger_{\calE_{\Delta}(L)} = \tcalO^\dagger_{\calE_\Delta(K)}, \]
we have
	\[ (M \otimes_{\tcalO^\dagger_{\calE_\Delta(K)}} \tcalO_{\calE_\Delta(K)})^{\phi_\Delta, \Gamma_{K, \Delta}} \subseteq M^{\phi_\Delta, \Gamma_{K, \Delta}}, \]
as desired.

To show that the functor is essentially surjective, let $M$ be an \'etale $(\phi_\Delta, \Gamma_{K, \Delta})$-module over $\tcalO_{\calE_\Delta(K)}$. We want to show $M$ descends to an \'etale $(\phi_\Delta, \Gamma_{K, \Delta})$-module over $\tcalO^\dagger_{\calE_\Delta(K)}$.

Suppose first that $M$ is killed by $p^m$ for some positive integer $m$. As noted in \thref{remark-same-torsion-quotient},
the rings $\tcalO^\dagger_{\calE_\Delta(K)}/p^m \tcalO^\dagger_{\calE_\Delta(K)}$ and $\tcalO{\calE_\Delta(K)}/p^m \tcalO{\calE_\Delta(K)}$ coincide, so $M$ may be viewed equally well as a module over both of them.

Suppose next that $M$ is flat over $\Z_p$; by applying \thref{phi-module-flat} to the quotients $M/p^m M$, we see that $M$ is projective over $\tcalO_{\calE_\Delta(K)}$.
As above, we can find a finite extension $L  / K$ for which $M \otimes_{\tcalO_{\calE_\Delta(K)}} \tcalO_{\calE_\Delta(L)}$ admits a basis that is fixed by $\varphi$ modulo $p$. 
By \thref{base-ext-mod-p-ess-surj}, $M \otimes_{\tcalO_{\calE_\Delta(K)}} \tcalO_{\calE_\Delta(L)}$ descends to an \'etale $\phi_\Delta$-module $M'$ over $\tcalO^\dagger_{\calE_\Delta(L)}$; by \thref{nu-commutes}, the action of $\Gamma_{K, \Delta}$ also descends. Now $\Gal(L  / K)$ acts on $\tcalO_{\calE_{\Delta}(L)}$, thus on $M \otimes_{\calE_\Delta(K)} \tcalO_{\calE_{\Delta}(L)} = M' \otimes_{\tcalO^\dagger_{\calE_{\Delta}(L)}} \tcalO_{\calE_{\Delta}(L)}$, and therefore on $M'$, since $\Gal(L  / K)$ preserves $\tcalO^\dagger_{\calE_{\Delta}(L)}$. By Galois descent, $M'$ descends to the $\tcalO^\dagger_{\calE_\Delta(K)}$-module $M^\dagger = (M')^{\Gal(L / K)}$.

Suppose finally that $M$ is general. By \thref{tilde-equivalence}, $M$ corresponds to a continuous representation of $G_{K,\Delta}$ on a finitely generated $\Z_p$-module $T$. Let $T_0$ be the torsion submodule of $T$ and put $T_1 = T/T_0$; we then have a short exact sequence
\[
0 \to T_0 \to T \to T_1 \to 0
\]
of continuous $\Z_p$-representations of $G_{K,\Delta}$ in which $T_0$ is killed by $p^m$ for some positive integer $m$ and
$T_1$ is projective. 
Let
\begin{equation} \label{eq:no dagger exact sequence}
0 \to M_0 \to M \to M_1 \to 0
\end{equation}
be the corresponding exact sequence of \'etale $(\phi_\Delta, \Gamma_{K, \Delta})$-modules over $\tcalO_{\calE_\Delta(K)}$.
By the previous paragraph, $M_1$ descends uniquely to an \'etale $(\phi_\Delta, \Gamma_{K, \Delta})$-module $M_1^\dagger$ over $\tcalO^\dagger_{\calE_\Delta(K)}$. 
Moreover, the splittings $M_1 \to M$ of \eqref{eq:no dagger exact sequence} in the category of $\tcalO_{\calE_\Delta(K)}$-modules
form a torsor for the group
\begin{align*}
\Hom_{\tcalO_{\calE_\Delta(K)}}(M_1, M_0)
&= \Hom_{\tcalO_{\calE_\Delta(K)}/p^m\tcalO_{\calE_\Delta(K)}}(M_1/p^m M_1, M_0) \\
&= \Hom_{\tcalO^\dagger_{\calE_\Delta(K)}/p^m\tcalO^\dagger_{\calE_\Delta(K)}}(M_1^\dagger/p^m M_1^\dagger, M_0) \\
&= \Hom_{\tcalO^\dagger_{\calE_\Delta(K)}}(M^\dagger_1, M_0).
\end{align*}
It follows that \eqref{eq:no dagger exact sequence} descends uniquely to an exact sequence
\[
0 \to M_0 \to M^\dagger \to M_1^\dagger \to 0
\]
of \'etale $(\phi_\Delta, \Gamma_{K, \Delta})$-modules over $\tcalO^\dagger_{\calE_\Delta(K)}$.
Thus $M$ descends to $M^\dagger$ as desired.
\end{pf}

\section{Descent for overconvergent power series}

In this section we complete the proof of \thref{mainthm1}
by establishing the following.
\begin{thm}[see \thref{full-equivalence}]
Base extension of \'etale $(\varphi_\Delta, \Gamma_{K, \Delta})$-modules from $\calO^\dagger_{\calE_\Delta(K)}$ to $\calO_{\calE_\Delta(K)}$
is an equivalence of categories. Consequently, both categories are equivalent to the category of continuous representations of $G_{K,\Delta}$ on finitely generated $\Z_p$-modules.
\end{thm}
We then deduce the corresponding equivalence for $\Qp$-representations
(see \thref{full-equivalence-Qp}).

Our arguments follow closely those in \cite[Sections 2.5--2.6]{KedlayaNewMethods}. As in the previous section,
we first prove everything for modules which are trivial modulo $p$, and then reduce to this case by faithfully flat descent.
For one $\alpha \in \Delta$ at a time, we use the action of $\Gamma_{K,\alpha}$ to eliminate fractional powers of $\varpi_\alpha$.

\begin{rmk}
By contrast with \thref{remark-perfectoid-overconvergent-descent}, the arguments of this section do not have a direct analogue for more general perfectoid fields. The extent to which the arguments can be generalized is in fact far from clear; see \cite{KedlayaFrobMod} for a detailed discussion.
\end{rmk}

\subsection{$\Zp$-representations}

In this subsection, we establish \thref{full-equivalence}. As noted above, for most of the proof we will consider only modules that are trivial modulo $p$; keeping in mind the final paragraph of the proof of \thref{tilde-dagger-to-tilde}, 
it will also be harmless to further restrict to projective modules.

Under these conditions, we first establish full faithfulness by applying what we already know about other base extensions.
\begin{lem}
Let $M$ be a projective \'etale $\phi$-module over $\calO^\dagger_{\calE_\Delta(K)}$ such that the action of $\phi$ is trivial modulo $p$. Then $M^{\phi} = (M \otimes_{\calO^\dagger_{\calE_\Delta(K)}} \calO_{\calE_\Delta(K)})^{\phi}$.
\label{fixed-points-mod-p}
\end{lem}
\begin{pf}
By \thref{tilde-fixed-points-mod-p}, we have
	\[ (M \otimes_{\calO^\dagger_{\calE_\Delta(K)}} \calO_{\calE_\Delta(K)})^{\phi} \subseteq (M \otimes_{\calO^\dagger_{\calE_\Delta(K)}} \tcalO_{\calE_\Delta(K)})^{\phi} = (M \otimes_{\calO^\dagger_{\calE_\Delta(K)} }\tcalO^\dagger_{\calE_\Delta(K)})^{\phi}. \]
Since $M$ is a projective module, we have
	\[ (M \otimes_{\calO^\dagger_{\calE_\Delta(K)}} \calO_{\calE_\Delta(K)}) \cap (M \otimes_{\calO^\dagger_{\calE_\Delta(K)}} \tcalO^\dagger_{\calE_\Delta(K)}) = M. \]
The result follows by taking fixed points on both sides.
\end{pf}

\begin{cor}
Base extension of projective \'etale $\phi_\Delta$-modules which are trivial modulo $p$ from $\calO_{\calE_\Delta(K)}^\dagger$ to $\calO_{\calE_\Delta(K)}$ is fully faithful.
\label{base-ext-mod-p-fully-faithful}
\end{cor}
\begin{pf}
This follows directly from \thref{fixed-points-mod-p}, by letting $\Hom_{\calO_{\calE_\Delta(K)}^\dagger}(M,N)$ play the role of $M$, as in the proofs of \thref{D-fully-faithful} and \thref{tilde-base-ext-mod-p-fully-faithful}, since by \thref{hom-phi-module} homomorphisms of $\phi_\Delta$-modules are exactly those elements of $\Hom_{\calO_{\calE_\Delta(K)}^\dagger}(M,N)$ that are fixed by $\phi_\Delta$.
\end{pf}

We now begin the real work of the proof, to establish essential surjectivity under our harmless extra restrictions.
We start by setting up an important decomposition of the ring $R_\Delta(K)$.

\begin{nota}
For $\alpha \in \Delta$, let $X_{\alpha} \in \calO_{\calE_\Delta(K)}$ be the copy of $[\epsilon] - 1$
in the copy of $\calO_\calE$ indexed by $\alpha$. Let $\bar{X}_\alpha$ be the reduction of $X_{\alpha}$ modulo $p$.
\end{nota}

\begin{prop}
For $i = 1, \ldots, n$, let $\bar{T}_i \subseteq R_\Delta(K)$ be the closure of the subgroup generated by
	\[(1 + \bar{X}_{\alpha_1})^{e_1} \cdots (1 + \bar{X}_{\alpha_i})^{e_i} \calO_{\calE_\Delta(K)} / p \calO_{\calE_\Delta(K)}\]
for $e_1, \ldots, e_i \in \Z[p^{-1}] \cap [0,1)$ with $e_i \neq 0$. The natural map
	\[\calO_{\calE_\Delta(K)}/p\calO_{\calE_\Delta(K)} \oplus \bar{T}_1 \oplus \cdots \oplus \bar{T}_n \to
	R_\Delta(K)\]
is an isomorphism of Banach spaces over $\calO_{\calE_\Delta(K)} / p \calO_{\calE_\Delta(K)}$.
\label{dir-sum}
\end{prop}
\begin{pf}
The case when $n = 1$ was proved in \cite[Lemma 2.5.3 and Corollary 2.5.4]{KedlayaNewMethods}, and an analogous proof establishes the general case.
\end{pf}

\begin{rmk}
Let us spell out what \thref{dir-sum} is saying when $K = K_0$.
In this case, we can write elements of $R_\Delta(K)$ as power series in fractional powers of the variables $\bar{\varpi}_{\alpha_i}$ for $i=1,\dots,n$. The summand $\bar{T}_i$ corresponds to those power series in which each monomial includes a nonintegral power of $\bar{\varpi}_{\alpha_i}$, but only integral powers of $\bar{\varpi}_{\alpha_j}$ for any $j > i$. The complement of these corresponds to power series in the $\bar{\varpi}_{\alpha_i}$ themselves,
which constitute precisely $\calO_{\calE_\Delta(K)}/p\calO_{\calE_\Delta(K)}$.
\end{rmk}

We next show that in a strong sense, the action of $\Gamma_{K,\alpha_i}$ has no fixed points on $\bar{T}_i$.
\begin{prop}
For every $\gamma \in \Gamma_{K, \alpha_i}$ of infinite order, with $\alpha_i \in \Delta$, there exists a $c > 0$ such that every $\bar{Y} \in R_\Delta(K)$ can be written uniquely in the form
	\[\bar{U} + (\gamma - 1)(\bar{V}_i) + \sum_{j \neq i} \bar{V}_j \]
with $\bar{U} \in \calO_{\calE_\Delta(K)} / p \calO_{\calE_\Delta(K)}$, all $\bar{V}_j \in \bar{T}_j$, and
	\[\max\{|\bar{U}|', |\bar{V}_1|', \ldots, |\bar{V}_n|'\} \leq c |\bar{Y}|'.\]
\label{decompose-mod-p}
\end{prop}
\begin{pf}
In analogy with \cite[Lemma 2.5.5 and Corollary 2.5.6]{KedlayaNewMethods}, which proved this result in the case $n = 1$, we can show that the map $(\gamma - 1)\colon \bar{T}_i \to \bar{T}_i$ is bijective with bounded inverse. The desired result then follows from \thref{dir-sum}.
\end{pf}

We next lift the preceding definitions and propositions out of characteristic $p$.

\begin{prop}
For $i = 1, \ldots, n$, let $T_i \subseteq \tcalO_{\calE_\Delta(K)}$ be the closure, with respect to the weak topology, of the subgroup generated by
	\[ (1 + X_{\alpha_1})^{e_1} \cdots (1 + X_{\alpha_i})^{e_i} \calO_{\calE_\Delta(K)} \]
for $e_1, \ldots, e_i \in \Z[p^{-1}] \cap [0,1)$ with $e_i \neq 0$. Then for all $\gamma \in \Gamma_{K, \alpha_i}$ of infinite order, there exist $c, r_0 > 0$ such that every $Y \in \tcalO_{\calE_\Delta(K)}$ can be written uniquely as
	\[ U + (\gamma - 1)(V_i) + \sum_{j \neq i} V_j \]
with $U \in \calO_{\calE_\Delta(K)}$, all $V_j \in T_j$, and
	\[\max \{|U|_r, |V_1|_r, \ldots, |V_n|_r\} \leq c^r |Y|_r\]
for $r \in (0, r_0]$.
\label{decompose-along-gamma}
\end{prop}
\begin{pf}
The case when $n = 1$ was proved in \cite[Corollary 2.5.7]{KedlayaNewMethods}, and an analogous proof establishes the general case.
\end{pf}

\begin{nota}
For $i=0,\dots,n$, let $\calO^{(i)}_{\calE_\Delta(K)}$ be the subring
$\calO_{\calE_\Delta(K)} \oplus T_1 \oplus \cdots \oplus T_{i}$ of $\tcalO_{\calE_\Delta(K)}$;
note that $\calO^{(0)}_{\calE_\Delta(K)} = \calO_{\calE_\Delta(K)}$ 
and $\calO^{(n)}_{\calE_\Delta(K)} = \tcalO_{\calE_\Delta(K)}$.
Also let $\calO^{(i),\dagger}_{\calE_\Delta(K)} := \calO^{(i)}_{\calE_\Delta(K)} \cap \tcalO^\dagger_{\calE_\Delta(K)}$. These rings are stable under $\varphi_\Delta$ and $\Gamma_{K, \Delta}$.
\end{nota}

For these subrings of $\tcalO_{\calE_\Delta(K)}$, we have the following analogue of \thref{tilde-fixed-points-mod-p}.
\begin{lem}
For $i \in \{0,\dots,n\}$,
let $M$ be an \'etale $\phi$-module over $\calO^{(i)}_{\calE_\Delta(K)}$. Then $M^{\phi} = (M \otimes_{\calO^{(i)}_{\calE_\Delta(K)}} \tcalO_{\calE_\Delta(K)})^{\phi}$.
\label{index-fixed-points-mod-p}
\end{lem}
\begin{pf}
This amounts to the full faithfulness of base extension of \'etale $\phi$-modules from $\calO^{(i)}_{\calE_\Delta(K)}$
to $\tcalO_{\calE_\Delta(K)}$. In fact this functor is an equivalence of categories, as may be shown following the proof of \cite[Corollary~5.4.6]{KedlayaLiu2}.
\end{pf}

This in turn yields an analogue of \thref{nu-commutes}.
\begin{lem}
Choose $i \in \{0,\dots,n\}$.
Let $\nu$ be an endomorphism of $\tcalO_{\calE_\Delta(K)}$ which commutes with $\phi$ and sends
$\calO^{(i)}_{\calE_\Delta(K)}$ into itself. 
Let $M$ be an \'etale $\phi$-module over $\tcalO_{\calE_\Delta(K)}$ with a commuting semilinear action of $\nu$.
Suppose that the action of $\phi$ on $M$ is trivial modulo $p$, and that $M$ has a basis $e_1, \ldots, e_d$ with respect to which the matrix $F$ of the action of $\phi$ has entries in $\calO^{(i)}_{\calE_\Delta(K)}$. Then the matrix of the action of $\nu$ with respect to $e_1,\dots,e_d$ also has entries in $\calO^{(i)}_{\calE_\Delta(K)}$.
\label{index-nu-commutes}
\end{lem}
\begin{pf}
Proceed as in the proof of \thref{nu-commutes}, using \thref{index-fixed-points-mod-p} in place of
\thref{tilde-fixed-points-mod-p}.
\end{pf}

At last, we are ready for the crucial calculation, which enables us to eliminate fractional powers one variable at a time.

\begin{lem}
Let $M$ be a projective \'etale $(\phi_\Delta, \Gamma_{K,\Delta})$-module over $\tcalO^\dagger_{\calE_\Delta(K)}$ 
admitting a basis $e_1,\dots,e_d$ fixed modulo $p$ by $\phi_\Delta$ and $\Gamma_{K, \Delta}$. 
For $i = 1, \ldots, n$, pick $\gamma_i \in \Gamma_{\alpha_i}$ with $\gamma_i \equiv 1 \mod p$ and $\gamma_i \neq 1$.
Then for $i=0,\dots,n$, there exists a basis $e_1^{(i)}, \ldots, e_d^{(i)}$ of $M$
congruent to $e_1,\dots,e_d$ modulo $p$,
such that the matrices of the actions of $\phi$ and $\gamma_i$ on this basis have entries in $\calO^{(i),\dagger}_{\calE_\Delta(K)}$.
\label{eliminate-fractional-powers}
\end{lem}
\begin{pf}
We proceed by descending induction on $i$. For the base case $i=n$, we take $e_1^{(n)}, \dots, e_d^{(n)} = e_1,\dots,e_d$. Given the statement for some $i>0$, 
let $G \in \GL_d(\calO^{(i),\dagger}_{\calE_\Delta(K)})$ be the matrix of 
action of $\gamma_i$ on $e_1^{(i)}, \dots, e_d^{(i)}$, i.e.,
	\[ \gamma_i(e_k^{(i)}) = \sum_{j = 1}^d G_{jk} e_j^{(i)}. \]
Since the $e_j^{(i)}$ are fixed by $\gamma_i$ modulo $p$ (by virtue of being congruent to the original $e_j$), we have $p \mid G - 1$. By \thref{abs-val-lim}, there exists an $r \in (0, r_0]$ such that
	\[ \epsilon := \abs[r]{G - 1}^{1/3} < \min\{ c^{-r}, 1 \}. \]

We now construct a sequence of matrices whose product will converge to a change of basis matrix converting $e_1^{(i)},\dots,e_d^{(i)}$ into a new basis $e_1^{(i-1)},\dots,e_d^{(i-1)}$. Let $U_0 := 1$. By \thref{decompose-along-gamma}, there exist matrices
	\[ Y_0 \in M_d(\calO^\dagger_{\calE_\Delta(K)}) \qquad \text{and} \qquad Z_{0j} \in M_d(T_j) \]
for $j = 1, \ldots, i$ such that 
	\[ G = 1 + Y_0 + Z_{01} + \cdots + Z_{0(i-1)} + (\gamma_i - 1)(Z_{0i}) \]
and
	\[ \abs[r]{Y_0}, \abs[r]{Z_{0j}} \leq c^r \abs[r]{G - 1} \leq \epsilon^2 \]
for all $i$. Now suppose that we have constructed matrices $U_l \in \GL_d(\calO^{(i),\dagger}_{\calE_\Delta(K)})$, $Y_l \in M_d(\calO^\dagger_{\calE_\Delta(K)})$, and $Z_{lj} \in M_d(T_j)$ such that
	\begin{gather*}
	U_l \equiv 1 \mod p, \\
	\abs[r]{Y_l}, \abs[r]{Z_{l1}}, \ldots, \abs[r]{Z_{l(i-1)}} \leq \epsilon^2, \\
	\abs[r]{Z_{li}} \leq \epsilon^{l + 2},
	\end{gather*}
and 
	\[ U_l^{-1} G \gamma_i(U_l) = 1 + Y_l + Z_{l1} + \cdots + Z_{l(i-1)} + (\gamma_i - 1)(Z_{li}). \]
Let $G_l := U_l^{-1} G \gamma_1(U_l)$. Let $U_{l + 1} := U_l(1 - Z_{li})$. Then
	\begin{align*}
	G_{l+1} &= (1 - Z_{li})^{-1} U_l G \gamma_i(U_l) \gamma_i(1 - Z_{li}) \\
		&= (1 - Z_{li})^{-1} (1 + Y_l + Z_{l1} + \cdots + Z_{l(i-1)} + (\gamma_i - 1)(Z_{li})(1 - \gamma_i(Z_{li})) \\
		&= 1 + Y_l + Z_{l1} + \cdots + Z_{l(i-1)} + Z_{li} Y_l - Y_l \gamma_i(Z_{li}) + E_l
	\end{align*}
for some $E_l \in M_d(\calO^{(i)}_{\calE_\Delta(K)})$ with $\abs[r]{E_l} \leq \epsilon^{2l + 4}$. We have $Z_{li} Y_l - Y_l \gamma_i(Z_i) \in M_d(\calO^{(i)}_{\calE_\Delta(K)})$ and
	\[ \abs[r]{Z_{li} Y_l - Y_l \gamma_i(Z_{li})} \leq \epsilon^{l + 4}. \]
Write $Z_{li} Y_l - Y_l \gamma_i(Z_{li})$ as
	\[ A_l + B_{l1} + \cdots + B_{l(i-1)} + (\gamma_i - 1)(Z_{li}) \]
with $A_l \in M_d(\calO^\dagger_{\calE_\Delta(K)})$, $B_{lj} \in M_d(T_j)$ for all $j$, and
	\[ \abs[r]{A_l}, \abs[r]{B_{lj}} \leq c^r \epsilon^{l + 4} \leq \epsilon^{l + 3}. \]
Let $Y_{l + 1} := Y_l + A_l$, $Z_{(l + 1)j} := Z_{lj} + B_{lj}$ for $j = 1, \ldots, i-1$, and $Z_{(l + 1)i} := B_{li}$. We then have
	\[ \abs[r]{Y_{l + 1}}, \abs[r]{Z_{(l + 1)1}}, \ldots, \abs[r]{Z_{(l + 1)(i-1)}} \leq \epsilon^2 \qquad \text{and} \qquad \abs[r]{Z_{(l + 1)i}} \leq \epsilon^{(l + 1) + 2}. \]
If $l > 1$, it follows from the fact that $\abs[r]{B_{li}} \leq \epsilon^{l + 3}$ that $U_{l + 1} \equiv 1 \mod p$. If $l = 1$, observe that since $p \mid G - 1$, we have $p \mid (\gamma_i - 1)(Z_{0i})$; it then follows that $p \mid Z_{0i}$ as a consequence of the fact that $\gamma_i - 1$ is bijective on $\bar{T}_i$, so $U_{l + 1} \equiv 1 \mod p$ in this case also.

The product $U_1 U_2 \cdots$ converges to a matrix $U \in \GL_d(\calO^{(i),\dagger}_{\calE_\Delta(K)})$; we thus obtain a new basis $e_1^{(i-1)}, \dots, e_d^{(i-1)}$ of $M$ by setting
	\[ e_k^{(i-1)} := \sum_{j = 1}^d U_{jk} e_j^{(i)} \]
for $k = 1, \ldots, d$. Let $A$ and $H$ be the matrices of action of $\phi$ and $\gamma_i$, respectively, on the $e_j^{(i-1)}$:
	\begin{align*}
	\phi(e_k^{(i-1)}) &= \sum_{j = 1}^d A_{jk} e_j^{(i-1)} \\
	\gamma_i(e_k^{(i-1)}) &= \sum_{j = 1}^d H_{jk} e_j^{(i-1)}.
	\end{align*}
Then $H = U^{-1} G \gamma_i(U) \in M_d(\calO^{(i-1),\dagger}_{\calE_\Delta(K)})$ with $H \equiv 1 \mod p$. Since $\phi$ and $\gamma_i$ commute, we have $A \phi(H) = H \gamma_i(A)$. Write $A = B + C_1 + \cdots + C_i$ with $B \in M_d(\calO^\dagger_{\calE_\Delta(K)})$ and $C_j \in M_d(T_j)$ for all $j$. Then
	\[ H^{-1} C_i \phi(H) - C_i = (\gamma_i - 1)(C_i). \]
If $C_i \neq 0$, then let $m$ be the largest integer such that $p^m \mid C_i$. Since $H \equiv 1 \mod p$, we have $p^{m+1} \mid H^{-1} C_i \phi(H) - C_i$, but, referring again to the fact that $\gamma_i - 1$ is bijective on $\bar{T}_i$, we have $p^{m+1} \nmid (\gamma_i - 1)(C_1)$. By contradiction, $C_i = 0$ and thus the matrix $A$ has entries in
$\calO^{(i),\dagger}_{\calE_\Delta(K)}$.
By \thref{nu-commutes} and \thref{index-nu-commutes},
the matrix of the action of $\gamma_{i-1}$ on the basis of the $e_j^{(i-1)}$ also has entries in $\calO^{(i-1),\dagger}_{\calE_\Delta(K)}$; this completes the induction.
\end{pf}

We finally tie everything together and eliminate our auxiliary hypotheses.
\begin{thm}\label{equivcatoverconv}
Base extension of \'etale $(\phi_\Delta, \Gamma_{K, \Delta})$-modules from $\calO^\dagger_{\calE_{\Delta}(K)}$ to $\calO_{\calE_\Delta(K)}$ is an equivalence of categories.
\end{thm}
\begin{pf}
As in the proof of \thref{tilde-dagger-to-tilde}, there is nothing to check for modules killed by a power of $p$,
and given the statement for projective modules we may easily obtain the general case. We thus consider hereafter
only projective \'etale $(\phi_\Delta, \Gamma_{K, \Delta})$-modules.

By \thref{o-to-o-tilde}, base extension from $\calO_{\calE_\Delta(K)}$ to $\tcalO_{\calE_\Delta(K)}$ is an equivalence; by \thref{tilde-dagger-to-tilde}, base extension from $\tcalO^\dagger_{\calE_{\Delta}(K)}$ to $\tcalO_{\calE_\Delta(K)}$ is an equivalence. It suffices to show that base extension from $\calO^\dagger_{\calE_{\Delta}(K)}$ to $\calO_{\calE_\Delta(K)}$ is fully faithful and base extension from $\calO^\dagger_{\calE_{\Delta}(K)}$ to $\tcalO^\dagger_{\calE_{\Delta}(K)}$ is essentially surjective.
	\begin{center}
	\begin{tikzpicture}
	\node (d) {$\calO^\dagger_{\calE_{\Delta}(K)}$};
	\node[above=of d] (o) {$\calO_{\calE_\Delta(K)}$};
	\node[right=of d] (td) {$\tcalO^\dagger_{\calE_{\Delta}(K)}$};
	\node[above=of td] (t) {$\tcalO_{\calE_\Delta(K)}$};
	\draw[->] (d) -- (o);
	\draw[->] (d) -- (td);
	\draw[->] (o) -- (t);
	\draw[->] (td) -- (t);
	\end{tikzpicture}
	\end{center}
The fact that base extension from $\calO^\dagger_{\calE_{\Delta}(K)}$ to $\calO_{\calE_\Delta(K)}$ is fully faithful follows from \thref{base-ext-mod-p-fully-faithful}, just as in the proof of \thref{tilde-dagger-to-tilde}.

For the essential surjectivity of base extension from $\calO^\dagger_{\calE_{\Delta}(K)}$ to $\tcalO^\dagger_{\calE_{\Delta}(K)}$, as in the proof of \thref{tilde-dagger-to-tilde} it suffices to consider an $\tcalO^\dagger_{\calE_{\Delta}(K)}$-module with a basis fixed modulo $p$ by $\phi_\Delta$ and $\Gamma_{K, \Delta}$: let $M$ be such a module. For $i = 1, \ldots, n$, pick $\gamma_i \in \Gamma_{\alpha_i}$ with $\gamma_i \equiv 1 \mod p$.
By \thref{eliminate-fractional-powers} (applied with $i=0$), there exists a basis $e_1^{(0)}, \ldots, e_r^{(0)}$ of $M$
on which the matrix of action of $\phi$ has entries in $\calO^\dagger_{\calE_{\Delta}(K)}$
and is congruent to the identity matrix modulo $p$.
Let $M^\dagger$ be the $\calO^\dagger_{\calE_\Delta}$-span of $e_1^{(0)}, \ldots, e_r^{(0)}$. This is an \'etale $\phi$-module over $\calO^\dagger_{\calE_\Delta}$ such that $M^\dagger \otimes_{\calO^\dagger_{\calE_\Delta}} \tcalO^\dagger_{\calE_\Delta} \simeq M$. As the matrix of the action of $\phi$ has entries in $\calO^\dagger_{\calE_\Delta(K)}$, 
by \thref{nu-commutes} and \thref{index-nu-commutes} the matrices of each $\phi_\alpha$ and of each $\gamma \in \Gamma_{K, \alpha}$ also have entries in $\calO^\dagger_{\calE_\Delta(K)}$ for each $\alpha \in \Delta$.
Hence $M^\dagger$ is indeed an \'etale $(\varphi_\Delta, \Gamma_{K,\Delta})$-module over $\calO^{\dagger}_{\calE_{\Delta(K)}}$, completing the proof of essential surjectivity.
\end{pf}

For convenience, we summarize everything we have established in a single theorem statement.
\begin{thm}
The category of continuous representations of $G_{K,\Delta}$ on finite free $\Z_p$-modules is equivalent to the category of projective \'etale $(\phi_\Delta, \Gamma_\Delta)$-modules over 
each of the rings $\calO_{\calE_{\Delta}(K)}$, $\tcalO_{\calE_{\Delta}(K)}$, $\calO^\dagger_{\calE_{\Delta}(K)}$, and $\tcalO^\dagger_{\calE_{\Delta}(K)}$.
\label{full-equivalence}
\end{thm}
\begin{pf}
The equivalence between representations and $(\phi_\Delta, \Gamma_\Delta)$-modules over $\calO_{\calE_{\Delta}(K)}$ and
$\tcalO_{\calE_{\Delta}(K)}$ is given by \thref{tilde-equivalence}. We add $\calO^\dagger_{\calE_{\Delta}(K)}$ to the equivalence using \thref{tilde-dagger-to-tilde}, and
$\tcalO^\dagger_{\calE_{\Delta}(K)}$ using 
\thref{equivcatoverconv}.
\end{pf}

\subsection{$\Q_p$-representations}

We next formulate a corresponding statement for $\Q_p$-representations. It is easy to obtain a single statement of the right form, but there are some subtleties around finding the strongest possible formulation; we touch briefly on these issues, but do not give a definitive resolution.

\begin{thm}
The category of continuous representations of $G_{K,\Delta}$ on finite dimensional $\Q_p$-vector spaces is equivalent to the category of projective \'etale $(\phi_\Delta, \Gamma_\Delta)$-modules over 
each of the rings $\calE_{\Delta}(K)$, $\tcalE_{\Delta}(K)$, $\calE_{\Delta}^\dagger(K),$ and $\tcalE_{\Delta}^\dagger(K)$.
\label{full-equivalence-Qp}
\end{thm}
\begin{pf}
Since $G_{K,\Delta}$ is a profinite topological group, it is compact; consequently, any finite dimensional $\Q_p$-vector space with a continuous $G_{K,\Delta}$-action admits a stable $\Z_p$-lattice. (For example, if one starts with any $\Z_p$-lattice, taking the sum of its images under all elements of $G_{K,\Delta}$ yields a stable lattice.) Consequently, \thref{full-equivalence} defines fully faithful functors from the category of continuous representations of $G_{K,\Delta}$ on finite dimensional $\Q_p$-vector spaces to the various categories of projective \'etale $(\phi_\Delta, \Gamma_\Delta)$-modules. 
In view of our definition of the \'etale condition in this setting, essential surjectivity of these functors also reduces at once to \thref{full-equivalence}.
\end{pf}

While \thref{full-equivalence-Qp} will be sufficient for our present purposes, as noted in
\thref{remarks-on-rational-phiGamma-def} it should be possible to formally weaken the definition of a $(\varphi_\Delta, \Gamma_{K,\Delta})$-module over a ring in which $p$ is invertible.
We give some limited evidence in this direction.

\begin{lem}
Let $M$ be an \'etale $\phi_\Delta$-module over $\calO_{\calE_\Delta}/p \calO_{\calE_\Delta}$. Then the underlying module of $M$ is projective over $\calO_{\calE_\Delta}/p \calO_{\calE_\Delta}$.
\label{Mod-p-PhiModule-projective}
\end{lem}
\begin{pf}
Since $\calO_{\calE_\Delta}$ is noetherian, $M$ is finitely presented. Consequently,
$\widetilde{M} := M \otimes_{\calO_{\calE_\Delta}/p\calO_{\calE_\Delta}} R_\Delta$ is a finitely presented \'etale $\varphi_\Delta$-module over 
$R_\Delta$. By \thref{phi-module-flat}, $\widetilde{M}$ is projective over $R_\Delta$. Now note that the morphism $\calO_{\calE_\Delta}/p \calO_{\calE_\Delta} \to R_\Delta$ is split in the category of $\calO_{\calE_\Delta}/p\calO_{\calE_\Delta}$-modules by the morphism
\[
\sum_{i_1,\dots,i_n \in \Z[p^{-1}]} (\bar{a}_{i_1} \bar{\pi}_{\alpha_1}^{i_1}) \otimes \cdots \otimes (\bar{a}_{i_n} \bar{\pi}_{\alpha_n}^{i_n})  \mapsto
\sum_{i_1,\dots,i_n \in \Z} (\bar{a}_{i_1} \bar{\pi}_{\alpha_1}^{i_1}) \otimes \cdots \otimes (\bar{a}_{i_n} \bar{\pi}_{\alpha_n}^{i_n})
\]
(that is, discarding non-integral powers of the series parameters).
This means that $\calO_{\calE_\Delta}/p \calO_{\calE_\Delta} \to R_\Delta$ is a \df{pure} morphism of rings,
so the fact that $\widetilde{M}$ is finite projective over $R_\Delta$
implies that $M$ is finite projective over $\calO_{\calE_\Delta}/p \calO_{\calE_\Delta}$ \cite[Tag~08XD]{StacksProject},
as desired.
\end{pf}

\begin{cor}
Let $M$ be an \'etale $\phi_\Delta$-module over $\calO_{\calE_\Delta}$. 
If the underlying module of $M$ is $p$-torsion-free, then it is projective over $\calO_{\calE_\Delta}$.
\label{TorsionFreePhiModule-projective}
\end{cor}
\begin{pf}
By \thref{Mod-p-PhiModule-projective}, $M/pM$ is projective of some finite rank $r$ over $\calO_{\calE_\Delta}/p \calO_{\calE_\Delta}$. 
We next check that for each positive integer $m$, $M/p^m M$ is projective of rank $r$ over $\calO_{\calE_\Delta}/p^m \calO_{\calE_\Delta}$.
Since $M$ is $p$-torsion-free, it is flat over $\Zp$; consequently, $M/p^m M$ is flat over $\Z/p^m \Z$.
Since $M/pM$ is flat over $\calO_{\calE_\Delta}/p \calO_{\calE_\Delta}$, we may apply \cite[Tag~06A5]{StacksProject} to deduce that $M/p^m M$ is flat over  $\calO_{\calE_\Delta}/p^m \calO_{\calE_\Delta}$; we may then apply \cite[Tag~05CG]{StacksProject} to deduce that $M/p^m M$ is projective over $\calO_{\calE_\Delta}/p^m \calO_{\calE_\Delta}$. The rank-$r$ condition is then enforced by Nakayama's lemma.

Now let $\mathrm{Fitt}_i(M)$ denote the sequence of Fitting ideals of $M$
(see for example \cite[Tag~07Z6]{StacksProject}). 
For each positive integer $m$, the fact that $M/p^m M$ is projective of rank $r$ over $\calO_{\calE_\Delta}/p^m \calO_{\calE_\Delta}$ implies that $\mathrm{Fitt}_i(M)$ is contained in $p^m \calO_{\calE_\Delta}$ for all $i<r$ and is equal to $\calO_{\calE_\Delta}$ for all $i \geq r$.
Running over all $m$, we deduce that $\mathrm{Fitt}_i(M) = 0$ for $i<r$ and is equal to $\calO_{\calE_\Delta}$ for all $i \geq r$; hence $M$ is projective of rank $r$.
\end{pf}

\begin{thm}
Let $M$ be a $(\varphi_\Delta, \Gamma_{K, \Delta})$-module over $\calE_\Delta$
satisfying the following conditions.
\begin{itemize}
\item[(i)]
The underlying $\varphi_\Delta$-module of $M$ is the base extension of an \'etale $\varphi_\Delta$-module over $\calO_{\calE_\Delta}$.
\item[(ii)]
The action of $\Gamma_{K, \Delta}$ is \df{bounded}: for some (and hence any) finitely generated $\calO_{\calE_\Delta}$-submodule $M_0$ of $M$ which generates $M$ over $\calO_{\calE_\Delta}$, the action of $\Gamma_{K, \Delta}$ carries $M_0$ into $p^{-m} M_0$ for some nonnegative integer $m$.
\end{itemize}
Then $M$ is an \'etale $(\varphi_\Delta, \Gamma_{K, \Delta})$-module over $\calE_\Delta$ in the sense of
\thref{phiGamma-rational}; in particular, by \thref{full-equivalence-Qp} it corresponds to a continuous
$\Qp$-representation of $G_{K, \Delta}$.
\label{weaker-phiGamma-rational}
\end{thm}
\begin{pf}
In condition (i), there is no loss of generality in assuming that there exists a finitely generated
$\varphi_\Delta$-stable submodule
$\calO_{\calE_\Delta}$-submodule $M_0$ of $M$ which generates $M$ over $\calO_{\calE_\Delta}$,
such that $M_0$ is an \'etale $\varphi_\Delta$-module over $\calO_{\calE_\Delta}$.
Let $M_1$ be the $\calO_{\calE_\Delta}$-submodule of $M$ generated by $\gamma(M_0)$ for all $\gamma \in \Gamma_{K, \Delta}$; by condition (ii), $M_1$ is again a finitely generated $\calO_{\calE_\Delta}$-module.
In particular, $M_1$ is an \'etale $\varphi_\Delta$-module; by \thref{TorsionFreePhiModule-projective}, $M_1$ is projective over $\calO_{\calE_\Delta}$.
Hence $M_1$ has the structure of a projective \'etale $(\varphi_\Delta, \Gamma_{K, \Delta})$-module over $\calO_{\calE_\Delta}$, and so $M$ is an \'etale $(\varphi_\Delta, \Gamma_{K, \Delta})$-module over $\calE_\Delta$ in the sense of \thref{phiGamma-rational}.
\end{pf}

\begin{rmk}
One may extend \thref{weaker-phiGamma-rational} to $(\varphi_\Delta, \Gamma_{K, \Delta})$-modules over $\calE^\dagger_\Delta$ by similar arguments. We do not know how to extend it to $\tcalE_\Delta$ or $\tcalE^\dagger_\Delta$.

We also do not know whether in condition (i) of \thref{weaker-phiGamma-rational}, one may consider the underlying $\phi$-module instead of the underlying $\phi_\Delta$-module; in particular, it is unclear whether an analogue of the bounded condition must be applied to the partial Frobenius actions. To decide this, one may need some of the slope theory for modules over relative Robba rings developed in \cite{KedlayaLiu}.
\end{rmk}

\section{Galois cohomology}

The goal of this section is to show that the group cohomology of $\GKD$ with values in a $p$-adic representation $V$ is computed by the Herr complex of the multivariate $(\varphi_\Delta,\Gamma_{K,\Delta})$-module associated to $V$
via \thref{full-equivalence} or \thref{full-equivalence-Qp}. The case $K=\Qp$ is proven in \cite{PalZabradi}; we deduce the general case by reducing to that case, using a generalization of Shapiro's Lemma 
for $(\varphi, \Gamma)$-modules \cite[Theorem~2.2]{LiuHerr} to this context.

In the following discussion, we write $\mathbb{D}$ for the functor taking
a $\Z_p$-representation (resp.\ a $\Q_p$-representation) to its corresponding \'etale $(\varphi_\Delta, \Gamma_{K,\Delta})$-module over 
$\calO_{\calE_\Delta(K)}$ (resp.\ over $\calE_{\Delta}(K)$),
and $\mathbb{D}^\dagger$ for the functor taking such a representation to its
corresponding \'etale $(\varphi_\Delta, \Gamma_{K,\Delta})$-module over 
$\calO^\dagger_{\calE_\Delta(K)}$ (resp.\ over $\calE^\dagger_{\Delta}(K)$).

\begin{defn}\label{phicomplexdef}
For any abelian group $D^?$ equipped with commuting operators $\varphi_\alpha$ ($\alpha\in\Delta$) we define the cochain complex
\begin{align*}
\Phi^\bullet(D^?)\colon 0\to D^?\to \bigoplus_{\alpha\in\Delta}D^?\to \dots\to \bigoplus_{\{\alpha_1,\dots,\alpha_r\}\in \binom{\Delta}{r}}D^?\to\dots \to D^?\to 0
\end{align*}
where for all $0\leq r\leq |\Delta|-1$, the map $d_{\alpha_1,\dots,\alpha_r}^{\beta_1,\dots,\beta_{r+1}}\colon D^?\to D^?$ from the component in the $r$th term corresponding to $\{\alpha_1,\dots,\alpha_r\}\subseteq \Delta$ to the component corresponding to the $(r+1)$-tuple $\{\beta_1,\dots,\beta_{r+1}\}\subseteq\Delta$ is given by
\begin{equation*}
d_{\alpha_1,\dots,\alpha_r}^{\beta_1,\dots,\beta_{r+1}}=\begin{cases}0&\text{if }\{\alpha_1,\dots,\alpha_r\}\not\subseteq\{\beta_1,\dots,\beta_{r+1}\}\\ (-1)^{\varepsilon}(\id-\varphi_\beta)&\text{if }\{\beta_1,\dots,\beta_{r+1}\}=\{\alpha_1,\dots,\alpha_r\}\cup\{\beta\}\ ,\end{cases}
\end{equation*}
where $\varepsilon=\varepsilon(\alpha_1,\dots,\alpha_r,\beta)$ is the number of elements in the set $\{\alpha_1,\dots,\alpha_r\}$ smaller than $\beta$. 
\end{defn}

Let $K  / \Qp$ be a finite extension. We denote by $C_{K,\Delta}$ the torsion subgroup of $\Gamma_{K,\Delta}\simeq \prod_{\alpha\in\Delta}\Gal(K(\mu_{p^\infty})  / K)$ and by $\HKD^*$ the kernel of the composite quotient map $\GKD\twoheadrightarrow \Gamma_{K,\Delta}\twoheadrightarrow \Gamma_{K,\Delta}^*:=\Gamma_{K,\Delta}/C_{K,\Delta}$. We choose topological generators $\gamma_{K,\alpha}\in \Gamma_{K,\alpha}^*:=\Gamma_{K,\alpha}/(\Gamma_\alpha\cap C_{K,\Delta})$ for each $\alpha\in\Delta$. 

\begin{defn}\label{gammacomplexdef}
If $A$ is an arbitrary (for now abstract) representation of the group $\Gamma_{K,\Delta}^*\simeq \prod_{\alpha\in\Delta}\Zp$ on a $\Zp$-module we denote by $\Gamma_{K,\Delta}^\bullet(A)$ the cochain complex
\begin{align*}
\Gamma_{K,\Delta}^\bullet(A)\colon 0\to A\to \bigoplus_{\alpha\in\Delta}A\to \dots\to \bigoplus_{\{\alpha_1,\dots,\alpha_r\}\in \binom{\Delta}{r}}A\to\dots \to A\to 0
\end{align*}
where for all $0\leq r\leq |\Delta|-1$, the map $d_{\alpha_1,\dots,\alpha_r}^{\beta_1,\dots,\beta_{r+1}}\colon A\to A$ from the component in the $r$th term corresponding to $\{\alpha_1,\dots,\alpha_r\}\subseteq \Delta$ to the component corresponding to the $(r+1)$-tuple $\{\beta_1,\dots,\beta_{r+1}\}\subseteq\Delta$ is given by
\begin{equation*}
d_{\alpha_1,\dots,\alpha_r}^{\beta_1,\dots,\beta_{r+1}}=\begin{cases}0&\text{if }\{\alpha_1,\dots,\alpha_r\}\not\subseteq\{\beta_1,\dots,\beta_{r+1}\}\\ (-1)^{\varepsilon}(\id-\gamma_{K,\beta})&\text{if }\{\beta_1,\dots,\beta_{r+1}\}=\{\alpha_1,\dots,\alpha_r\}\cup\{\beta\}\ ,\end{cases}
\end{equation*}
where $\varepsilon=\varepsilon(\alpha_1,\dots,\alpha_r,\beta)$ is the number of elements in the set $\{\alpha_1,\dots,\alpha_r\}$ smaller than $\beta$.
\end{defn}

\begin{defn}\label{phigammacomplexdef}
Let $D$ be an \'etale $(\varphi_\Delta,\Gamma_{K,\Delta})$-module over any of the rings $\mathcal{O}_{\mathcal{E}_\Delta(K)}$, $\mathcal{O}_{\mathcal{E}_\Delta(K)}^\dagger$, $\mathcal{E}_\Delta(K)$, or $\mathcal{E}_\Delta^\dagger(K)$. We define the cochain complex $\Phi\Gamma_{K,\Delta}^\bullet(D)$ as the total complex of the double complex $\Gamma_{K,\Delta}^\bullet(\Phi^\bullet(D^{C_{K,\Delta}}))$ and call it the \df{Herr complex} of $D$.
\end{defn}

\begin{defn}
If $\Qp\leq F\leq K$ are finite extensions and $V$ is a continuous finite dimensional representation of $\GKD$ either over $\Qp$ or $\Zp$, then we denote by $\Ind_K^FV:=\Zp[\GFD]\otimes_{\Zp[\GKD]}V$ the representation of $\GFD$ induced from the representation $V$ of the finite index subgroup $\GKD$. 
\end{defn}

\begin{rmk}
Since $\GKD$ has finite index in $\GFD$ the induced and coinduced representations are isomorphic, so we may use the latter.
\end{rmk}

By Shapiro's Lemma for continuous cohomology of profinite groups
\cite[Proposition~I.2.5.10]{SerreGalois}, we have a natural isomorphism of cohomological $\delta$-functors $H^\bullet(\GFD,\Ind_K^F(\cdot))\simeq H^\bullet(\GKD,\cdot)$.

\begin{defn}\label{inducedphigammadef}
Let $D$ be an \'etale $(\varphi_\Delta,\Gamma_{K,\Delta})$-module over any of the rings $\mathcal{O}_{\mathcal{E}_\Delta(K)}$, $\mathcal{E}_\Delta(K)$, $\mathcal{O}_{\mathcal{E}_\Delta(K)}^\dagger$, or $\mathcal{E}^\dagger_\Delta(K)$. We define the induced $(\varphi_\Delta,\Gamma_{F,\Delta})$-module over the analogous ring with base $F$ instead of $K$ as $$\Ind_K^FD:=\Zp[\Gamma_{F,\Delta}]\otimes_{\Zp[\Gamma_{K,\Delta}]}D\ .$$ We equip $\Ind_K^F D$ with the obvious $\Gamma_{F,\Delta}$-action on the left factor, and we put $\varphi_\alpha(\gamma\otimes x):=\gamma\otimes(\varphi_\alpha(x))$ and $\lambda\cdot(\gamma\otimes x):=\gamma\otimes(\gamma^{-1}(\lambda)\cdot x)$ for $\gamma\in\Gamma_{F,\Delta}$, $\alpha\in\Delta$, $\lambda\in \mathcal{O}_{\mathcal{E}_\Delta(F)}$ (resp.\ $\mathcal{E}_\Delta(F)$, $\mathcal{O}_{\mathcal{E}_\Delta(F)}^\dagger$, $\mathcal{E}^\dagger_\Delta(F)$). Here note that $\mathcal{O}_{\mathcal{E}_\Delta(F)}$ (resp.\ $\mathcal{E}_\Delta(F)$, $\mathcal{O}^\dagger_{\mathcal{E}_\Delta(F)}$, $\mathcal{E}^\dagger_\Delta(F)$) is naturally a subring of $\mathcal{O}_{\mathcal{E}_\Delta(K)}$ (resp.\ of $\mathcal{E}_\Delta(K)$, $\mathcal{O}^\dagger_{\mathcal{E}_\Delta(K)}$, $\mathcal{E}^\dagger_\Delta(K)$) and in the definition of $\Ind_K^F D$, we regard $D$ as a module over this subring.
\end{defn}

\begin{prop}\label{inducedphigamma}
Let $\Qp\leq F\leq K$ be finite extensions and let $V$ be a continuous representation of $\GKD$ on a finite dimensional $\Qp$-vector space or finitely generated $\Zp$-module.
Then we have $\mathbb{D}(\Ind_K^FV)\simeq \Ind_K^F\mathbb{D}(V)$ and $\mathbb{D}^\dagger(\Ind_K^FV)\simeq \Ind_K^F\mathbb{D}^\dagger(V)$.
\end{prop}
\begin{pf}
This is completely analogous to the proof of \cite[Proposition~2.1] {LiuHerr}. Choose a set of representatives $\overline{U}$ of cosets of $\Gamma_{F,\Delta}/\Gamma_{K,\Delta}$ and lift it to a subset $U$ of $\GFD$. Similarly, let $W\subset \HFD$ be a set of representatives of the cosets of $\HFD/\HKD$; then $UW=\{uw\mid u\in U,w\in W\}$ is a set of representatives of the cosets of $\GFD/\GKD$. We thus compute 
\begin{align*}
\mathbb{D}(\Ind_K^FV)&=(\Ind_K^FV\otimes_{\Zp}\mathcal{O}_{\mathcal{E}_{\Delta}^{\mathrm{nr}}})^{\HFD}\\
&=\left((\Zp[\GFD]\otimes_{\Zp[\GKD]}V)\otimes_{\Zp}\mathcal{O}_{\mathcal{E}_{\Delta}^{\mathrm{nr}}}\right)^{\HFD}\\
&=\left(\left(\bigoplus_{u\in U}\bigoplus_{w\in W}uw\otimes V\right)\otimes_{\Zp}\mathcal{O}_{\mathcal{E}_{\Delta}^{\mathrm{nr}}}\right)^{\HFD}\\
&=\bigoplus_{u\in U}\left(\left(\bigoplus_{w\in W}uw\otimes V\right)\otimes_{\Zp}\mathcal{O}_{\mathcal{E}_{\Delta}^{\mathrm{nr}}}\right)^{\HFD}\ .
\end{align*}
Now an element $x=\sum_i\sum_{w\in W}uw\otimes v_{i,w}\otimes\lambda_{i,w}\in (\bigoplus_{w\in W}uw\otimes V)\otimes_{\Zp}\mathcal{O}_{\mathcal{E}_{\Delta}^{\mathrm{nr}}}$ lies in the $\HFD$-invariant part if and only if $$u^{-1}x =\sum_i\sum_{w\in W}w\otimes v_{i,w}\otimes u^{-1}(\lambda_{i,w})$$ does, since $\HFD$ is normalized by $u\in\GFD$. Using the invariance under the multiplication by $w'w^{-1}\in\HFD$ we deduce $$x_u:=\sum_i v_{i,w}\otimes w^{-1}u^{-1}(\lambda_{i,w})=\sum_i v_{i,w'}\otimes w'^{-1}u^{-1}(\lambda_{i,w'})$$ for all $w,w'\in W$. Further, $x_u$ must be $\HKD$-invariant, i.e.\ it belongs to $(V\otimes_{\Zp}\mathcal{O}_{\mathcal{E}_{\Delta}^{\mathrm{nr}}})^{\HKD}$. So the isomorphism $\mathbb{D}(\Ind_K^FV)\to \Ind_K^F\mathbb{D}(V)$ is given by sending $x$ to $\sum_{u\in U}\bar{u}\otimes x_u$ where $\bar{u}$ is the image of $u$ under the quotient map $\GFD\to\Gamma_{F,\Delta}$. This is bijective since the ranks of the two sides are equal.

The statement on $\mathbb{D}^\dagger$ follows from \thref{equivcatoverconv} using the result on $\mathbb{D}$.
\end{pf}

\begin{lem}\label{gammaind}
Let $A$ be a $\Zp[\Gamma_{K,\Delta}^*]$-module and assume the topological generators $\gamma_{K,\alpha}$ and $\gamma_{F,\alpha}$ are chosen so that $\gamma_{K,\alpha}=\gamma_{F,\alpha}^{p^r}$ for all $\alpha\in\Delta$ where $p^r:=[\Gamma_{F,\alpha}^*:\Gamma_{K,\alpha}^*]$. Then the complex $\Gamma_{K,\Delta}^\bullet(A)$ is quasi-isomorphic to the complex $\Gamma_{F,\Delta}^\bullet(\Zp[\Gamma_{F,\Delta}^*]\otimes_{\Zp[\Gamma_{K,\Delta}^*]}A)$. The quasi-isomorphism is functorial in $A$.
\end{lem}
\begin{pf}
We proceed by induction on $|\Delta|$. If $|\Delta|=1$ then by our assumption on the topological generators the diagram
\begin{align*}
\xymatrix{
\Gamma_K^\bullet\colon & 0\ar[r] & A\ar[r]^{\gamma_K-\id}\ar[d]_{(\sum_{j=0}^{p^r-1}\gamma_F^j)\otimes\id}& A\ar[d]^{1\otimes\id}\ar[r] & 0\\
\Gamma_F^\bullet\colon & 0\ar[r] & \Zp[\Gamma_{F,\Delta}^*]\otimes_{\Zp[\Gamma_{K,\Delta}^*]}A\ar[r]^{\gamma_F-\id} & \Zp[\Gamma_{F,\Delta}^*]\otimes_{\Zp[\Gamma_{K,\Delta}^*]}A\ar[r] & 0
}
\end{align*}
commutes. Since $(\sum_{j=0}^{p^r-1}\gamma_F^j)\otimes\id$ is injective, so is the induced map on $h^0$. Further, any element $x\in \Zp[\Gamma_{F,\Delta}^*]\otimes_{\Zp[\Gamma_{K,\Delta}^*]}A$ can uniquely be written as $x=\sum_{j=0}^{p^r-1}\gamma_F^j\otimes x_j$ with $x_j\in A$ and $x$ is fixed by $\gamma_F$ if and only if $x_0=x_1=\dots=x_{p^r-1}$ is fixed by $\gamma_F^{p^r}=\gamma_K$. Hence we deduce $h^0\Gamma_K^\bullet(A)\simeq h^0\Gamma_F^\bullet(\Zp[\Gamma_{F,\Delta}^*]\otimes_{\Zp[\Gamma_{K,\Delta}^*]}A)$.

On the other hand if $1\otimes x=(\gamma_F-1)(\sum_{j=0}^{p^r-1}\gamma_F^j\otimes y_j)$ then $y_0=y_1=\dots=y_{p^r-1}$ and $x$ lies in the image of $\gamma_K-1$ since the $\gamma_F^j$-component of the right hand side has to vanish for $j\geq 1$. This shows the injectivity on $h^1$. Finally, $\sum_{j=0}^{p^r-1}\gamma_F^j\otimes x_j-1\otimes(\sum_{j=0}^{p^r-1}x_j)$ lies in the image of $\gamma_F-1$ showing that the induced map on $h^1$ is onto. 

The induction step follows from the spectral sequences associated to the double complexes $\Gamma_{K,\alpha}^\bullet(\Gamma_{K,\Delta\setminus\{\alpha\}}^\bullet(A))$ and $\Gamma_{F,\alpha}^\bullet(\Gamma_{F,\Delta\setminus\{\alpha\}}^\bullet(\Zp[\Gamma_{F,\Delta}^*]\otimes_{\Zp[\Gamma_{K,\Delta}^*]}A))$.
\end{pf}

\begin{rmk}
Whenever the action of $\Zp[\Gamma_{K,\Delta}^*]$ on $A$ extends to the Iwasawa algebra $\Zp\bs \Gamma_{K,\Delta}^*\js$, we may relax our assumption that $\gamma_{K,\alpha}=\gamma_{F,\alpha}^{p^r}$: for any other topological generator $\gamma_{K,\alpha}'$ of the group $\Gamma_{K,\alpha}^*$, the element $\frac{\gamma_{K,\alpha}'-1}{\gamma_{K,\alpha}-1}$ is a unit in the Iwasawa algebra, and therefore the complex defined using $\gamma_{K,\alpha}'$ instead of $\gamma_{K,\alpha}$ ($\alpha\in\Delta$) is quasi-isomorphic to $\Gamma_{K,\Delta}^\bullet(A)$.
\label{independence-of-generator}
\end{rmk}

\begin{thm}\label{shapirophigamma}
We have a natural isomorphism of cohomological $\delta$-functors $$h^i\Phi\Gamma_{F,\Delta}^\bullet(\Ind_K^F(\cdot))\simeq h^i\Phi\Gamma_{K,\Delta}^\bullet(\cdot)$$ on the category of \'etale $(\varphi_\Delta,\Gamma_{K,\Delta})$-modules over $\mathcal{O}_{\mathcal{E}_\Delta(K)}$ (resp.\ over $\mathcal{E}_\Delta(K)$, $\mathcal{O}^\dagger_{\mathcal{E}_\Delta(K)}$, and $\mathcal{E}^\dagger_\Delta(K)$).
\end{thm}
\begin{pf}
Let $D$ be an \'etale $(\varphi_\Delta,\Gamma_{K,\Delta})$-module over any of these rings. Since $\Zp[\Gamma_{F,\Delta}]$ is free as a module over $\Zp[\Gamma_{K,\Delta}]$, we have a natural identification $$h^i\Phi^\bullet(\Ind_K^F D)\simeq \Zp[\Gamma_{F,\Delta}]\otimes_{\Zp[\Gamma_{K,\Delta}]}h^i\Phi^\bullet(D)$$ of cohomological $\delta$-functors. Taking $C_{F,\Delta}$-invariants we obtain $$h^i\Phi^\bullet((\Ind_K^F D)^{C_{F,\Delta}})\simeq \Zp[\Gamma^*_{F,\Delta}]\otimes_{\Zp[\Gamma^*_{K,\Delta}]}h^i\Phi^\bullet(D^{C_{K,\Delta}})\ .$$ In case $\gamma_{K,\alpha}=\gamma_{F,\alpha}^{p^r}$ for all $\alpha\in\Delta$ the result follows from \thref{gammaind} (with $A:=h^i\Phi^\bullet(D^{C_{K,\Delta}})$) using the spectral sequences associated to the double complexes $\Gamma^\bullet_{K,\Delta}(\Phi^\bullet(D^{C_{K,\Delta}}))$ and $\Gamma^\bullet_{F,\Delta}(\Phi^\bullet((\Ind_K^FD)^{C_{F,\Delta}}))$. Note that for $F=\Qp$, the topological generator $\gamma_{\Qp,\alpha}$ of $\Gamma_{\Qp,\alpha}^*$ can be chosen arbitrarily by \cite[Theorem~2.6.2, Corollary~3.5.10]{PalZabradi}. Hence by applying from \thref{independence-of-generator} with $F=\Qp$, we deduce that the Herr complex $\Phi\Gamma_{K,\Delta}^\bullet(D)$ does not depend on the choices of topological generators of $\Gamma_{K,\alpha}^*$ for $\alpha\in\Delta$ up to quasi-isomorphism.
\end{pf}

\begin{cor}
We have a natural isomorphism of cohomological $\delta$-functors $H^i(\GKD,\cdot)\simeq h^i\Phi\Gamma_{K,\Delta}^\bullet(\mathbb{D}(\cdot))\simeq h^i\Phi\Gamma_{K,\Delta}^\bullet(\mathbb{D}^\dagger(\cdot))$ on the categories of continuous $\Zp$- or $\Qp$-representations of $\GKD$.
\label{Galois-cohomology-comparison}
\end{cor}
\begin{pf}
This is a combination of \thref{inducedphigamma}, \thref{shapirophigamma} (with $F:=\Qp$), and  \cite[Theorem~2.6.2]{PalZabradi}. The statement on $\mathbb{D}^\dagger$ follows in a similar fashion from \cite[Corollary~3.5.10]{PalZabradi}.
\end{pf}

Now we turn to the discussion of the Iwasawa cohomology $$H^i_{\mathrm{Iw}}(\GKD,V):=\varprojlim_{\HKD\leq H\leq \GKD}H^i(H,V)\ .$$ By Shapiro's Lemma we have $H^i_{\mathrm{Iw}}(\GKD,V)\simeq H^i(\GKD,\Zp\bs\Gamma_{K,\Delta}\js\otimes_{\Zp}V)$ where the right-hand side refers to continuous cochains via the diagonal action of $\GKD$ on the coefficients (see \cite[Lemma~2.5.1]{PalZabradi}). In particular, $H^i_{\mathrm{Iw}}(\GKD,V)$ is a module over the Iwasawa algebra $\Zp\bs\Gamma_{K,\Delta}\js$. On $\Zp\bs \Gamma_{K,\Delta}\js$-modules, the functor $\Zp\bs \Gamma_{F,\Delta}\js\otimes_{\Zp\bs\Gamma_{K,\Delta}\js}\cdot$ is naturally isomorphic to the functor $\Zp[\Gamma_{F,\Delta}]\otimes_{\Zp[\Gamma_{K,\Delta}]}\cdot$.

\begin{lem}
For any continuous finite dimensional representation $V$ of $\GKD$ over $\Zp$ or $\Qp$, we have $\Zp\bs \Gamma_{F,\Delta}\js\otimes_{\Zp\bs\Gamma_{K,\Delta}\js}H^i_{\mathrm{Iw}}(\GKD,V)\simeq H^i_{\mathrm{Iw}}(\GFD,\Ind_K^FV)$ for all $i\geq 0$.
\end{lem}
\begin{pf}
Let $H\leq \GKD\HFD$ be a subgroup containing $\HFD$ so that we have $\GKD\HFD=\GKD H$. Since the quotient $\GFD/\HFD\simeq \Gamma_{F,\Delta}$ is abelian, $H$ is automatically normal both in $\GFD$ and $\GKD\HFD$. Therefore taking $H$-cohomologies commutes with $\Ind_{\GKD\HFD}^{\GFD}$. In particular, we compute
\begin{align*}
H^i(H,\Ind_K^FV)&\simeq H^i(H,\Zp[\GFD]\otimes_{\Zp[\GKD\HFD]}\Zp[\GKD\HFD]\otimes_{\Zp[\GKD]}V) \\
&\simeq \Zp[\GFD]\otimes_{\Zp[\GKD\HFD]}H^i(H,\Zp[\GKD H]\otimes_{\Zp[\GKD]}V).
\end{align*}
Now any set of representatives of $H/(H\cap \GKD)$ is also a set of representatives of $\GKD H/\GKD$, whence we deduce $$\Zp[\GKD H]\otimes_{\Zp[\GKD]}V\simeq \Zp[H]\otimes_{\Zp[H\cap\GKD]}V.$$ Using Shapiro's Lemma, we obtain
\begin{align*}
H^i(H,\Ind_K^FV)\simeq \Zp[\GFD]\otimes_{\Zp[\GKD\HFD]}H^i(H\cap \GKD,V)\ .
\end{align*}
Taking the projective limit with respect to $H$, we deduce 
\begin{align*}
\Zp[\Gamma_{F,\Delta}]\otimes_{\Zp[\Gamma_{K,\Delta}]}H^i_{\mathrm{Iw}}(\GKD,V)&\simeq \Zp[\GFD]\otimes_{\Zp[\GKD\HFD]}H^i_{\mathrm{Iw}}(\GKD,V) \\
&\simeq H^i_{\mathrm{Iw}}(\GFD,\Ind_K^FV)
\end{align*}
as $\GKD$ acts on $H^i_{\mathrm{Iw}}(\GKD,V)$ via its quotient $\Gamma_{K,\Delta}$ and $\varprojlim_H$ commutes with the functor $\Zp[\GFD]\otimes_{\Zp[\GKD\HFD]}\cdot$ since $\Zp[\GFD]$ is finite free over $\Zp[\GKD\HFD]$.
\end{pf}

Recall \cite{PalZabradi} that $\GQpD$ is a Poincar\'e group at $p$ of dimension $2|\Delta|$. The dualizing module is  $I=\mu_{p^\infty,\Delta}$ which is by definition the $\GQpD$-module isomorphic abstractly to $\mu_{p^\infty}$ (i.e.\ to $\Qp/\Zp$) on which each component $\GQpa$ ($\alpha\in\Delta$) acts as on $\mu_{p^\infty}$ (i.e.\ via the cyclotomic character $\chi_\alpha\colon \GQpa\twoheadrightarrow\Gamma_{\Qp,\alpha}\to\Zp^\times$). Let $\Zp(\mathbf{1}_\Delta):=T_p(\mu_{p^\infty,\Delta})=\varprojlim_n \mu_{p^n,\Delta}$ be the $p$-adic Tate module of $\mu_{p^\infty,\Delta}$. For a $p$-primary discrete $\GQpD$-module $A$, we define the Tate twist $A(\mathbf{1}_\Delta):=A\otimes_{\Zp}\Zp(\TD)$ and the Cartier dual $\Hom(A,\mu_{p^\infty,\Delta})=A^\vee(\TD)$. 

\begin{prop}\label{tate}
For any discrete $p$-primary $\GKD$-module $A$, the cup product pairing induces an isomorphism $H^i(\GKD,A)\simeq H^{2d-i}(\GKD,A^\vee(\TD))^\vee$ for every $i$, where $(\cdot)^\vee=\Hom_{\Zp}(\cdot,\Qp/\Zp)$ stands for the Pontryagin dual.
\end{prop}
\begin{pf}
In case $K=\Qp$ this is \cite[Theorem~2.3.1]{PalZabradi}. For an arbitrary finite extension $K  / \Qp$, the statement follows from Shapiro's Lemma by inducing $A$ from $\GKD$ to $\GQpD$.
\end{pf}

By \thref{tate}, we may further identify these cohomology groups using the Cartier dual $A^\vee(\TD)$ as follows: 
\begin{align*}
H^i(\GKD,\Zp[\GKD/H]\otimes_{\Zp}A)& \simeq
 H^{2d-i}(\GKD,(\Zp[\GKD/H]\otimes_{\Zp}A)^\vee(\TD))^\vee\\
&\simeq H^{2d-i}(\GKD,\Zp[\GKD/H]\otimes_{\Zp}(A^\vee(\TD)))^\vee \\
&\simeq H^{2d-i}(H,A^\vee(\TD))^\vee
\end{align*}
since the index $|\GKD:H|$ is finite. The duals of the corestriction maps are the restriction maps, so we deduce
\begin{equation*}
H^i_{\mathrm{Iw}}(\GKD,A)\simeq \left(\varinjlim_H H^{2d-i}(H,A^\vee(\TD))\right)^\vee=H^{2d-i}(\HKD,A^\vee(\TD))^\vee\ .
\end{equation*}

Moreover, the complex $\Phi^\bullet (\mathbb{D}(A^\vee(\TD)))$ computes the $\HKD$-cohomology of $A^\vee(\TD)$ by \cite[Proposition~2.1.4]{PalZabradi}. (That result is only stated for $K=\Qp$, but the proof---including \cite[Proposition~4.1]{Zabradi}---goes over unchanged to the case of finite extensions $K  / \Qp$.) In particular, this shows
\begin{equation*}
H^i_{\mathrm{Iw}}(\GKD,A)\simeq (h^{2d-i}\Phi^\bullet(\mathbb{D}(A^\vee(\TD)))^\vee)\ .
\end{equation*}

Recall \cite{PalZabradi} that in case $K=\Qp$ and $D_{\Qp}$ is an \'etale $(\varphi_\Delta,\Gamma_\Delta)$-module over $\OED$ killed by a power of $p$, we have the residue pairing
\begin{eqnarray}
\{\cdot,\cdot\}\colon D_{\Qp}\times D_{\Qp}^*(\TD)&\to&\Qp/\Zp\notag\\
(x,y)&\mapsto&\{x,y\}:=\res(y(x)) \ .\label{respairing}
\end{eqnarray}

Here $D_{\Qp}^*:=\Hom_{\OED}(D_{\Qp},\mathcal{E}_\Delta/\OED)$ is the dual $(\varphi_\Delta,\Gamma_\Delta)$-module. Further, $$\res\colon \mathcal{E}_\Delta/\OED(\TD)\to \Qp/\Zp$$ sends an element an element $F(X_\bullet)e$ in $\mathcal{E}_\Delta/\OED(\TD)=\Qp/\Zp\otimes_{\Zp}\OED e$ (with $\varphi_\alpha(e)=e$, $\gamma_\alpha(e)=\chi_\alpha(\gamma_\alpha)e$ where $\chi_\alpha\colon \Gamma_\alpha\to\Zp^\times$ is the cyclotomic character) to the coefficient $a_{-1_\bullet}\in\Qp/\Zp$ of $\frac{1}{X_\Delta}=\prod_{\alpha\in\Delta}X_\alpha^{-1}$ in the expansion of $\frac{F(X_\bullet)}{\prod_{\alpha\in\Delta}(1+X_\alpha)}$ as 
\begin{equation*}
\frac{F(X_\bullet)}{\prod_{\alpha\in\Delta}(1+X_\alpha)}=\sum_{i_\alpha\geq -N_F,\alpha\in\Delta}a_{i_\bullet}\prod_{\alpha\in\Delta}X_\alpha^{i_\alpha}
\end{equation*}
with $a_{i_\bullet}\in\Qp/\Zp$ for $i_\bullet=(i_\alpha)_{\alpha\in\Delta}\in\mathbb{Z}^\Delta$ and some integer $N_F\in \mathbb{Z}$ depending on $F$. Moreover, \eqref{respairing} is $\Gamma_\Delta$- and $\prod_{\alpha\in\Delta}(1+X_\alpha)^{\Zp}$-equivariant with respect to which the adjoint of $\varphi_\alpha$ is $\psi_\alpha$ ($\alpha\in\Delta$).

Now if $A$ is a continuous mod-$p^n$ representation of $\GKD$, then the residue pairing \eqref{respairing} applies to $D_{\Qp}:=\Ind_K^{\Qp}\mathbb{D}(A)$.
\begin{lem}\label{gammaorth}
We have $D^*_{\Qp}\simeq \Zp[\Gamma_{\Qp,\Delta}]\otimes_{\Zp[\Gamma_{K,\Delta}]}\mathbb{D}(A)^*$ where we put $\mathbb{D}(A)^*:=\Hom_{\mathcal{O}_{\mathcal{E}_\Delta(K)}}(\mathbb{D}(A),\mathcal{E}_\Delta(K)/\mathcal{O}_{\mathcal{E}_\Delta(K)})$. Moreover, under this identification $1\otimes \mathbb{D}(A)$ is orthogonal to $\gamma\otimes \mathbb{D}(A)^*(\TD)$ for all $\gamma\in \Gamma_{\Qp,\Delta}$ not lying in $\Gamma_{K,\Delta}$. In particular the residue pairing \eqref{respairing} descends to a pairing
\begin{align*}
\{\cdot,\cdot\}\colon \mathbb{D}(A)\left(\simeq 1\otimes \mathbb{D}(A)\right)\quad\times\quad \left(1\otimes \mathbb{D}(A)^*(\TD)\simeq\right) \mathbb{D}(A)^*(\TD)\to\Qp/\Zp
\end{align*}
such that $\{\sum_{\gamma\in U}\gamma\otimes x_\gamma,\sum_{\gamma\in U}\gamma\otimes y_\gamma\}=\sum_{\gamma\in U}\{x_\gamma,y_\gamma\}$ for any choice $U\subset\Gamma_{\Qp,\Delta}$ of coset representatives of $\Gamma_{\Qp,\Delta}/\Gamma_{K,\Delta}$.
\end{lem}
\begin{pf}
We define the map $F$ as
\begin{eqnarray*}
\Zp[\Gamma_{\Qp,\Delta}]\otimes_{\Zp[\Gamma_{K,\Delta}]}\mathbb{D}(A)^*&\to& D^*_{\Qp}\\
\sum_{\gamma\in U}\gamma\otimes f_\gamma &\mapsto& \left(\sum_{\gamma\in U}\gamma\otimes x_\gamma\mapsto \sum_{\gamma\in U}\Tr_{\HQpD/\HKD}(f_\gamma(x_\gamma))\right)
\end{eqnarray*}
where $\Tr_{\HQpD/\HKD}=\sum_{u\in \HQpD/\HKD}u\colon \mathcal{E}_\Delta(K)/\mathcal{O}_{\mathcal{E}_\Delta(K)}\to \mathcal{E}_\Delta/\mathcal{O}_{\mathcal{E}_\Delta}$ is the trace map. The map $F$ is $\OED$-linear and $(\varphi_\Delta,\Gamma_{\Qp,\Delta})$-equivariant by construction. For the bijectivity of $F$ assume first that $pA=0$. Note that $$\Tr_{\HQpD/\HKD}\colon E_{\Delta}(K)\simeq p^{-1}\mathcal{O}_{\mathcal{E}_\Delta(K)}/\mathcal{O}_{\mathcal{E}_\Delta(K)} \to E_{\Delta}(\Qp)\simeq p^{-1}\OED/\OED$$ is onto since it is the composite of the trace maps $\Tr_{\HQpa/\HKa}\colon E_{\alpha}(K)\to E_\alpha(\Qp)$ for all $\alpha\in\Delta$ that are each onto since $E_\alpha(K)  / E_\alpha(\Qp)$ is a Galois extension with Galois group $\HQpa/\HKa$. In particular, $\Ker(\Tr_{\HQpa/\HKa})$ does not contain any nonzero ideal of $E_\Delta(K)$ as its rank over $E_\Delta(\Qp)$ equals $|\HQpD/\HKD|-1<|\HQpD/\HKD|$. Now if $\sum_{\gamma\in U}\gamma\otimes f_\gamma\neq 0$ then $f_\gamma\neq 0$ for at least one choice of $\gamma\in U$, so we may put $x_{\gamma'}:=0$ for all $\gamma\neq \gamma'\in U$ and choose $x_\gamma$ so that $\Tr_{\HQpD/\HKD}(f_\gamma(x_\gamma))=\Tr_{\HQpD/\HKD}(x_\gamma f_\gamma(1))\neq 0$ (as we have just seen that $\Tr_{\HQpD/\HKD}(E_\Delta(K) f_\gamma(1))\neq \{0\}$). Hence $F$ is injective if $pA=0$. On the other hand, the ranks of the domain and codomain of $F$ are equal, therefore it is an isomorphism by \cite[Proposition~2.2]{Zabradi}. The case of general $A$ follows by devissage.

The second statement is deduced from the first one using the $\Gamma_{\Qp,\Delta}$-invariance of the pairing $\{\cdot,\cdot\}$.
\end{pf}

\begin{defn}\label{psicomplexdefinition}
Let $D$ be an \'etale $(\varphi_\Delta,\Gamma_{K,\Delta})$-module over any of the rings $\mathcal{O}_{\mathcal{E}_\Delta(K)}$, $\mathcal{E}_\Delta(K)$, $\mathcal{O}^\dagger_{\mathcal{E}_\Delta(K)}$, or $\mathcal{E}^\dagger_\Delta(K)$. We define the cochain complex
\begin{align}\label{psicomplexdef}
\Psi^\bullet(D)\colon 0\to D\to \bigoplus_{\alpha\in\Delta}D\to \dots\to \bigoplus_{\{\alpha_1,\dots,\alpha_r\}\in \binom{\Delta}{r}}D\to\dots \to D\to 0
\end{align}
where for all $0\leq r\leq |\Delta|-1$ the map $d_{\alpha_1,\dots,\alpha_r}^{\beta_1,\dots,\beta_{r+1}}\colon D\to D$ from the component in the $r$th term corresponding to $\{\alpha_1,\dots,\alpha_r\}\subseteq \Delta$ to the component corresponding to the $(r+1)$-tuple $\{\beta_1,\dots,\beta_{r+1}\}\subseteq\Delta$ is given by
\begin{equation*}
d_{\alpha_1,\dots,\alpha_r}^{\beta_1,\dots,\beta_{r+1}}=\begin{cases}0&\text{if }\{\alpha_1,\dots,\alpha_r\}\not\subseteq\{\beta_1,\dots,\beta_{r+1}\}\\ (-1)^{\eta}(\id-\psi_\beta)&\text{if }\{\beta_1,\dots,\beta_{r+1}\}=\{\alpha_1,\dots,\alpha_r\}\cup\{\beta\}\ ,\end{cases}
\end{equation*}
where $\eta=\eta(\alpha_1,\dots,\alpha_r,\beta)$ is the number of elements in the set $\Delta\setminus\{\alpha_1,\dots,\alpha_r\}$ smaller than $\beta$, and $\psi_\beta$ is the reduced trace of $\phi_\beta$, ie.\ for an element $x=\sum_{j=0}^{p-1}(1+X_\beta)^j\varphi_\beta(x_j)\in D$ we put $\psi_\beta(x)=x_0$. (Note that the sign convention here is different from the one defining the complex $\Phi^\bullet(D)$. The reason for this is that with this choice of signs, the differentials are adjoint to each other under the residue pairing \eqref{respairing} defined above.)
\end{defn}

Since the $\psi$-operators commute with the action of $\Gamma_{K,\Delta}$, we have $\Psi^\bullet(\Ind_K^F D)\simeq \Zp[\Gamma_{F,\Delta}]\otimes_{\Zp[\Gamma_{K,\Delta}]}\Psi^\bullet(D)$. Further, $\Zp[\Gamma_{F,\Delta}]$ is finite free over $\Zp[\Gamma_{K,\Delta}]$, so we deduce $h^i\Psi^\bullet(\Ind_K^F D)\simeq \Zp[\Gamma_{F,\Delta}]\otimes_{\Zp[\Gamma_{K,\Delta}]}h^i\Psi^\bullet(D)$ for all $i\geq 0$. Similarly, we have $\Phi^\bullet(\Ind_K^F D)\simeq \Zp[\Gamma_{F,\Delta}]\otimes_{\Zp[\Gamma_{K,\Delta}]}\Phi^\bullet(D)$ and $h^i\Phi^\bullet(\Ind_K^F D)\simeq \Zp[\Gamma_{F,\Delta}]\otimes_{\Zp[\Gamma_{K,\Delta}]}h^i\Phi^\bullet(D)$ for all $i\geq 0$.

\begin{thm}
We have a natural isomorphism of cohomological $\delta$-functors $H^i_{\mathrm{Iw}}(\GKD,\cdot)\simeq h^{i-d}\Psi^\bullet(\mathbb{D}(\cdot))\simeq h^{i-d}\Psi^\bullet(\mathbb{D}^\dagger(\cdot))$ on the categories of continuous $\Zp$- or $\Qp$-representations of $\GKD$.
\label{comparison-Iwasawa-cohomology}
\end{thm}
\begin{pf}
The case $K=\Qp$ is proven in \cite[Corollary~3.5.10]{PalZabradi}. Let $A$ be a $p$-power torsion representation of $\GKD$ and put $D:=\mathbb{D}(A)$, $D_{\Qp}:=\Ind_K^{\Qp}D$. The isomorphism is constructed via the residue pairing which induces a pairing between $\Psi^{d-r}(D_{\Qp})$ and $\Phi^r(D_{\Qp}^*(\TD))$, and hence between $h^{d-r}\Psi^\bullet(D_{\Qp})$ and $h^r\Phi^\bullet(D_{\Qp}^*(\TD))$ ($0\leq r\leq d$). By \thref{gammaorth}, the isomorphism $$\eta_r(A)\colon h^{d-r}\Psi^\bullet(D_{\Qp})\to h^r\Phi^\bullet(D_{\Qp}^*(\TD))^\vee\simeq H^r_{\mathrm{Iw}}(\GQpD,\Ind_K^{\Qp}A)$$ descends to a $\Zp\bs\Gamma_{K,\Delta}\js$-linear map
\begin{align*}
\widetilde{\eta_r}(A)\colon h^{d-r}\Psi^\bullet(D)\to h^r\Phi^\bullet(D^*(\TD))^\vee\simeq H^r_{\mathrm{Iw}}(\GKD,A)
\end{align*}
such that $\eta_r=\Zp\bs\Gamma_{\Qp,\Delta}\js\otimes_{\Zp\bs\Gamma_{K,\Delta}\js}\widetilde{\eta_r}$. Since $\eta_r$ is an isomorphism and the ring extension $\Zp\bs\Gamma_{K,\Delta}\js\hookrightarrow \Zp\bs\Gamma_{\Qp,\Delta}\js$ is faithfully flat, we deduce that $\widetilde{\eta_r}$ is also an isomorphism. If $T$ is an arbitrary continuous representation of $\GKD$ on a finitely generated $\Zp$-module, then the result on $\mathbb{D}(T)$ follows by taking the projective limit of the isomorphisms $\widetilde{\eta_r}(T/p^nT)$.

Now if $V$ is a continuous representation of $\GKD$ over $\Qp$, then it contains a $\GKD$-invariant $\Zp$-lattice $T\leq V$ by the compactness of $\GKD$ and the isomorphism for $\mathbb{D}(V)$ follows from that for $\mathbb{D}(T)$ by inverting $p$. Finally, the inclusion of $\mathbb{D}^\dagger(T)\hookrightarrow \mathbb{D}(T)$ (resp.\ $\mathbb{D}^\dagger(V)\hookrightarrow \mathbb{D}(V)$) induces a morphism $\Psi^\bullet(\mathbb{D}^\dagger(T))\to \Psi^\bullet(\mathbb{D}(T))$ (resp.\ $\Psi^\bullet(\mathbb{D}^\dagger(V))\to \Psi^\bullet(\mathbb{D}(V))$) and, by taking cohomologies, a $\Zp\bs\Gamma_{K,\Delta}\js$-linear map $\iota_i\colon h^i\Psi^\bullet(\mathbb{D}^\dagger(T))\to h^i\Psi^\bullet(\mathbb{D}(T))$ (resp.\ $\iota_i\colon h^i\Psi^\bullet(\mathbb{D}^\dagger(V))\to h^i\Psi^\bullet(\mathbb{D}(V))$) for $i\geq 0$ which becomes an isomorphism after base change to $\Zp\bs \Gamma_{\Qp,\Delta}\js$ by \cite[Corollary~3.5.10]{PalZabradi}. Hence $\iota_i$ is an isomorphism for all $i\geq 0$ using again the faithfully flat property of the ring extension $\Zp\bs\Gamma_{K,\Delta}\js\hookrightarrow \Zp\bs \Gamma_{\Qp,\Delta}\js$.
\end{pf}

\section{Distinct factors}
\label{distinct-factors}

We now formulate the corresponding results for products of Galois groups of distinct finite extensions of $\Q_p$.
\begin{nota}
Again let $\Delta$ be a finite set, but now let $\underline{K} = (K_\alpha: \alpha \in \Delta)$ be a tuple of finite extensions of $\Qp$ and put
\begin{align*}
G_{\underline{K}, \Delta} &:= \prod_{\alpha \in \Delta} \Gal(K_\alpha\alg  / K_\alpha) \\
H_{\underline{K}, \Delta} &:= \prod_{\alpha \in \Delta} \Gal(K_\alpha\alg  / K_\alpha(\mu_{p^\infty})) \\
\Gamma_{\underline{K}, \Delta} &:= \prod_{\alpha \in \Delta} \Gal(K_\alpha(\mu_{p^\infty})  / K_\alpha).
\end{align*}
We then follow \thref{notation-for-product-rings-1},
\thref{notation-for-Gauss-norms}, and \thref{notation-for-product-rings-2}, but taking $\tcalO_{\calE_\alpha}$ to be a copy of $\tcalO_{\calE}$ for $\calE$ associated to the field $K_\alpha$. This yields rings
\begin{gather*}
\calO_{\calE_\Delta(\underline{K})}, \quad \tcalO_{\calE_\Delta(\underline{K})}, \quad
\calO^\dagger_{\calE_\Delta(\underline{K})}, \quad \tcalO^\dagger_{\calE_\Delta(\underline{K})}, \\
\calE_\Delta(\underline{K}), \quad \tcalE_\Delta(\underline{K}), \quad
\calE^\dagger_{\Delta}(\underline{K}), \quad \tcalE^\dagger_\Delta(\underline{K}).
\end{gather*}
\end{nota}

We then have the following extensions of \thref{full-equivalence} and \thref{full-equivalence-Qp}.

\begin{thm}
The category of continuous representations of $G_{\underline{K},\Delta}$ on finite free $\Z_p$-modules is equivalent to the category of projective \'etale $(\phi_\Delta, \Gamma_\Delta)$-modules over 
each of the rings $\calO_{\calE_{\Delta}(\underline{K})}$, $\tcalO_{\calE_{\Delta}(\underline{K})}$, $\calO^\dagger_{\calE_{\Delta}(\underline{K})}$, and $\tcalO^\dagger_{\calE_{\Delta}(\underline{K})}$.
\label{full-equivalence-multiple}
\end{thm}
\begin{pf}
Since \thref{cor-4.3.16-integral} already includes the case where the fields $F_i$ need not be equal, we may apply it to deduce the analogue of \thref{tilde-equivalence}; this gives the equivalence between continuous $\Zp$-representations of $G_{\underline{K},\Delta}$ and $(\phi_\Delta, \Gamma_\Delta)$-modules over $\tcalO_{\calE_\Delta(\underline{K})}$.
To add the other rings, it suffices to do so after replacing each $K_\alpha$ with a single mutual extension $K$;
we thus reduce to \thref{full-equivalence}.
\end{pf}

\begin{thm}
The category of continuous representations of $G_{\underline{K},\Delta}$ on finite dimensional $\Q_p$-vector spaces is equivalent to the category of projective \'etale $(\phi_\Delta, \Gamma_\Delta)$-modules over 
each of the rings $\calE_{\Delta}(\underline{K})$, $\tcalE_{\Delta}(\underline{K})$, $\calE_{\Delta}^\dagger(\underline{K})$, and $\tcalE_{\Delta}^\dagger(\underline{K})$.
\label{full-equivalence-Qp-multiple}
\end{thm}
\begin{pf}
This follows from \thref{full-equivalence-multiple} just as \thref{full-equivalence-Qp} follows from \thref{full-equivalence}.
\end{pf}

In order to extend the results on Galois cohomology to distinct finite extensions of $\Qp$, note that \thref{phicomplexdef} of $\Phi^\bullet(\cdot)$ carries over to this situation. Further, putting $C_{\underline{K},\Delta}$ for the torsion subgroup of $\Gamma_{\underline{K},\Delta}$ and choosing generators $\gamma_\alpha\in \Gamma_{K_\alpha}/(\Gamma_{K_\alpha}\cap C_{\underline{K},\Delta})$ (where $\Gamma_{K_\alpha}:=\Gal(K_\alpha(\mu_{p^\infty})  / K_\alpha)$) we may also form the complex $\Phi\Gamma^\bullet_{\underline{K},\Delta}(D)$ as in \thref{gammacomplexdef} and \thref{phigammacomplexdef} for a projective \'etale $(\varphi_\Delta,\Gamma_\Delta)$-module $D$ over any of the rings $\calO_{\calE_{\Delta}(\underline{K})}$, $\calO^\dagger_{\calE_{\Delta}(\underline{K})}$, $\calE_{\Delta}(\underline{K})$, or $\calE_{\Delta}^\dagger(\underline{K})$. Moreover, if $\underline{F} = (F_\alpha: \alpha \in \Delta)$ is another tuple of finite extensions of $\Qp$ satisfying $F_\alpha\leq K_\alpha$ for all $\alpha\in\Delta$ then the induced $(\varphi_\Delta,\Gamma_\Delta)$-module $\Ind_{\underline{K}}^{\underline{F}}D := \Zp[\Gamma_{\underline{F},\Delta}]\otimes_{\Zp[\Gamma_{\underline{K},\Delta}]}D$ (as in \thref{inducedphigammadef}) also makes sense. So we have the following version of Shapiro's lemma.

\begin{thm}
We have a natural isomorphism of cohomological $\delta$-functors $$h^i\Phi\Gamma_{\underline{F},\Delta}^\bullet(\Ind_{\underline{K}}^{\underline{F}}(\cdot))\simeq h^i\Phi\Gamma_{\underline{K},\Delta}^\bullet(\cdot)$$ on the category of projective \'etale $(\varphi_\Delta,\Gamma_{\Delta})$-modules over $\mathcal{O}_{\mathcal{E}_\Delta(\underline{K})}$ (resp.\ over $\mathcal{E}_\Delta(\underline{K})$, $\mathcal{O}^\dagger_{\mathcal{E}_\Delta(\underline{K})}$, and $\mathcal{E}^\dagger_\Delta(\underline{K})$).
\end{thm}
\begin{pf}
The proof of \thref{shapirophigamma} goes through unchanged.
\end{pf}

Applying this in case $F_\alpha=\Qp$ for all $\alpha\in\Delta$ we deduce:
\begin{cor}
We have a natural isomorphism of cohomological $\delta$-functors $H^i(G_{\underline{K},\Delta},\cdot)\simeq h^i\Phi\Gamma_{\underline{K},\Delta}^\bullet(\mathbb{D}(\cdot))\simeq h^i\Phi\Gamma_{\underline{K},\Delta}^\bullet(\mathbb{D}^\dagger(\cdot))$ on the categories of continuous $\Zp$- or $\Qp$-representations of $G_{\underline{K},\Delta}$.
\end{cor}

Finally, \thref{psicomplexdefinition} of $\Psi^\bullet(\cdot)$ also carries over to this case and we have:
\begin{thm}
We have a natural isomorphism of cohomological $\delta$-functors $H^i_{\mathrm{Iw}}(G_{\underline{K},\Delta},\cdot)\simeq h^{i-d}\Psi^\bullet(\mathbb{D}(\cdot))\simeq h^{i-d}\Psi^\bullet(\mathbb{D}^\dagger(\cdot))$ on the categories of continuous $\Zp$- or $\Qp$-representations of $G_{\underline{K},\Delta}$.
\end{thm}
\begin{pf}
The proof of \thref{comparison-Iwasawa-cohomology} goes through unchanged.
\end{pf}

\section{Future directions}

We end by recording some possible directions in which this work could be continued.
\begin{itemize}
\item
Extend the comparison of Galois cohomology and the Herr complex to
the rings $\tcalO_{\calE_\Delta}$, $\tcalO^\dagger_{\calE_\Delta}$, $\tcalE_\Delta$, and $\tcalE^\dagger_\Delta$, for which the maps $\phi_\alpha$ are bijective and so the construction of the reduced trace $\psi_\alpha$ is not relevant.
\item
Construct the analogue of the Robba ring and relate $(\phi_\Delta, \Gamma_\Delta)$-modules over it to vector bundles on a product of Fargues--Fontaine curves.
\item
Extend to representations with values in a coherent sheaf on a rigid analytic space over $\Q_p$, as in \cite{KedlayaPottharstXiao}.
\item
Consider Iwasawa cohomology in terms of cyclotomic deformations from the point of view
of \cite{KedlayaPottharst}. This potentially allows for the use of towers other than the cyclotomic one.
\item
Extend to representations of a product of \'etale fundamental groups of rigid analytic spaces, as in \cite{KedlayaLiu, KedlayaLiu2}. Note that \thref{o-plus-of-Y} is already written at a suitable level of generality for this purpose.
\item
Apply Drinfeld's lemma for perfectoid spaces to other constructions of multivariate $(\varphi, \Gamma)$-modules, such as that of Berger \cite{BergerMulti1,BergerMulti2}. In that construction, one starts with a representation of a single copy of $G_K$, but it should be possible to interpret the resulting objects in terms of representations of a suitable power of $G_K$.
\end{itemize}

\bibliographystyle{plain}
\bibliography{RefList}

\def\cprime{$'$}
\begin{thebibliography}{10}

\bibitem{Andre2}
Y.~Andr{\'e}.
\newblock La conjecture du facteur direct.
\newblock {\em Publ. Math. IH\'ES}, 127:71--93, 2018.

\bibitem{Andre1}
Y.~Andr{\'e}.
\newblock Le lemme d'{A}bhyankar perfectoide.
\newblock {\em Publ. Math. IH\'ES}, 127:1--70, 2018.

\bibitem{Berger1}
Laurent Berger.
\newblock Repr\'esentations {$p$}-adiques et \'equations diff\'erentielles.
\newblock {\em Invent. Math.}, 148(2):219--284, 2002.

\bibitem{BergerIntro}
Laurent Berger.
\newblock An introduction to the theory of {$p$}-adic representations.
\newblock In {\em Geometric aspects of {D}work theory. {V}ol. {I}, {II}}, pages
  255--292. Walter de Gruyter, Berlin, 2004.

\bibitem{Berger2}
Laurent Berger.
\newblock \'{E}quations diff\'erentielles {$p$}-adiques et {$(\phi,N)$}-modules
  filtr\'es.
\newblock {\em Ast\'erisque}, (319):13--38, 2008.
\newblock Repr\'esentations $p$-adiques de groupes $p$-adiques. I.
  Repr\'esentations galoisiennes et $(\phi,\Gamma)$-modules.

\bibitem{BergerMulti1}
Laurent Berger.
\newblock Multivariable {L}ubin-{T}ate {$(\varphi,\Gamma)$}-modules and
  filtered {$\varphi$}-modules.
\newblock {\em Math. Res. Lett.}, 20(3):409--428, 2013.

\bibitem{BergerMulti2}
Laurent Berger.
\newblock Multivariable {$(\varphi,\Gamma)$}-modules and locally analytic
  vectors.
\newblock {\em Duke Math. J.}, 165(18):3567--3595, 2016.

\bibitem{Bhatt}
Bhargav Bhatt.
\newblock On the direct summand conjecture and its derived variant.
\newblock {\em Invent. Math.}, 212(2):297--317, 2018.

\bibitem{BhattMorrowScholze}
Bhargav Bhatt, Matthew Morrow, and Peter Scholze.
\newblock Integral {$p$}-adic {H}odge theory.
\newblock {\em Publ. Math. Inst. Hautes \'{E}tudes Sci.}, 128:219--397, 2018.

\bibitem{CherbonnierColmez}
F.~Cherbonnier and P.~Colmez.
\newblock Repr\'esentations {$p$}-adiques surconvergentes.
\newblock {\em Invent. Math.}, 133(3):581--611, 1998.

\bibitem{Colmez}
Pierre Colmez.
\newblock Repr\'esentations de {${\rm GL}_2(\mathbf{Q}_p)$} et
  {$(\phi,\Gamma)$}-modules.
\newblock {\em Ast\'erisque}, (330):281--509, 2010.

\bibitem{deJong}
A.~J. de~Jong.
\newblock \'{E}tale fundamental groups of non-{A}rchimedean analytic spaces.
\newblock {\em Compositio Math.}, 97(1-2):89--118, 1995.
\newblock Special issue in honour of Frans Oort.

\bibitem{DrinfeldICM}
V.~G. Drinfel{\cprime}d.
\newblock Langlands' conjecture for {${\rm GL}(2)$}\ over functional fields.
\newblock In {\em Proceedings of the {I}nternational {C}ongress of
  {M}athematicians ({H}elsinki, 1978)}, pages 565--574. Acad. Sci. Fennica,
  Helsinki, 1980.

\bibitem{DrinfeldCohomology}
V.~G. Drinfel{\cprime}d.
\newblock Cohomology of compactified moduli varieties of {$F$}-sheaves of rank
  {$2$}.
\newblock {\em Zap. Nauchn. Sem. Leningrad. Otdel. Mat. Inst. Steklov. (LOMI)},
  162(Avtomorfn. Funkts. i Teor. Chisel. III):107--158, 189, 1987.

\bibitem{FarguesFontaine}
Laurent Fargues and Jean-Marc Fontaine.
\newblock Courbes et fibr\'{e}s vectoriels en th\'{e}orie de {H}odge
  {$p$}-adique.
\newblock {\em Ast\'{e}risque}, (406):xiii+382, 2018.
\newblock With a preface by Pierre Colmez.

\bibitem{Fontaine}
Jean-Marc Fontaine.
\newblock Repr\'esentations {$p$}-adiques des corps locaux. {I}.
\newblock In {\em The {G}rothendieck {F}estschrift, {V}ol.\ {II}}, volume~87 of
  {\em Progr. Math.}, pages 249--309. Birkh\"auser Boston, Boston, MA, 1990.

\bibitem{FontaineWintenberger2}
Jean-Marc Fontaine and Jean-Pierre Wintenberger.
\newblock Extensions alg\'ebrique et corps des normes des extensions {APF} des
  corps locaux.
\newblock {\em C. R. Acad. Sci. Paris S\'er. A-B}, 288(8):A441--A444, 1979.

\bibitem{FontaineWintenberger1}
Jean-Marc Fontaine and Jean-Pierre Wintenberger.
\newblock Le ``corps des normes''\ de certaines extensions alg\'ebriques de
  corps locaux.
\newblock {\em C. R. Acad. Sci. Paris S\'er. A-B}, 288(6):A367--A370, 1979.

\bibitem{Herr1}
Laurent Herr.
\newblock Sur la cohomologie galoisienne des corps {$p$}-adiques.
\newblock {\em Bull. Soc. Math. France}, 126(4):563--600, 1998.

\bibitem{Herr2}
Laurent Herr.
\newblock Une approche nouvelle de la dualit\'e locale de {T}ate.
\newblock {\em Math. Ann.}, 320(2):307--337, 2001.

\bibitem{Huber}
Roland Huber.
\newblock {\em \'Etale cohomology of rigid analytic varieties and adic spaces}.
\newblock Aspects of Mathematics, E30. Friedr. Vieweg \& Sohn, Braunschweig,
  1996.

\bibitem{kaplansky1}
Irving Kaplansky.
\newblock Maximal fields with valuations.
\newblock {\em Duke Math. J.}, 9:303--321, 1942.

\bibitem{KedlayaNewMethods}
Kiran~S. Kedlaya.
\newblock New methods for {$(\varphi,\Gamma)$}-modules.
\newblock {\em Res. Math. Sci.}, 2:Art. 20, 31, 2015.

\bibitem{KedlayaSSS}
Kiran~S. Kedlaya.
\newblock Sheaves, stacks, and shtukas.
\newblock In {\em Perfectoid Spaces: Lectures from the 2017 Arizona Winter
  School}, pages 58--205. Amer. Math. Soc., 2019.

\bibitem{KedlayaFrobMod}
Kiran~S. Kedlaya.
\newblock Frobenius modules over multivariate {R}obba rings.
\newblock arXiv:1311.7468v3, 2020.

\bibitem{KedlayaSimple}
Kiran~S. Kedlaya.
\newblock Simple connectivity of {Fargues}-{Fontaine} curves.
\newblock 2021.
\newblock arXiv:1806.11528v4.

\bibitem{KedlayaLiu}
Kiran~S. Kedlaya and Ruochuan Liu.
\newblock Relative {$p$}-adic {H}odge theory: {F}oundations.
\newblock {\em Ast\'erisque}, (371):239, 2015.

\bibitem{KedlayaLiu2}
Kiran~S. Kedlaya and Ruochuan Liu.
\newblock Relative {$p$}-adic {H}odge theory, {II}: Imperfect period rings.
\newblock arXiv:1602.06899v3, 2019.

\bibitem{KedlayaPottharst}
Kiran~S. Kedlaya and Jonathan Pottharst.
\newblock On categories of {$(\varphi,\Gamma)$}-modules.
\newblock In {\em Algebraic geometry: {S}alt {L}ake {C}ity 2015}, volume~97 of
  {\em Proc. Sympos. Pure Math.}, pages 281--304. Amer. Math. Soc., Providence,
  RI, 2018.

\bibitem{KedlayaPottharstXiao}
Kiran~S. Kedlaya, Jonathan Pottharst, and Liang Xiao.
\newblock Cohomology of arithmetic families of {$(\varphi,\Gamma)$}-modules.
\newblock {\em J. Amer. Math. Soc.}, 27(4):1043--1115, 2014.

\bibitem{LLafforgue}
Laurent Lafforgue.
\newblock Chtoucas de {D}rinfeld et conjecture de {R}amanujan-{P}etersson.
\newblock {\em Ast\'erisque}, (243):ii+329, 1997.

\bibitem{VLafforgue1}
V.~Lafforgue.
\newblock Introduction to chtoucas for reductive groups and to the global
  {Langlands} parameterization.
\newblock arXiv:1404.6416v2, 2015.

\bibitem{VLafforgue2}
V.~Lafforgue.
\newblock Chtoucas pour les groupes r\'eductifs et param\'etrisation de
  {Langlands} globale.
\newblock {\em J. Amer. Math. Soc.}, 31(3):719--891, 2018.

\bibitem{Lau}
Eike Lau.
\newblock {\em On generalized {$\mathcal{D}$}-shtukas}.
\newblock PhD thesis, Universit\"at Bonn, 2004.
\newblock retrieved February 2017.

\bibitem{LiuHerr}
Ruochuan Liu.
\newblock Cohomology and duality for {$(\phi,\Gamma)$}-modules over the {R}obba
  ring.
\newblock {\em Int. Math. Res. Not. IMRN}, (3):Art. ID rnm150, 32, 2008.

\bibitem{PalZabradi}
Aprameyo Pal and Gergely Z\'{a}br\'{a}di.
\newblock Cohomology and overconvergence for representations of powers of
  {G}alois groups.
\newblock {\em J. Inst. Math. Jussieu}, 20(2):361--421, 2021.

\bibitem{Schneider}
Peter Schneider.
\newblock {\em Galois representations and {$(\varphi,\Gamma)$}-modules}, volume
  164 of {\em Cambridge Studies in Advanced Mathematics}.
\newblock Cambridge University Press, Cambridge, 2017.

\bibitem{SchneiderVignerasZabradi}
Peter Schneider, Marie-France Vigneras, and Gergely Zabradi.
\newblock From \'etale {$P_+$}-representations to {$G$}-equivariant sheaves on
  {$G/P$}.
\newblock In {\em Automorphic forms and {G}alois representations. {V}ol. 2},
  volume 415 of {\em London Math. Soc. Lecture Note Ser.}, pages 248--366.
  Cambridge Univ. Press, Cambridge, 2014.

\bibitem{ScholzeHodge}
Peter Scholze.
\newblock {$p$}-adic {H}odge theory for rigid-analytic varieties.
\newblock {\em Forum Math. Pi}, 1:e1, 77, 2013.

\bibitem{Scholze_TorsionCohomology}
Peter Scholze.
\newblock On torsion in the cohomology of locally symmetric varieties.
\newblock {\em Ann. of Math. (2)}, 182(3):945--1066, 2015.

\bibitem{ScholzeDiamonds}
Peter Scholze.
\newblock {\'E}tale cohomology of diamonds.
\newblock 2021.
\newblock arXiv:1709.07343v2.

\bibitem{Scholze}
Peter Scholze and Jared Weinstein.
\newblock {\em Berkeley lectures on {$p$}-adic geometry}.
\newblock Annals of Mathematics Studies. Princeton University Press, Princeton,
  2020.

\bibitem{SerreGalois}
Jean-Pierre Serre.
\newblock {\em Galois cohomology}.
\newblock Springer Monographs in Mathematics. Springer-Verlag, Berlin, english
  edition, 2002.
\newblock Translated from the French by Patrick Ion and revised by the author.

\bibitem{StacksProject}
The {Stacks Project Authors}.
\newblock \textit{Stacks Project}.
\newblock \url{https://stacks.math.columbia.edu}, 2018.

\bibitem{WeinsteinGQp}
Jared Weinstein.
\newblock {${\rm Gal}(\overline{\bf Q}_p/{\bf Q}_p)$} as a geometric
  fundamental group.
\newblock {\em Int. Math. Res. Not. IMRN}, (10):2964--2997, 2017.

\bibitem{Weinstein}
Jared Weinstein.
\newblock Adic spaces.
\newblock In {\em Perfectoid Spaces: Lectures from the 2017 Arizona Winter
  School}, pages 14--57. Amer. Math. Soc., 2019.

\bibitem{Zabradi}
Gergely Z\'abr\'adi.
\newblock Multivariable {$(\phi, {\Gamma})$}-modules and products of {G}alois
  groups.
\newblock {\em Math. Research Letters}, 25(2):687--721, 2018.

\bibitem{ZabradiFunctor}
Gergely Z\'abr\'adi.
\newblock Multivariable {$(\varphi,{\Gamma})$}-modules and smooth $o$-torsion
  representations.
\newblock {\em Selecta Math.}, 24(2):935--995, 2018.

\end{thebibliography}

\end{document}